\documentclass[11pt]{article}
\usepackage[T1]{fontenc}
\usepackage[utf8]{inputenc}
\usepackage{lmodern}
\usepackage{amsmath}
\usepackage{amssymb}
\usepackage{amsfonts}
\usepackage{amsthm}
\usepackage{mathtools}
\usepackage{mathrsfs} 
\usepackage{enumitem} 
\usepackage[a4paper,margin=30mm]{geometry}
\usepackage{microtype}

\usepackage[hidelinks]{hyperref}
\usepackage[nameinlink,noabbrev]{cleveref}

\numberwithin{equation}{section}

\newtheorem{theorem}{Theorem}[section]
\newtheorem{proposition}[theorem]{Proposition}
\newtheorem{lemma}[theorem]{Lemma}
\newtheorem{corollary}[theorem]{Corollary}

\theoremstyle{definition}
\newtheorem{definition}[theorem]{Definition}

\theoremstyle{remark}
\newtheorem{remark}[theorem]{Remark}

\newcommand{\ckgcasebeta}{\beta_{\mathrm{in}}}
\newcommand{\ckgnewtonbeta}{\beta_{\mathrm{Newt}}}
\newcommand{\iu}{\mathrm{i}}

\newcommand{\Id}{\mathrm{Id}}
\newcommand{\T}{\mathbb{T}}
\newcommand{\Z}{\mathbb{Z}}
\newcommand{\cP}{\mathcal{P}}
\newcommand{\cZ}{\mathcal{Z}}

\newcommand{\tR}{{\tt R}}
\newcommand{\cJ}{\mathcal{J}}
\newcommand{\cIR}{\mathcal{IR}}
\newcommand{\lip}{\mathrm{lip}}
\newcommand{\C}{\mathbb{C}}
\newcommand{\cF}{\mathcal{F}}
\newcommand{\cE}{\mathcal{E}}
\title{Global-in-Time Nonrelativistic Limit and KAM Tori for the Complex Klein-Gordon Equation on the Circle}

\author{A. Belloni\footnote{Dipartimento di Matematica, Politecnico di Milano, Piazza Leonardo Da Vinci 32, 20133, Milano, Italy email\,\,:\,\,\texttt{andrea.belloni@polimi.it}}}

\date{}
\begin{document}
\maketitle
\vspace{-2em}
\begin{abstract}
This paper establishes a global-in-time nonrelativistic limit for families of small-amplitude quasi-periodic solutions of the cubic complex Klein--Gordon equation on the circle. We first construct, by a KAM block scheme, small-amplitude two-frequency invariant tori and show that their linear stability depends on whether the two excited Fourier indices coincide or are distinct: coincident indices yield linearly stable tori, whereas distinct indices yield partially hyperbolic tori. Moreover, controlling the dependence of the KAM scheme on the speed of light parameter $c$, we construct quasi-periodic solutions whose frequencies, after a removal of the fast rest-energy phase, coincide exactly with those of quasi-periodic solutions of the limiting coupled NLS system, obtaining an approximation uniform for all times in the nonrelativistic regime. By this result, we further show that there exists a positive-measure set of frequencies such that, for each frequency, we can fix a quasi-periodic solution of the coupled NLS system and a suitable subsequence of Klein–Gordon solutions converging globally in time to the fixed NLS solution, as $c \to\infty$.
\end{abstract}

\tableofcontents

\section{Introduction}
\label{sec:introduction}

\subsection{The problem and the main result}

We consider the cubic complex Klein--Gordon equation on the one-dimensional torus $\mathbb T=\mathbb R/(2\pi\mathbb Z)$,
\begin{equation}\label{eq:intro-CKG}
\frac{1}{c^2}\partial_{tt}u-\partial_{xx}u+c^2u+{\rm a}_1|u|^2u=0,\qquad {\rm a}_1\in\{\pm1\},\qquad u:\mathbb R\times\mathbb T\to\mathbb C,
\end{equation}
where $c>1$ denotes the speed of light. Our main interest is to study the nonrelativistic regime $c\to\infty$. The nonrelativistic limit of nonlinear Klein--Gordon equations has been extensively studied. Classical results on convergence to a single nonlinear Schr\"odinger equation typically establish convergence in norms weaker than the regularity assumed for the initial data, and include for instance \cite{tsutsumi1984nonrelativistic,najman1990nonrelativistic,machihara2001nonrelativistic,machihara2002nonrelativistic}. In these settings, the leading-order dynamics can be described by a single envelope evolving on the slow time scale relative to the fast relativistic oscillations of frequency $c^2$. A more general description is obtained by retaining two independent envelopes associated with the positive- and negative-frequency temporal oscillations. This leads to a system for $ (\xi, \eta) : \mathbb{R} \times \mathbb{T} \to \mathbb{C}^2$ of coupled nonlinear Schr\"odinger equations 
\begin{equation}\label{eq:intro-NLS}
 \begin{cases}
 \displaystyle 2\iu\partial_t\xi+\partial_{xx}\xi
   -{\rm a}_1\bigl(|\xi|^2+2|\eta|^2\bigr)\xi=0,\\[0.35em]
 \displaystyle 2\iu\partial_t\eta+\partial_{xx}\eta -{\rm a}_1\bigl(|\eta|^2+2|\xi|^2\bigr)\eta=0.
 \end{cases}
\end{equation}
This coupled system was identified by Masmoudi and Nakanishi in \cite{masmoudi2002nonlinear}.  These classical approximation results with respect to nonlinear Schr\"odinger equations are typically valid on compact time intervals. Extending such approximations uniformly to arbitrarily long times is a substantially different problem.

A first global-in-time result was obtained by Bambusi, Belloni and Giuliani \cite{BBG} for the real cubic Klein--Gordon equation, in the setting of small-amplitude quasi-periodic solutions. The authors construct families of quasi-periodic solutions which, after the fast temporal oscillations have been removed, remain uniformly close for all times to suitable quasi-periodic solutions of the cubic NLS. Their approach is based on a KAM construction that simultaneously controls the Klein--Gordon and Schr\"odinger dynamics, allowing their frequencies to be matched exactly, up to a removal of the fast rest-energy phase. Subsequently, Cong, Li and Wu in \cite{cong2025existencefulldimensionalkamtori} extend this approach to almost-periodic solutions of the real nonlinear Klein--Gordon equation with a convolution potential.

Both works draw on KAM techniques previously developed for nonlinear wave, Klein-Gordon and Schr\"odinger equations (see, for instance, \cite{poschel1996kam,  kuksin1996invariant, poschel1996quasi, berti2013kam, congYuan2021, biasco2023small, corsi2024almost} ), adapted to investigate the nonrelativistic regime. The situation is different for the complex Klein--Gordon equation \eqref{eq:intro-CKG}, since to the best of our knowledge, no KAM theorem has previously been established for the autonomous cubic complex Klein--Gordon equation on the circle, even for a fixed value of $c$.

The purpose of this paper is therefore not only to establish a global-in-time nonrelativistic approximation, but also to develop a KAM construction for the complex Klein--Gordon equation itself. We show that the interaction between the two components leads to different normal geometries, including both elliptic and partially hyperbolic invariant tori. The coexistence of these two types of tori is known for related systems of coupled cubic Schr\"odinger equations, as shown by Gr\'ebert and Vila\c ca da Rocha \cite{GVdR}. However, their construction cannot be applied directly to the present Klein--Gordon problem, whose relativistic dispersion relation and resonant structure require a separate analysis. A further difficulty is to obtain the necessary KAM estimates uniformly as $c\to\infty$ and to compare the resulting transformations and frequencies with those of the limiting coupled Schr\"odinger system.

\bigskip

To formulate our results precisely, we first recall the diagonal variables associated with the complex Klein--Gordon equation. We perform a linear canonical transformation $(u,u_t/c^2)\mapsto(\xi,\eta)$ by diagonalizing the quadratic Hamiltonian associated with \eqref{eq:intro-CKG}, and identify $(\xi,\eta)$ with its Fourier sequence $\{(\xi_j,\eta_j)\}_{j\in\mathbb Z}$. The resulting nonlinear dynamics consist of two coupled complex evolution equations with the same positive linear frequencies 
\begin{equation}\label{eq:intro-linear-frequencies}
\lambda_j(c)=c\sqrt{c^2+j^2}=c^2+\frac{j^2}{2}+O(c^{-2}j^4),\qquad j\in\mathbb Z.
\end{equation}
The original complex field is recovered through
\begin{equation}\label{eq:intro-reconstruction}
u=A_c^{-1}(\xi+\bar\eta),\qquad A_c=(1-c^{-2}\partial_{xx})^{1/4}.
\end{equation}
For an invariant two-torus of the Klein--Gordon equation, we write its frequency vector in the form $c^2(1,1)+\omega$ and call $\omega\in\mathbb R^2$ its slow frequency vector. Our comparison concerns Klein--Gordon and NLS tori with exactly the same slow frequencies.  This exact matching is what makes the nonrelativistic approximation uniform over the whole time axis.

We study the problem on the exponentially weighted Sobolev space $\ell^{a,q}$ of pairs of Fourier sequences, for $a\ge0$ and $q\ge0$, with norm
\begin{equation}\label{norm:aq}
 \|(\xi,\eta)\|_{a,q}^2:=\sum_{j\in\mathbb Z}e^{2a|j|}\langle j\rangle^{2q}\bigl(|\xi_j|^2+|\eta_j|^2\bigr),\qquad \langle j\rangle:=\max\{1,|j|\},
\end{equation}
where $\xi_j,\eta_j$ are the Fourier coefficients of the two components. The same notation is used for the corresponding scalar norm. We consider invariant two-tori obtained by exciting one tangential Fourier mode in each component, with indices $m,n\in\mathbb Z$, in the following sense. Let 
\begin{equation}\label{Q:def}
(\rho_1,\rho_2) \in Q:=r_0^{4/3}\left[\tfrac{1}{2},\tfrac{3}{2}\right]^2,
\end{equation}
for some $0< r_0 \ll 1$.  The linear systems obtained by \eqref{eq:intro-CKG} and \eqref{eq:intro-NLS} without the nonlinear terms, expressed in diagonal coordinates, admit invariant two-dimensional tori supported on the Fourier modes $m$ of $\xi$ and $n$ of $\eta$. These tori are parametrized by the embeddings $\mathcal I_{m,n}:\mathbb T^2\times Q\longrightarrow\ell^{a,p}$ with
\begin{equation}\label{eq:intro-reference-torus}
 \mathcal I_{m,n}(\vartheta,\rho):=(\xi,\eta),\qquad \xi_j=\sqrt{\rho_1}e^{\iu\vartheta_1}\delta_{j,m},\quad \eta_j=\sqrt{\rho_2}e^{\iu\vartheta_2}\delta_{j,n},
\end{equation}
where $\delta_{j,k}$ denotes the Kronecker delta.

\smallskip

Our main result is the following.
\begin{theorem}\label{thm:intro-CKG}
Fix $m,n\in\mathbb Z$, $a\ge0$, $p>9/2$ and $0<\varsigma<1/18$. There exist constants $r_*,h_*,\kappa_0,C,\kappa>0$ such that, for $0<r_0<r_*$ and $0<h:=c^{-2}<\min\{h_*,\kappa_0r_0^{4/3}\}$, there exist a measurable set $\mathcal W_h\subset\mathbb R^2$ of positive measure and two bi-Lipschitz maps onto their images $\rho_h,\rho_0:\mathcal W_h\to Q$, where $Q$ is defined in \eqref{Q:def}, such that the following holds.

For each $\omega\in\mathcal W_h$, the complex Klein--Gordon equation \eqref{eq:intro-CKG} and the coupled NLS system \eqref{eq:intro-NLS} admit real-analytic embeddings of invariant two-tori
\[
 \Psi_h(\cdot,\rho_h(\omega)),\quad \Psi_0(\cdot,\rho_0(\omega)):\mathbb T^2\longrightarrow\ell^{a,p},
\]
whose internal frequency vectors are $c^2(1,1)+\omega$ and $\omega$, respectively. The corresponding actions satisfy
\begin{equation}\label{eq:intro-action-comparison}
 |\rho_h(\omega)-\rho_0(\omega)|\le Ch,
\end{equation}
and the selected parameter sets have relative measure
\begin{equation}\label{eq:intro-common-measure}
 \frac{|Q\setminus\rho_h(\mathcal W_h)|}{|Q|}+\frac{|Q\setminus\rho_0(\mathcal W_h)|}{|Q|}\le C\left(r_0^{2\varsigma/3}+\frac{h}{r_0^{4/3}}\right).
\end{equation}
Furthermore, these tori are small deformations of the linear two-mode tori in \eqref{eq:intro-reference-torus}: uniformly in $\omega\in\mathcal W_h$, they satisfy
\begin{equation}\label{eq:intro-torus-closeness}
\begin{aligned}
 \sup_{\vartheta\in\mathbb T^2}\biggl(&\frac{\|\Psi_h(\vartheta,\rho_h(\omega))-\mathcal I_{m,n}(\vartheta,\rho_h(\omega))\|_{a,p}}{\|\mathcal I_{m,n}(\vartheta,\rho_h(\omega))\|_{a,p}}\\
 &+\frac{\|\Psi_0(\vartheta,\rho_0(\omega))-\mathcal I_{m,n}(\vartheta,\rho_0(\omega))\|_{a,p}}{\|\mathcal I_{m,n}(\vartheta,\rho_0(\omega))\|_{a,p}}\biggr)\le Cr_0^{1/9-2\varsigma}.
\end{aligned}
\end{equation}
If $m=n$, both tori are linearly stable. If $m\ne n$, the normal linearized dynamics of both tori are partially hyperbolic, with exponential growth and decay rates bounded below by $\kappa r_0^{4/3}$. Finally, for any common initial angle $\vartheta_0\in\mathbb T^2$, the corresponding exact quasi-periodic solutions in diagonal coordinates are
\begin{align*}
 (\xi_h(t),\eta_h(t))&:=\Psi_h(\vartheta_0-(c^2(1,1)+\omega)t,\rho_h(\omega)),\\
 (\xi_0(t),\eta_0(t))&:=\Psi_0(\vartheta_0-\omega t,\rho_0(\omega)),
\end{align*}
and writing $\mathscr R_\phi(\xi,\eta):=(e^{-\iu\phi}\xi,e^{-\iu\phi}\eta)$, for every $0<s\le1$ we have, uniformly in $\omega$ and $\vartheta_0$,
\begin{equation}\label{eq:intro-uniform-comparison}
 \sup_{t\in\mathbb R}\bigl\|\mathscr R_{-c^2t}(\xi_h(t),\eta_h(t))-(\xi_0(t),\eta_0(t))\bigr\|_{a,p-4s}\le Ch^s r_0^{\frac23(1-s)-\frac32s}.
\end{equation}
For the physical Klein--Gordon solution $u_h=A_c^{-1}(\xi_h+\overline{\eta_h})$ (as in \eqref{eq:intro-reconstruction}), let
\begin{equation}\label{eq:intro-two-carrier}
u_{\mathrm{app},c}(t) :=e^{-\iu c^2t}\xi_0(t) +e^{\iu c^2t}\overline{\eta_0(t)}.
\end{equation}
Then
\begin{equation}\label{eq:intro-physical-comparison}
 \sup_{t\in\mathbb R}\|u_h(t)-u_{\mathrm{app},c}(t)\|_{a,p-4s}\le Ch^s r_0^{\frac23(1-s)-\frac32s},\qquad 0<s\le1.
\end{equation}
\end{theorem}

By linear stability we mean stability of the normal linearized dynamics of the invariant torus; the tangential dynamics are the quasi-periodic motion on the torus. Moreover, the selected action sets have asymptotically full relative measure in the joint limit $r_0\to0$, $h/r_0^{4/3}\to0$. For fixed admissible $r_0$ and $s>0$, the global-in-time errors \eqref{eq:intro-uniform-comparison} and \eqref{eq:intro-physical-comparison} vanish as $c\to\infty$, at the expense of a loss of $4s$ derivatives.

\smallskip

We observe that the common frequency set and the corresponding KG invariant tori in Theorem~\ref{thm:intro-CKG} are selected for each value of $c$, and change with it, whereas the NLS embedding is constructed independently of \(c\), and the maps \(\rho_0\) are restrictions of the same \(c\)-independent inverse frequency map. We can further strengthen this comparison by fixing a single NLS torus along a suitable subsequence, by proving the existence of a positive-measure set of common frequencies\footnote{After completing the proof of Theorem~\ref{cor:intro-CKG-subsequence} (see Section~\ref{proofcorollary}), we learned through correspondence with one of the authors of \cite{cong2025existencefulldimensionalkamtori} that they are developing independently a similar strategy, to be included in a further revision of \cite{cong2025existencefulldimensionalkamtori}, still in preparation. We jointly agreed to emphasize that the results were obtained independently.}.
\begin{theorem}[Non-relativistic limit of KAM solutions.]\label{cor:intro-CKG-subsequence} Under the assumptions of Theorem~\ref{thm:intro-CKG}, fix a sufficiently small $r_0$ and let $c_j\to+\infty$ be any sequence. Set $h_j=c_j^{-2}$ and, after discarding finitely many initial terms, assume that the hypotheses of Theorem~\ref{thm:intro-CKG} hold for every $j$. By setting
\[
 \mathcal W_*:=\limsup_{j\to\infty}\mathcal W_{h_j} =\bigcap_{N\ge1}\bigcup_{j\ge N}\mathcal W_{h_j},
\]
there exists a constant $L\ge1$, independent of the sequence and of $r_0$, such that, for $r_0$ sufficiently small,
\begin{equation}\label{eq:intro-fixed-NLS-measure}
 |\mathcal W_*|\ge L^{-2}r_0^{8/3}(1-Cr_0^{2\varsigma/3})>0.
\end{equation}
Then, for every $\omega\in\mathcal W_*$ and any fixed initial angle $\vartheta_0\in\mathbb T^2$, one can select a subsequence $c_{j_k}\to\infty$, corresponding quasi-periodic Klein--Gordon solutions in diagonal coordinates $(\xi_{h_{j_k}}(t),\eta_{h_{j_k}}(t))$ with the same slow frequency $\omega$, and a \emph{single} exact solution $(\xi_0(t),\eta_0(t))$ of the coupled NLS system \eqref{eq:intro-NLS}, independent of $c_{j_k}$, such that, for every $0<s\le1$,
\begin{equation}\label{eq:intro-fixed-NLS-limit}
 \lim_{k\to\infty}\sup_{t\in\mathbb R}\bigl\|\mathscr R_{-c_{j_k}^2t}(\xi_{h_{j_k}}(t),\eta_{h_{j_k}}(t))-(\xi_0(t),\eta_0(t))\bigr\|_{a,p-4s}=0.
\end{equation}
The corresponding physical Klein--Gordon solutions satisfy the analogous uniform-in-time convergence with the two-carrier approximation \eqref{eq:intro-two-carrier}, constructed using this same fixed NLS solution, as in \eqref{eq:intro-physical-comparison}.
\end{theorem}

\begin{remark}
The two results presented above concern two-dimensional invariant tori, obtained by exciting one tangential Fourier mode in each component. As we stated in Theorem \ref{thm:intro-CKG}, even in this setting the interaction between the components produces two qualitatively different normal geometries for the resulting tori. We believe that the construction developed in this paper already addresses most of the main difficulties which would arise also in the construction of invariant tori with an arbitrary number of tangential modes. However, that case would not be a straightforward generalization, as additional resonant interactions may lead to more complicated coupling structures. We do not pursue this question here, since our main purpose is the study of the non-relativistic limit. We expect that the mechanism underlying the global-in-time KG--NLS comparison, developed in this paper, remains applicable to higher-dimensional invariant tori, provided a suitable uniform KAM construction is available (see in particular Remark \ref{gen}).
\end{remark}

\subsection{Related literature}

The non-relativistic limit of nonlinear Klein--Gordon equations has a long history, beginning with the classical convergence results recalled above \cite{tsutsumi1984nonrelativistic,najman1990nonrelativistic,machihara2001nonrelativistic,machihara2002nonrelativistic, masmoudi2002nonlinear}. Another important result about the non-relativistic limit was obtained in the framework of scattering theory in \cite{nakanishi2002scattering}. Finally, we observe that non-relativistic limits have also been established for coupled relativistic systems, including for instance the Klein--Gordon--Maxwell equations (see \cite{bechouche2004nonrelativistic} and the references therein). More recently, considerable progress has been made in quantifying the dependence of non-relativistic approximation errors on time. For the real cubic Klein--Gordon equation on $\mathbb R^3$, Bao, Lu and Zhang \cite{baoLuZhang2024} obtained error estimates growing linearly in time, valid on time intervals increasing with the speed of light. Lei and Wu \cite{leiWu2023} established global-in-time error bounds in dimensions two and three; in particular, they obtained a uniform-in-time estimate when the approximating profile solves a \emph{Schr\"odinger--wave equation}, rather than the limiting Schr\"odinger equation usually studied. Higher-order expansions have subsequently been studied in \cite{shenWangZheng2026}. For complex Klein--Gordon equations, the recent work \cite{JiaYang2026} gives quantitative estimates for the approximation by a system of coupled NLS equations in dimensions two and three, with error constants growing linearly in time. Very recently Lu and Zheng \cite{luZheng2026} obtained a uniform-in-time approximation for a Klein--Gordon--Schr\"odinger system on $\mathbb R^d$, $d\ge4$, with respect to a \emph{Schr\"odinger--wave semi-limiting} model. This concerns a different coupled relativistic system from the classical one considered here. 

\medskip

In the periodic setting, Lin and Wu \cite{linWu2010} established the non-relativistic limit from a modulated defocusing nonlinear Klein--Gordon equation to the defocusing NLS on $\mathbb T^d$, on finite time intervals. Faou and Schratz \cite{faouSchratz2014} developed asymptotic expansions for nonlinear Klein--Gordon equations with periodic boundary conditions, including complex-valued solutions whose leading profiles satisfy a coupled Schr\"odinger system, and used these expansions to construct asymptotic-preserving numerical schemes. A different approach to the long-time dynamics is provided by Hamiltonian normal-form and KAM methods. In this direction, we refer to \cite{pasquali2019dynamics, pasquali2018almost}, while the KAM-based results \cite{BBG,cong2025existencefulldimensionalkamtori}, discussed above, yield uniform-in-time non-relativistic approximations along invariant tori. The present work follows this latter approach, with the additional difficulty of constructing the invariant tori for the complex Klein--Gordon equation while controlling their normal stability uniformly in the non-relativistic parameter.

\medskip

Our KAM construction builds on the classical results for Hamiltonian PDEs \cite{Kuk87, way90, kuksin1996invariant, poschel1996kam, poschel1996quasi, bourgain1998Hamiltonian} and on subsequent developments (see, for instance, \cite{chierchiaYou2000, bertiBiasco2011, berti2013kam, EGK, bertiKappelerMontalto2018, eliassonKuksin2010, procesiProcesi2015, bertiBolle2013, bambusiBertiMagistrelli2011, bambusiGraffi2001, bertiFeolaProcesiTerracina2024}). Of particular relevance is \cite{GVdR}, where both elliptic and hyperbolic invariant tori are constructed for a coupled NLS system. A natural further direction would be to extend this comparison to infinite-dimensional invariant tori and almost-periodic solutions for the system without convolutional potential. In this respect, recent constructions for nonlinear Schr\"odinger equations without external parameters \cite{ bernier2024infinite, bernierProcesi2026} provide an interesting perspective. Developing comparable almost-periodic KAM constructions for autonomous Klein--Gordon equations, even in the real case, and investigating their uniform non-relativistic limit remain open problems.

Another natural direction is to investigate whether analogous uniform-in-time results can be established for other singular limits, in particular those involving quasi-linear PDEs. A fundamental contribution in this direction is due to Franzoi and Montalto \cite{franzoi2024kam}, who established the uniform-in-time inviscid limit from quasi-periodic solutions of the forced two-dimensional Navier--Stokes equations to quasi-periodic solutions of the corresponding Euler equations. More broadly, KAM and Nash--Moser methods have been extensively developed for quasi-linear and fully nonlinear PDEs, including Airy and KdV type equations \cite{baldiBertiMontalto2014Airy, baldiBertiMontalto2016KdV, giuliani2017generalizedKdV, bertiKappelerMontalto2021KdV, corsiMontaltoProcesi2021Airy}, Schr\"odinger equations \cite{feolaProcesi2015, feolaProcesi2017autonomousNLS}, Kirchhoff equations \cite{montalto2017Kirchhoff, corsiMontalto2018Kirchhoff}, the Degasperis--Procesi equation \cite{feolaGiulianiProcesi2020DP}, incompressible Euler flows \cite{baldiMontalto2021Euler}, forced rotating fluids \cite{bianchiniFranzoiMontaltoTerracina2025} and water waves \cite{baldiBertiHausMontalto2018WW, bertiMontalto2020WW, bertiFranzoiMaspero2021WW, bertiFranzoiMaspero2024WW, feolaGiuliani2024WW, feolaMontaltoTerracina2025WW}.

\subsection{Organization of the paper and outline of the proof}

In Section~\ref{sec:CKG-Hamiltonian} we start by introducing the canonical Hamiltonian variables, the phase spaces and the limiting coupled-NLS system. After separating the fast linear oscillations, we decompose the Klein--Gordon Hamiltonian into the NLS part, terms that break the gauge symmetries of the single components, and remainders of order $h=c^{-2}$. As in \cite{BBG}, tracking this decomposition through the subsequent step in the algorithm is fundamental for the nonrelativistic comparison.

In Section~\ref{sec:first-BNF} we perform a partial Birkhoff normal form in parallel for KG and NLS, preserving the comparison between the two Hamiltonians. Besides the usual integrable terms, the resonant quartic terms include exchanges of Fourier modes between the two components. Section~\ref{sec:two-mode-normal-forms} analyzes how the choice of the tangential indices affects the normal dynamics. When $m=n$, the normal form is elliptic, whereas for $m\neq n$ the interaction between the two components produces a finite hyperbolic block, with eigenvalues uniformly separated from the imaginary axis at the action scale. In the latter case, a cubic term in the normal variables is too large for the KAM iteration. We eliminate it by a finite homological equation, using a hyperbolic resolvent bound of order the inverse of the action scale, and control the corresponding KG--NLS transformations.

In Sections~\ref{sec:block-KAM-setting}--\ref{sec:CKG-measure-estimates} we prove the block normal forms, solve the homological equation and establish the required Melnikov measure estimates. The finite hyperbolic block is retained in the normal form, while the remaining equations are treated by block-operator estimates.  We perform the small-divisor analysis using momentum conservation, gauge symmetries, and the asymptotic behavior of the normal frequencies. Since the Klein--Gordon frequencies are asymptotically linear, we follow the approach of \cite{poschel1996kam,BBG} and implement in Section~\ref{sec:CKG-KAM-theorem} a KAM iteration with a weighted preliminary step, whose normal-frequency correction is retained as a fixed reference. This step is fundamental to obtain frequency sets of asymptotically full measure.

Finally, Sections~\ref{sec:common-frequencies}--\ref{sec:final-comparison} select KG and NLS tori with exactly matching slow frequencies and prove the uniform-in-time approximation. Exact frequency matching prevents a secular phase drift, and reconstruction of the complex KG field yields the two opposite fast carriers in \eqref{eq:intro-two-carrier}. Identifying the NLS tori with the autonomous NLS KAM family then gives subsequential convergence to a single fixed NLS torus.

\medskip
\noindent\textbf{Notation.}
We write $X\lesssim Y$ if $X\le CY$ with a constant independent of the non-relativistic parameter $c$, and $X\asymp Y$ if both $X\lesssim Y$ and $Y\lesssim X$.  Additional parameter dependence of implicit constants is specified when needed.
Let $\mathcal A$ be a finite or countable index set. For sequences $v=(v_a)_{a\in\mathcal A}$ and $w=(w_a)_{a\in\mathcal A}$ we use the (at this level) formal complex-bilinear pairing $\langle v,w\rangle :={\sum}_{a\in\mathcal A}v_aw_a.$ For matrices, ${}^t$ denotes transposition without conjugation and ${}^*$ the Hermitian adjoint.
\bigskip

	\noindent
	{\bf Acknowledgements:} A. Belloni has been supported by INdAM/GNAMPA.

\vspace{0.2cm}

\noindent\textbf{AI declaration:} The author takes full responsibility for the mathematical ideas, originality, and content of the paper. During the preparation and revision of this manuscript, ChatGPT (OpenAI) was used to assist with language editing, literature discovery and consistency checks. All the outputs were independently verified by the author.

\vspace{0.2cm}

\noindent\textbf{Statements and Declarations.} Data sharing is not applicable to this article as no datasets were generated or analyzed during the current study.

\vspace{0.2cm}

\noindent\textbf{Competing Interests}: The author has no conflicts of interest to declare.
\section{The complex Klein-Gordon equation and its Hamiltonian structure}
\label{sec:CKG-Hamiltonian}

For clarity, we present the detailed construction for the defocusing choice ${\rm a}_1=1$ in \eqref{eq:intro-CKG}, namely
\begin{equation}\label{eq:CKG}
  \frac1{c^2}u_{tt}-u_{xx}+c^2u+|u|^2u=0, \qquad u=u(t,x)\in\mathbb C.
\end{equation}
By setting $v:=u_t/c^2$, the PDE has Hamiltonian given by
\begin{equation}\label{eq:HKG-uv}
  H(u,v) =\int_{\mathbb T} \left( \frac{c^2}{2}|v|^2+\frac12|u_x|^2 +\frac{c^2}{2}|u|^2+\frac14|u|^4 \right)\,dx ,
\end{equation}
with real symplectic form
\begin{equation}\label{eq:sympl-uv}
  \Omega_{u,v} =-\Re\int_{\mathbb T}du\wedge d\bar v\,dx.
\end{equation}
Throughout the paper we use the Hamiltonian convention
\begin{equation}\label{eq:Hamiltonian-convention}
 \Omega(X_F,Y)=-dF[Y],
\end{equation}
for every Hamiltonian $F$ and tangential vector $Y$, and we define the Poisson bracket as
\begin{equation*}
 \{F,G\}:=dF[X_G].
\end{equation*}

For the real symplectic form \eqref{eq:sympl-uv}, let $\nabla_u F,\nabla_v F$ denote the real $L^2$ gradients, namely
\[
 dF[\widehat u,\widehat v] = \Re\int_{\mathbb T} \left( \nabla_uF\,\overline{\widehat u} +\nabla_vF\,\overline{\widehat v} \right)\,dx, \qquad \forall \hat{u},\hat{v} \in L^2.
\]
Then
\begin{equation*}
 X_F(u,v) = \bigl(\nabla_vF,-\nabla_uF\bigr),
\end{equation*}
and therefore for any Hamiltonian $G$ it holds
\begin{equation*}
 \{F,G\} = \Re\int_{\mathbb T} \left( \nabla_uF\,\overline{\nabla_vG} - \nabla_vF\,\overline{\nabla_uG} \right)\,dx.
\end{equation*}
The equation has two symmetries which will be used throughout the paper. It is invariant under spatial translations,
\begin{equation*}
  u(x)\mapsto u(x+\varphi), \qquad v(x)\mapsto v(x+\varphi), \qquad \varphi\in\mathbb T,
\end{equation*}
and under the physical gauge action
\begin{equation}\label{eq:physical-gauge-u}
  u\mapsto e^{\mathrm i\vartheta}u, \qquad v\mapsto e^{\mathrm i\vartheta}v, \qquad \vartheta\in\mathbb T.
\end{equation}

Following the direct diagonalization in \cite[Section~2]{BBG}, we diagonalize the quadratic Hamiltonian by setting
\begin{equation*}
  \langle\nabla\rangle_c:=\sqrt{c^2-\partial_{xx}}, \qquad A_c:=\left(\frac{\langle\nabla\rangle_c}{c}\right)^{1/2}.
\end{equation*}
Thus $A_c$ is the Fourier multiplier with symbol $(1+j^2/c^2)^{1/4}$ and $A_c^\alpha$ is the operator with Fourier multiplier $(1+j^2/c^2)^{\alpha/4}$. We introduce the complex variables
\begin{equation}\label{eq:xi-eta-direct}
   \begin{aligned}
  \xi&:=\frac12\left(A_cu+\mathrm i A_c^{-1}v\right), \qquad \eta&:=\frac12\left(A_c\bar u+\mathrm i A_c^{-1}\bar v\right)
  \end{aligned}
\end{equation}
and they satisfy the following properties.
\begin{lemma}
For every $q\ge1$, $a \geq 0$ and $c>1$, the transformation \eqref{eq:xi-eta-direct} defines a bounded real-linear isomorphism
\begin{equation}\label{domain}
 \Phi_c: \ell^{a,q}\times \ell^{a,q-1} \longrightarrow \ell^{a,q-\frac12} \times \ell^{a,q-\frac12},
\end{equation}
with $\ell^{a,q}$ introduced by \eqref{norm:aq}. It is canonical and its inverse is given by
\begin{equation}\label{eq:xi-eta-inverse}
  u=A_c^{-1}(\xi+\bar\eta), \qquad v=-\mathrm i A_c(\xi-\bar\eta).
\end{equation}
In the variables $(\xi,\eta)$ the symplectic form becomes
\begin{equation}\label{eq:sympl-xieta}
  \Omega =\mathrm i\int_{\mathbb T} \left(d\bar\xi\wedge d\xi+d\bar\eta\wedge d\eta\right)\,dx,
\end{equation}
and the Hamiltonian \eqref{eq:HKG-uv} is written as
\begin{equation}\label{H:xi-eta}
  H(\xi,\bar\xi,\eta,\bar\eta)=\Lambda+P,
\end{equation}
with
\begin{align}
  \Lambda &:=\int_{\mathbb T} \left( \bar\xi\,c\langle\nabla\rangle_c\xi +\bar\eta\,c\langle\nabla\rangle_c\eta
  \right)\,dx,\label{eq:Lambda-field}\\
  P &:=\frac14\int_{\mathbb T} \left|A_c^{-1}(\xi+\bar\eta)\right|^4\,dx.\label{eq:P-field}
\end{align}
\end{lemma}

\begin{proof}
First, \eqref{domain} follows by
\[
\begin{aligned}
 A_c:\ell^{a,q}\to \ell^{a,q-\frac12}, \qquad A_c^{-1}:\ell^{a,q-1}\to \ell^{a,q-\frac12}, \\
 A_c^{-1}:\ell^{a,q-\frac12}\to \ell^{a,q}, \qquad A_c: \ell^{a,q-\frac12} \to \ell^{a,q-1}.
\end{aligned}
\]
The inverse formula \eqref{eq:xi-eta-inverse} follows directly by adding and subtracting $\xi$ to $\bar\eta$.  Since $A_c$ is real, positive and self-adjoint, substitution of \eqref{eq:xi-eta-inverse} into \eqref{eq:sympl-uv} gives \eqref{eq:sympl-xieta}.  Then, we explicitly compute the Hamiltonians.  By \eqref{eq:xi-eta-inverse},
\[
  c^2\|v\|_{L^2(\T)}^2 =c^2\|A_c(\xi-\bar\eta)\|_{L^2(\T)}^2, \qquad \|u_x\|_{L^2(\T)}^2+c^2\|u\|_{L^2(\T)}^2 =\|\langle\nabla\rangle_c A_c^{-1}(\xi+\bar\eta)\|_{L^2(\T)}^2.
\]
Since $A_c^2=\langle\nabla\rangle_c/c$, both quadratic terms carry the same operator $c\langle\nabla\rangle_c$.  Expanding the squares, the mixed terms cancel, and we obtain
\[
  \frac{c^2}{2}\|v\|_{L^2(\T)}^2 +\frac12\|\langle\nabla\rangle_c u\|_{L^2(\T)}^2 =\langle\xi,c\langle\nabla\rangle_c\xi\rangle_{L^2(\T)} +\langle\eta,c\langle\nabla\rangle_c\eta\rangle_{L^2(\T)}.
\]
The quartic part follows directly from the expression of $u$ in \eqref{eq:xi-eta-inverse}.
\end{proof}

On the complexified phase space, barred and unbarred variables are treated as independent.  With the convention \eqref{eq:Hamiltonian-convention}, the Hamiltonian vector field associated with \eqref{eq:sympl-xieta} is
\begin{equation}\label{eq:XF-xieta}
 X_F = \left( -\iu\,\partial_{\bar\xi}F,\, \iu\,\partial_\xi F,\, -\iu\,\partial_{\bar\eta}F,\, \iu\,\partial_\eta F \right),
\end{equation}
and
\[
 \iu\,\dot\xi=\partial_{\bar\xi}F, \qquad \iu\,\dot\eta=\partial_{\bar\eta}F, \qquad \iu\,\dot{\bar{\xi}}=-\partial_{\xi}F, \qquad \iu\,\dot{\bar{\eta}}=-\partial_{\eta}F.
\]
The corresponding Poisson bracket is
\begin{equation}\label{eq:Poisson-xieta}
 \{F,G\} = \iu\int_{\mathbb T} \left( \partial_{\bar\xi}F\,\partial_\xi G -\partial_\xi F\,\partial_{\bar\xi}G +\partial_{\bar\eta}F\,\partial_\eta G -\partial_\eta F\,\partial_{\bar\eta}G \right)\,dx.
\end{equation}
The Hamilton equations in the diagonal variables can be written explicitly. With the symplectic convention \eqref{eq:sympl-xieta} and the vector-field formula \eqref{eq:XF-xieta}, the system \eqref{eq:CKG} becomes
\begin{equation*}
 \begin{cases}
 \displaystyle \iu\partial_t\xi = c\langle\nabla\rangle_c\,\xi +\frac12 A_c^{-1} \left[ \left|A_c^{-1}(\xi+\bar\eta)\right|^2 A_c^{-1}(\xi+\bar\eta) \right],
 \\[1.1em]
 \displaystyle \iu\partial_t\eta = c\langle\nabla\rangle_c\,\eta +\frac12 A_c^{-1} \left[ \left|A_c^{-1}(\xi+\bar\eta)\right|^2 A_c^{-1}(\bar\xi+\eta) \right].
 \end{cases}
\end{equation*}

We remark that in the new coordinates the physical gauge action \eqref{eq:physical-gauge-u} is
\begin{equation*}
  \xi\mapsto e^{\mathrm i\vartheta}\xi, \qquad \eta\mapsto e^{-\mathrm i\vartheta}\eta, \qquad \forall \vartheta \in \T,
\end{equation*}
and then the two components have opposite physical charges, even though their linear Klein-Gordon frequencies are identical.

\subsection{Phase space and the non-relativistic approximation}
\label{sec:phase-NR}

We now introduce an adapted phase space. This section follows the corresponding construction in \cite{BBG} with the difference given by the presence of a finite species index $s\in\{\pm 1\}$. Indeed, in \cite{BBG}, since $u=\bar u$ and $v=\bar v$, it holds $\xi=\eta$. Since this extra index is finite, most of the results about the phase space introduced in \cite{BBG} remain unchanged, up to absolute constants.

Set $h:=c^{-2}$. We introduce
\begin{equation}\label{jei}
\jmath(\alpha):=j, \quad \mathfrak s(\alpha):=s, \quad \widetilde w_j:=\sqrt{1+h j^2}, \qquad \forall \alpha=(j,s) \in \mathcal{Z}:=\Z \times \{\pm 1\}.
\end{equation}
Given $a\ge0$, $p>1/2$ and $\beta\in[0,1]$, we define
\begin{equation*}
  \ell^{a,p,\beta}(\mathcal Z) :=\left\{z=(z_\alpha)_{\alpha\in\mathcal Z}: \|z\|_{a,p,\beta}<\infty\right\},
\end{equation*}
with
\begin{equation}\label{eq:norm-apb}
  \|z\|_{a,p,\beta}^2 :=\sum_{\alpha\in\mathcal Z} |z_\alpha|^2 e^{2a|\jmath(\alpha)|} \langle\jmath(\alpha)\rangle^{2p} \widetilde w_{\jmath(\alpha)}^{2\beta}, \qquad \langle j\rangle:=\max\{1,|j|\}.
\end{equation}
We use the Fourier convention
\begin{equation*}
  \xi(x)=\frac1{\sqrt{2\pi}}\sum_{j\in\mathbb Z}\xi_j e^{\mathrm i jx}, \qquad \bar{\xi}(x)=\frac1{\sqrt{2\pi}}\sum_{j\in\mathbb Z}\bar{\xi}_j e^{-\mathrm i jx},
\end{equation*}
and the same for $\eta$, recalling that on the complexified space, barred and unbarred variables are considered independent. Then we set
\begin{equation*}
  \zeta_{(j,+1)}:=\xi_j, \qquad \zeta_{(j,-1)}:=\eta_j.
\end{equation*}
The phase space is
\begin{equation}\label{Papb}
  \mathcal P^{a,p,\beta} :=\ell^{a,p,\beta}(\mathcal Z)\times \ell^{a,p,\beta}(\mathcal Z),
\end{equation}
and we are going to consider $(\zeta,\bar\zeta) \in  \mathcal P^{a,p,\beta}$ with norm given by
\[
\|(\zeta, \bar\zeta)\|_{a,p,\beta}^2:=\|\zeta\|_{a,p,\beta}^2+\|\bar\zeta\|_{a,p,\beta}^2.
\]
For $\tR>0$ we denote the open ball of radius $\tR$ as
\begin{equation*}
 B_{a,p,\beta}(\tR):= \left\{(\zeta,\bar\zeta)\in\mathcal P^{a,p,\beta}: \|\zeta\|_{a,p,\beta}+\|\bar\zeta\|_{a,p,\beta}<\tR\right\},
\end{equation*}
and $B_{a,p}(\tR):=B_{a,p,0}(\tR)$. We discuss now some symmetries. For $s\in\{\pm 1\}$ we set the partial masses as
\begin{equation}\label{localmasses}
 \mathcal N_s :=\sum_{j\in\mathbb Z} \zeta_{(j,s)}\bar\zeta_{(j,s)}.
\end{equation}
We then define the total mass, the physical charge and the momentum by
\begin{equation}\label{eq:NQM-zeta}
 \begin{aligned}
 \mathcal N &:=\mathcal N_{+1}+\mathcal N_{-1}
 =\sum_{\alpha\in\mathcal Z}\zeta_\alpha\bar\zeta_\alpha,\\
 \mathcal Q &:=\mathcal N_{+1}-\mathcal N_{-1} =\sum_{\alpha\in\mathcal Z}
 \mathfrak s(\alpha)\zeta_\alpha\bar\zeta_\alpha,\\
 \mathcal M &:=\sum_{\alpha\in\mathcal Z} \jmath(\alpha)\zeta_\alpha\bar\zeta_\alpha.
 \end{aligned}
\end{equation}
Thus a Hamiltonian $F$ is
\[
\begin{aligned}
\emph{translational invariant if} \qquad &\{F,\mathcal M\}=0,&\\
\emph{physically gauge invariant if} \qquad &\{F,\mathcal Q\}=0,&\\
\emph{locally gauge invariant if} \qquad &\{F,\mathcal N_s\}=0, \qquad \forall s \in \{\pm1\}.&
\end{aligned}
\]
Equivalently, local gauge invariance is
\[
 \{F,\mathcal N\}=\{F,\mathcal Q\}=0.
\]
 For $\sigma\in\{+1,-1\}$ we put $\zeta_\alpha^{+1}:=\zeta_\alpha$ and $\zeta_\alpha^{-1}:=\bar\zeta_\alpha$.  If $n \in \mathbb{N}$, $\boldsymbol\alpha=(\alpha_1,\dots,\alpha_n) \in \cZ^n$ and $\boldsymbol\sigma=(\sigma_1,\dots,\sigma_n) \in \{ \pm 1 \}^n$, we write a monomial as
\begin{equation}\label{mon-xi-eta}
 \zeta_{\boldsymbol\alpha}^{\boldsymbol\sigma} :=\prod_{\ell=1}^n \zeta_{\alpha_\ell}^{\sigma_\ell}.
\end{equation}
Then we define
\begin{align}
 \pi_{\mathcal M}(\boldsymbol\alpha,\boldsymbol\sigma) &:=\sum_{\ell=1}^n \sigma_\ell\jmath(\alpha_\ell),
  \label{pi:moment}\\
 \pi_{\mathcal Q}(\boldsymbol\alpha,\boldsymbol\sigma) &:=\sum_{\ell=1}^n \sigma_\ell\mathfrak s(\alpha_\ell),
 \label{GaugeG}\\
 \pi_s(\boldsymbol\alpha,\boldsymbol\sigma) &:=\sum_{\ell=1}^n \sigma_\ell \mathbf 1_{\{\mathfrak s(\alpha_\ell)=s\}}, \qquad \forall s \in \{\pm1\}.
 \label{GaugeL}
\end{align}
Then by a direct computation it holds
\begin{align*}
 \{\mathcal M, \zeta_{\boldsymbol\alpha}^{\boldsymbol\sigma}\} =\iu\, \pi_{\mathcal M}(\boldsymbol\alpha,\boldsymbol\sigma)
 \zeta_{\boldsymbol\alpha}^{\boldsymbol\sigma}, \quad \{\mathcal Q, \zeta_{\boldsymbol\alpha}^{\boldsymbol\sigma}\} =\iu\, \pi_{\mathcal Q}(\boldsymbol\alpha,\boldsymbol\sigma)
 \zeta_{\boldsymbol\alpha}^{\boldsymbol\sigma},\quad  \{\mathcal N_s, \zeta_{\boldsymbol\alpha}^{\boldsymbol\sigma}\} =\iu\, \pi_s(\boldsymbol\alpha,\boldsymbol\sigma) \zeta_{\boldsymbol\alpha}^{\boldsymbol\sigma}.
\end{align*}
Hence a monomial is
\[
\begin{aligned}
\emph{translation invariant} \, \, \,  \iff \, \, \,  &\pi_{\mathcal M}(\boldsymbol\alpha,\boldsymbol\sigma)=0,& \\
\emph{physically gauge invariant} \, \, \,  \iff \, \, \,  &\pi_{\mathcal Q}(\boldsymbol\alpha,\boldsymbol\sigma)=0,& \\
\emph{locally gauge invariant} \, \, \,  \iff \, \, \, &\pi_s(\boldsymbol\alpha,\boldsymbol\sigma)=0, \qquad \forall s \in \{\pm1\}.& \\
\end{aligned}
\]
For a physically gauge invariant monomial, the latter condition is equivalent to
\begin{equation}\label{eq:fast-gauge-condition}
\pi_*(\boldsymbol\sigma):= \sum_{\ell=1}^n\sigma_\ell=0.
\end{equation}
Indeed, the sum and the difference of the two identities in \eqref{GaugeL} are precisely \eqref{eq:fast-gauge-condition} and \eqref{GaugeG}. We now rewrite the Hamiltonian in Fourier variables.  From \eqref{eq:Lambda-field},
\begin{equation*}
  \Lambda =\sum_{\alpha\in\mathcal Z}\lambda_{\jmath(\alpha)} \zeta_\alpha\bar\zeta_\alpha, \qquad \lambda_j:=c\sqrt{c^2+j^2}=h^{-1}\sqrt{1+h j^2}.
\end{equation*}
The quartic Hamiltonian is translation invariant and physically gauge invariant.  If we use the convention \eqref{mon-xi-eta}, it can be written as
\begin{equation}\label{PKG}
  P=\sum_{(\boldsymbol\alpha,\boldsymbol\sigma)\in \mathcal Z^4\times\{\pm1\}^4} P_{\boldsymbol\alpha,\boldsymbol\sigma} \zeta_{\boldsymbol\alpha}^{\boldsymbol\sigma},
\end{equation}
where
\begin{equation*}
  P_{\boldsymbol\alpha,\boldsymbol\sigma} =\begin{cases} \displaystyle \frac{1}{48\pi} \prod_{r=1}^4\widetilde w_{\jmath(\alpha_r)}^{-1/2},
  &\text{if }\pi_{\mathcal M}(\boldsymbol\alpha,\boldsymbol\sigma)=0\text{ and } \pi_{\mathcal Q}(\boldsymbol\alpha,\boldsymbol\sigma)=0,\\[1ex]
  0,&\text{otherwise}.
  \end{cases}
\end{equation*}

We consider now the limit system
\begin{equation}\label{NLS:coupledcubic}
  \begin{cases}
  \mathrm i\xi_t+\frac12\xi_{xx}
  -\left(\frac12|\xi|^2+|\eta|^2\right)\xi=0,\\[0.4em]
  \mathrm i\eta_t+\frac12\eta_{xx} -\left(\frac12|\eta|^2+|\xi|^2\right)\eta=0.
  \end{cases}
\end{equation}
Its Hamiltonian is 
\begin{equation}\label{HNLS}
  H^{\rm NLS} = \Lambda^{\rm NLS}+P^{\rm NLS} = \frac12\int_{\mathbb T} \left(|\xi_x|^2+|\eta_x|^2\right)\,dx +\frac14\int_{\mathbb T} \left(|\xi|^4+|\eta|^4+4|\xi|^2|\eta|^2\right)\,dx.
\end{equation}

If we use the convention \eqref{mon-xi-eta}, it can be written as
\begin{equation*}
  P^{\rm NLS} =\sum_{\boldsymbol\alpha,\boldsymbol\sigma} P^{\rm NLS}_{\boldsymbol\alpha,\boldsymbol\sigma} \zeta_{\boldsymbol\alpha}^{\boldsymbol\sigma},
\end{equation*}
with
\begin{equation}\label{P:NLS}
  P^{\rm NLS}_{\boldsymbol\alpha,\boldsymbol\sigma} =\begin{cases} \displaystyle\frac1{48\pi}, &\text{if } \pi_{\mathcal M}(\boldsymbol\alpha,\boldsymbol\sigma)=0,\ \pi_{\mathcal Q}(\boldsymbol\alpha,\boldsymbol\sigma)=0, \, \pi_*(\boldsymbol\sigma)=0 \\ 0, &\text{otherwise},
  \end{cases}
\end{equation}
while
\begin{equation*}
  \Lambda^{\rm NLS} =\sum_{\alpha\in\mathcal Z}\lambda^{\rm NLS}_{\jmath(\alpha)} \zeta_\alpha\bar\zeta_\alpha, \qquad \lambda^{\rm NLS}_j:=j^2/2.
\end{equation*}

\medskip

We want to study the relation between the two Hamiltonians \eqref{HNLS} and \eqref{H:xi-eta}. Starting with the linear frequencies we set
\begin{equation*}
  \lambda_j=h^{-1}+\nu_j(h), \qquad \nu_j(h):=\frac{\sqrt{1+h j^2}-1}{h} =\frac{j^2}{1+\sqrt{1+h j^2}}.
\end{equation*}
If we denote
\begin{equation*}
  \Lambda^R :=\sum_{\alpha\in\mathcal Z} \left(\nu_{\jmath(\alpha)}(h)-\frac{\jmath(\alpha)^2}{2}\right) \zeta_\alpha\bar\zeta_\alpha,
\end{equation*}
we have the following.
\begin{remark}\label{expaHlin}
For every $a\ge0$ and $p>9/2$ and $\beta \in [0,1]$ we have 
\begin{equation}\label{eq:Lambda-decomp}
  \Lambda=h^{-1}\mathcal N+\Lambda^{\rm NLS}+\Lambda^R,
\end{equation}
where
\begin{equation*}
  \|X_{\Lambda^R}(\zeta,\bar\zeta)\|_{a,p-4} \lesssim h\|(\zeta,\bar\zeta)\|_{a,p}.
\end{equation*}
The estimate follows from the weighted norm \eqref{eq:norm-apb} and \eqref{eq:XF-xieta}, as in \cite{BBG}.
\end{remark}
Concerning the perturbations, we have the following lemma
\begin{lemma}
  \label{decompP}
 For the perturbations $P$ and $P^{\rm NLS}$ introduced in \eqref{eq:P-field} and \eqref{HNLS}, the following decomposition holds
  \begin{align}
\label{decompP-eq}
&P=P^{\rm NLS}+P^{\mathcal{NG}}+P^{R}, \\
\label{P:NG0}
    &P^{\mathcal{NG}} :=\int_{\mathbb T}\Bigg[ \frac12\left(|X|^2+|Y|^2\right) \left(XY+\bar X\bar Y\right) +\frac14\left(X^2Y^2+\bar X^2\bar Y^2\right) \Bigg]dx,
  \end{align}
  where $X=A_c^{-1}\xi$, $Y=A_c^{-1}\eta$, $\bar X=A_c^{-1}\bar \xi$, $\bar Y=A_c^{-1}\bar \eta$, while $P^R$ is defined by \eqref{decompP-eq}. Moreover, for all $\tR >0$ we have
  \begin{align} 
 \label{P:in}
&\sup_{{(\zeta,\bar \zeta)} \in B_{a,q,0}(\tR)}\left\|X_{P}{(\zeta,\bar \zeta)}
  \right\|_{a,q,1}\lesssim\tR^3,& \qquad  &\sup_{{(\zeta,\bar \zeta)} \in B_{a,p}(\tR)}\left\|X_{P^{R}}{(\zeta,\bar \zeta)}
  \right\|_{a,p-2}\lesssim h\tR^3,& \\
 \label{NG:P}
 &\sup_{{(\zeta,\bar \zeta)} \in B_{a,q}(\tR)}\left\|X_{P^{\mathcal{NG}}}{(\zeta,\bar \zeta)}
  \right\|_{a,q}\lesssim\tR^3,& \qquad &\sup_{{(\zeta,\bar \zeta)} \in B_{a, q}(\tR)}\left\|X_{P^{\rm NLS}}{(\zeta,\bar \zeta)} \right\|_{a,q}\lesssim\tR^3,&
  \end{align}
for both $q=p$ and $q=p-4$. Moreover, $P^{\rm NLS}$ is locally gauge invariant, whereas $P^{\mathcal{NG}}$ contains no locally gauge invariant monomial.
\end{lemma}

\begin{proof}
The proof follows \cite[Lemma~2.9]{BBG}, using \eqref{eq:XF-xieta} and the weighted norm \eqref{eq:norm-apb}. We only discuss the different expression of the Hamiltonians with respect to the real case. By the definition of $(X,Y,\bar X, \bar Y)$ in the statement,  \eqref{eq:P-field} reads as
\begin{equation}\label{eq:Pinf}
  P =\frac14\int_{\mathbb T}|X+\bar Y|^4\,dx
\end{equation}
and the projection onto monomials satisfying $\sum_r\sigma_r\neq 0$ gives exactly \eqref{P:NG0}. On the other hand, the projection onto monomials satisfying $\sum_r\sigma_r= 0$ gives
\begin{equation}\label{PG}
  P^{\mathcal{G}} :=\frac14\int_{\mathbb T} \left(|X|^4+|Y|^4+4|X|^2|Y|^2\right)\,dx.
\end{equation}
The proof of the estimates  in \eqref{P:in} and \eqref{NG:P} follows by Lemma \ref{sti.abs.2}. In particular, the second estimate in \eqref{P:in} follows by applying Lemma \ref{sti.abs.2} to $P^R=P^{\mathcal{G}}-P^{\rm NLS}$ since, by \eqref{HNLS} and \eqref{PG}, this difference has Fourier coefficients
\[
  \frac1{48\pi} \left( \prod_{r=1}^4\widetilde w_{j_r}^{-1/2}-1 \right).
\]
\end{proof}
Collecting \eqref{eq:Lambda-decomp} and \eqref{decompP-eq}, we have obtained
\begin{equation*}
  H=h^{-1}\mathcal N+H^{\rm NLS} +P^{\mathcal{NG}}+\Lambda^R+P^R.
\end{equation*}
This is the decomposition which will be used in the first Birkhoff step.  It is the exact two-component analogue of the decomposition in \cite{BBG}; the new feature is that the resonant quartic part of $H^{\rm NLS}$ is not a function of the actions alone.
\section{The unperturbed two-mode tori}
\subsection{First Birkhoff normal form for the complex Klein-Gordon Hamiltonian}
\label{sec:first-BNF}

Let $\mathcal J\subset\mathcal Z$ be a fixed finite set of internal modes. For $(\zeta,\bar\zeta) \in \cP^{a,p,\beta}$ we define 
\begin{equation}\label{pro}
\Pi_{\cJ}(\zeta,\bar\zeta):=(\zeta_{\alpha},\bar\zeta_{\alpha})_{\alpha \in \cJ}, \qquad \Pi_{\cJ}^\perp(\zeta,\bar\zeta):=(\zeta_{\alpha},\bar\zeta_{\alpha})_{\alpha \in \cZ \setminus  \cJ},
\end{equation}
and we denote
\begin{equation}\label{eq:Jfourier}
  J:=\{\jmath(\alpha):\alpha\in\mathcal J\}\subset\mathbb Z.
\end{equation}

For $(\boldsymbol\alpha,\boldsymbol\sigma)\in \mathcal Z^4\times\{\pm1\}^4$, we define
\begin{equation}\label{eq:OmegaKG}
  \Omega(\boldsymbol\alpha,\boldsymbol\sigma) :=\sum_{r=1}^4 \sigma_r\lambda_{\jmath(\alpha_r)}, \qquad  \Omega^{\rm NLS}(\boldsymbol\alpha,\boldsymbol\sigma) :=\sum_{r=1}^4 \sigma_r\lambda^{\rm NLS}_{\jmath(\alpha_r)}.
\end{equation}

We start by characterizing exactly the set of $4-$resonances.

\begin{definition}\label{resonant}
Recalling \eqref{jei}, we denote
\[
  \alpha_r=(j_r,s_r)\in\mathcal Z, \qquad j_r:=\jmath(\alpha_r), \qquad s_r:=\mathfrak{s}(\alpha_r), \qquad r=1,\ldots,4.
\]
Then we define
\[
  \mathcal{IR}\subset \mathcal Z^4\times\{\pm1\}^4
\]
as the set of indices
\[
  (\boldsymbol\alpha,\boldsymbol\sigma) = \bigl( (\alpha_1,\ldots,\alpha_4), (\sigma_1,\ldots,\sigma_4) \bigr)
\]
such that there exists a permutation $\pi\in\mathcal S_4$ for which
\begin{equation*}
  j_{\pi(1)}=j_{\pi(2)}, \qquad j_{\pi(3)}=j_{\pi(4)}, \qquad \sigma_{\pi(1)}=-\sigma_{\pi(2)}, \qquad \sigma_{\pi(3)}=-\sigma_{\pi(4)},
\end{equation*}
Recalling \eqref{pi:moment}, we also set
\begin{equation}\label{LJ}
  \mathcal L_{\mathcal J} := \left\{ (\boldsymbol\alpha,\boldsymbol\sigma) \in\mathcal Z^4\times\{\pm1\}^4 \mid \pi_{\mathcal M}=0, \ \exists r\in\{1,2,3,4\} \ \text{s.t.} \, \, \, \alpha_r\in\mathcal J \right\}.
\end{equation}

Finally, for any quartic Hamiltonian on $\mathcal P^{a,q,\beta}$,
\[
  F = \sum_{(\boldsymbol\alpha,\boldsymbol\sigma) \in\mathcal Z^4\times\{\pm1\}^4} F_{\boldsymbol\alpha,\boldsymbol\sigma} \zeta_{\boldsymbol\alpha}^{\boldsymbol\sigma},
\]
we set
\[
  \Pi_{\mathcal{IR}}F := \sum_{(\boldsymbol\alpha,\boldsymbol\sigma)\in\mathcal{IR}} F_{\boldsymbol\alpha,\boldsymbol\sigma} \zeta_{\boldsymbol\alpha}^{\boldsymbol\sigma}, \qquad \Pi_{\mathcal{IR}\cap\mathcal L_{\mathcal J}}F := \sum_{(\boldsymbol\alpha,\boldsymbol\sigma) \in\mathcal{IR}\cap\mathcal L_{\mathcal J}} F_{\boldsymbol\alpha,\boldsymbol\sigma} \zeta_{\boldsymbol\alpha}^{\boldsymbol\sigma}.
\]
\end{definition}
The set $\cIR$ is precisely the set of momentum-preserving indices for which \eqref{eq:OmegaKG} vanishes as in \cite{berti2013kam}, after forgetting the species index. We discuss this generalization in Lemma~\ref{APP:lem:quartic-resonances}.

\medskip

Following the classical approach in \cite{kuksin1996invariant}, \cite{poschel1996quasi}, \cite{berti2013kam} and \cite{BBG}, the following initial normal form reduction allows us to isolate the amplitudes as internal parameters for the subsequent application of KAM theory.

\begin{theorem}[First partial Birkhoff normal form]\label{Teorema0}
Assume $p>9/2$ and fix a finite internal set $\mathcal J\subset\mathcal Z$. There exist $h_*,\tR_*,C>0$, depending only on $a,p$ and the fixed internal set such that, for every $0<h<h_*$ and $0<\tR<\tR_*$, there exist real analytic canonical maps
\begin{align}
  \mathcal T_0(h;\cdot)&:B_{a,q,\beta}(\tR)\longrightarrow
  B_{a,q,\beta}(2\tR),\label{eq:T0-map}\\
  \mathcal T_0^{\rm NLS}&:B_{a,q}(\tR)\longrightarrow B_{a,q}(2\tR),
  \label{eq:T0NLS-map}
\end{align}
for $q=p,p-4$ and $\beta\in[0,1]$, such that
\begin{align}
  H\circ\mathcal T_0 &=\Lambda+\Lambda_++\widehat P+P_{\rm BNF},
  \label{NFS}\\
  H^{\rm NLS}\circ\mathcal T_0^{\rm NLS} &=\Lambda^{\rm NLS}+\Lambda_+^{\rm NLS} +\widehat P^{\rm NLS}+P_{\rm BNF}^{\rm NLS},
  \label{NFS:NLS}
\end{align}
with all the Hamiltonians translation invariant and physically gauge invariant, while the NLS Hamiltonians are also locally gauge invariant.

The quartic normal forms are
\begin{align}
  \Lambda_+ &:=  \Pi_{\cIR\cap \mathcal{L}_{\mathcal{J}}}P= \sum_{(\boldsymbol\alpha,\boldsymbol\sigma) \in\mathcal L_{\mathcal J}\cap\mathcal{IR}} P_{\boldsymbol\alpha,\boldsymbol\sigma} \zeta_{\boldsymbol\alpha}^{\boldsymbol\sigma},
  \label{res.nij}\\
  \Lambda_+^{\rm NLS} &:= \Pi_{\cIR\cap \mathcal{L}_{\mathcal{J}}} P^{\rm NLS}= \sum_{(\boldsymbol\alpha,\boldsymbol\sigma) \in\mathcal L_{\mathcal J}\cap\mathcal{IR}} P^{\rm NLS}_{\boldsymbol\alpha,\boldsymbol\sigma} \zeta_{\boldsymbol\alpha}^{\boldsymbol\sigma},
  \label{L+NLS}
\end{align}
and satisfy
\begin{align}
  \sup_{B_{a,q}(\tR)}\|X_{\Lambda_+}\|_{a,q,1} +\sup_{B_{a,q}(\tR)}\|X_{\Lambda_+^{\rm NLS}}\|_{a,q}
  &\le C\tR^3,\label{res.0.2}\\
  \sup_{B_{a,p}(\tR)} \|X_{\Lambda_+-\Lambda_+^{\rm NLS}}\|_{a,p-2} &\le Ch\tR^3.\label{stimaLR}
\end{align}
About the change of coordinates we have
\begin{align}
  &\sup_{B_{a,q,0}(\tR)} \|\mathcal T_0-\mathrm{Id}\|_{a,q,1}
  \le C \tR^3, \qquad \sup_{B_{a,q}(\tR)} \|\mathcal T_0^{\rm NLS}-\mathrm{Id}\|_{a,q}
  \le C \tR^3,& \label{sti.tra.1}\\
  &\sup_{B_{a,p}(\tR)} \|\mathcal T_0-\mathcal T_0^{\rm NLS}\|_{a,p-4} \le C h\tR^3.&\label{res.tr.0}
\end{align}
The perturbations satisfy
\begin{align}
  \sup_{B_{a,q,0}(\tR)}\|X_{\widehat P}\|_{a,q,1} &\le C\tR^3,& \sup_{B_{a,q,0}(\tR)}\|X_{P_{\rm BNF}}\|_{a,q,1}
  &\le C\tR^5,\label{res.0.0}\\
  \sup_{B_{a,q}(\tR)}\|X_{\widehat P^{\rm NLS}}\|_{a,q} &\le C\tR^3,& \sup_{B_{a,q}(\tR)}\|X_{P_{\rm BNF}^{\rm NLS}}\|_{a,q} &\le C\tR^5,\label{res.NLS}
\end{align}
where the quartic remainders contain external variables, i.e. 
\begin{equation}\label{eq:Phat-external}
  \widehat P, \widehat P^{\rm NLS}= \mathcal{O}(\|\Pi_{\cJ}^\perp(\zeta,\bar\zeta)\|_{a,q}^4),
\end{equation}
with $\Pi_{\cJ}$ defined in \eqref{pro}. Finally, by setting
\begin{equation}\label{comp:PP}
  \widehat P^R:=\widehat P-\widehat P^{\rm NLS} -\widehat P^{\mathcal{NG}}, \qquad P_{\rm BNF}^R:=P_{\rm BNF}-P_{\rm BNF}^{\rm NLS}-P_{\rm BNF}^{\mathcal{NG}},
\end{equation}
(where the perturbations with the superscript $\mathcal{NG}$ are indexed on $\sum_r\sigma_r\ne0$) it holds
\begin{align}
 &\sup_{(z,\bar z)\in B_{a,q}(\tR)}\left\|X_{\widehat{P}^{\mathcal{NG}}}(z,\bar z)\right\|_{a,q}\leq C\tR^3& \qquad &\sup_{B_{a,p}(\tR)}\|X_{\widehat P^R}\|_{a,p-4}
  \le Ch\tR^3,&\label{eq:PhatR-est}\\
&\sup_{(z,\bar z)\in B_{a,q}(\tR)}\left\|X_{P_{\rm BNF}^{\mathcal{NG}}}(z,\bar z)\right\|_{a,q}\leq C \tR^5,& \qquad &\sup_{B_{a,p}(\tR)}\|X_{P_{\rm BNF}^R}\|_{a,p-4} \le Ch\tR^5.& \label{eq:PBNFR-est}
\end{align}
\end{theorem}

\begin{proof}
We follow the analytic argument in \cite[Section~3]{BBG}, presenting the general ideas and emphasizing the modifications required by the two-species structure. Recalling \eqref{jei}, for any $\alpha \in \cZ$ the linear frequency for both the equations depends only on its Fourier component $\jmath(\alpha)$. Hence, by the Poisson formula \eqref{eq:Poisson-xieta} and \eqref{eq:OmegaKG}, for every quartic monomial,
\[
 \{\Lambda,\zeta_{\boldsymbol\alpha}^{\boldsymbol\sigma}\} = \iu\Omega(\boldsymbol\alpha,\boldsymbol\sigma) \zeta_{\boldsymbol\alpha}^{\boldsymbol\sigma}, \qquad  \{\Lambda^{\rm NLS},\zeta_{\boldsymbol\alpha}^{\boldsymbol\sigma}\} = \iu\Omega^{\rm NLS}(\boldsymbol\alpha,\boldsymbol\sigma) \zeta_{\boldsymbol\alpha}^{\boldsymbol\sigma}.
\]
Thus the homological divisors are exactly the same as in the one-species argument of \cite{BBG} and \cite{berti2013kam}, with Lemma~\ref{lem:small-divisors-CKG}  allowing us to define
\[
 G:= -\sum_{(\boldsymbol\alpha,\boldsymbol\sigma) \in\mathcal L_{\mathcal J}\setminus\mathcal{IR}} \frac{P_{\boldsymbol\alpha,\boldsymbol\sigma}} {\iu\Omega(\boldsymbol\alpha,\boldsymbol\sigma)} \zeta_{\boldsymbol\alpha}^{\boldsymbol\sigma}, \qquad  G^{\rm NLS}:= -\sum_{\substack{(\boldsymbol\alpha,\boldsymbol\sigma)
 \in\mathcal L_{\mathcal J}\setminus\mathcal{IR}\\ \sum_r \sigma_r=0}} \frac{P^{\rm NLS}_{\boldsymbol\alpha,\boldsymbol\sigma}} {\iu\Omega^{\rm NLS}(\boldsymbol\alpha,\boldsymbol\sigma)} \zeta_{\boldsymbol\alpha}^{\boldsymbol\sigma},
\]
The algebraic difference with respect to the real case of \cite{BBG} is the resonant projection: by Lemma~ \ref{APP:lem:quartic-resonances}, the resonant quartic monomials are precisely those retained in \eqref{res.nij}--\eqref{L+NLS}, containing also the exchange monomials \eqref{APP:eq:exchange-monomial}. 

It holds 
\begin{equation}\label{om}
 \{\Lambda,G\}+P=\Lambda_++\widehat P, \qquad \{\Lambda^{\rm NLS},G^{\rm NLS}\}+P^{\rm NLS} =\Lambda_+^{\rm NLS}+\widehat P^{\rm NLS}.
\end{equation}
By the definition of $\mathcal L_{\mathcal J}$ in \eqref{LJ}, every non-resonant quartic monomial containing an internal variable is removed, whereas the resonant ones are retained in $\Lambda_+$ and $\Lambda_+^{\rm NLS}$. Consequently, the remaining quartic terms contain external variables only, which proves \eqref{eq:Phat-external}. Moreover, by the resonant projections \eqref{res.nij}--\eqref{L+NLS}, the bound
\[
|\widetilde w_j^{-1/2}-1|\lesssim h\langle j\rangle^2,
\]
and Lemma~\ref{sti.abs.2} we obtain \eqref{stimaLR}.

\medskip

Firstly, we observe that Lemma~\ref{sti.abs.2} yields, for $q=p,p-4$,
\begin{equation}\label{GNLS}
 \sup_{B_{a,q,0}(\tR)} \|X_G\|_{a,q,1}\lesssim \tR^3, \qquad \sup_{B_{a,q}(\tR)} \|X_{G^{\rm NLS}}\|_{a,q}\lesssim \tR^3.
\end{equation}
The same argument applied directly to the resonant projections gives the bounds in \eqref{res.0.2}.

Moreover, the decomposition of the coefficients into the non-locally-gauge-invariant part, the multiplier remainder, and the difference between the KG and NLS divisors is the same as in \cite{BBG}. Using Lemmas~\ref{decompP}, \ref{lem:small-divisors-CKG} and Remark \ref{expaHlin}, we obtain
\begin{equation}\label{GR:stim}
 \sup_{B_{a,p}(\tR)} \|X_{G-G^{\rm NLS}}\|_{a,p-4} \lesssim h\tR^3.
\end{equation}
We now set $\mathcal T_0:=\Phi_G^1,$ $\mathcal T_0^{\rm NLS}:=\Phi_{G^{\rm NLS}}^1$, which are canonical since they are Hamiltonian flows with respect to \eqref{eq:sympl-xieta}. The standard flow estimates used in \cite[Section~3]{BBG}, together with the bounds above, show that for $\tR$ sufficiently small they are real analytic on the domains stated in \eqref{eq:T0-map}-\eqref{eq:T0NLS-map}, and give the estimates \eqref{sti.tra.1}-\eqref{res.tr.0}.

Using the homological equation \eqref{om} and Taylor's formula with integral remainder, we obtain
\begin{equation}\label{h_0}
H\circ\phi_G^1 =\Lambda+\Lambda_++\widehat P+P_{\rm BNF},
\end{equation}
which gives  \eqref{NFS}, with
\begin{equation}\label{P0:int}
\begin{aligned}
P_{\rm BNF} &=\int_0^1\{P,G\}\circ\phi_G^s\,ds +\int_0^1(1-s)\{\{\Lambda,G\},G\}
     \circ\phi_G^s\,ds\\
&=\int_0^1s\{P,G\}\circ\phi_G^s\,ds +\int_0^1(1-s)\{\Lambda_++\widehat P,G\} \circ\phi_G^s\,ds.
\end{aligned}
\end{equation}
The corresponding estimate in \eqref{res.0.0} follows from \eqref{P:in}, \eqref{GNLS}, \eqref{res.0.2}, and the flow bounds established above, by the standard Lie-transform estimates. In particular, these estimates use the identities
\begin{equation}\label{stime}
\begin{aligned}
X_{\{P,G\}\circ\phi_G^s} &=\bigl( d\phi_G^{-s}[X_P,X_G]
  \bigr)\circ\phi_G^s,\\
[X_P,X_G] &=dX_P[X_G]-dX_G[X_P],
\end{aligned}
\end{equation}
together with Cauchy estimates for the vector fields and the differentials of the flows. The detailed analytic argument is the one given in \cite[Section~3]{BBG}. The same argument applied to $G^{\rm NLS}$ gives \eqref{NFS:NLS} and \eqref{res.NLS}.

We next establish the comparison estimates. Starting with the first decomposition in \eqref{comp:PP}, recalling \eqref{decompP-eq} we observe that $\widehat{P}^R$ and $\widehat{P}^{\mathcal{NG}}$ are the projections on purely external quartic monomials of $P^R$ and $P^{\mathcal{NG}}$ respectively. Then \eqref{eq:PhatR-est} follow by \eqref{P:in} and \eqref{NG:P}.

For the sixth-order remainders, we use the same decomposition of the Lie expansions as in the proof of the normal-form theorem in \cite[Section~3]{BBG}, by defining
\begin{equation}\label{P:0NG}
P_{\rm BNF}^{\mathcal{NG}} :=\int_0^1s  \{P^{\mathcal{NG}},G^{\rm NLS}\} \circ\phi_{G^{\rm NLS}}^s\,ds+\int_0^1(1-s)\{\widehat P^{\mathcal{NG}},G^{\rm NLS}\}\circ\phi_{G^{\rm NLS}}^s\,ds.
\end{equation}
This choice is motivated by the second identity in \eqref{P0:int}. Since $G^{\rm NLS}$ is locally gauge invariant, its flow and the corresponding Poisson brackets preserve the local gauge degree. Consequently, $P_{\rm BNF}^{\mathcal{NG}}$ is supported only on monomials with $\sum_r\sigma_r\ne0$. Its vector-field bound in \eqref{eq:PBNFR-est} follows from \eqref{NG:P}, \eqref{GNLS}, \eqref{sti.tra.1}, and the same Lie-transform estimates used for $P_{\rm BNF}$. Then, we set
\begin{equation}\label{P0R}
P_{\rm BNF}^R :=P_{\rm BNF}  -P_{\rm BNF}^{\mathcal{NG}}  -P_{\rm BNF}^{\rm NLS}.
\end{equation}
To estimate its vector field, we compare the corresponding terms in the integral representations of the KG and NLS remainders. Writing $G^R:=G-G^{\rm NLS}$, a representative difference admits the decomposition
\begin{align*}
\{P,G\}\circ\phi_G^s -\{P^{\rm NLS}+P^{\mathcal{NG}},G^{\rm NLS}\}\circ\phi_{G^{\rm NLS}}^s&= \{P^R,G\}\circ\phi_G^s\\
&+\{P^{\rm NLS}+P^{\mathcal{NG}},G^R\}\circ\phi_G^s\\
&+\{P^{\rm NLS}+P^{\mathcal{NG}},G^{\rm NLS}\}\circ\phi_G^s\\
&-\{P^{\rm NLS}+P^{\mathcal{NG}},G^{\rm NLS}\}\circ\phi_{G^{\rm NLS}}^s.
\end{align*}
The first two contributions are controlled by the bounds for $P^R$ and $G^R$, respectively, while the last two are handled by the comparison of the KG and NLS flows. The other terms in $P_{\rm BNF}^R$, related to the contribution of the second integral \eqref{P0:int}, are decomposed analogously, with differences involving $\Lambda_+-\Lambda_+^{\rm NLS}$, $\widehat P^R$, $G^R$, and the corresponding KG and NLS flows. Thus every term contains either a remainder satisfying a $O(h)$ bound in the weak topology or a difference of flows controlled by \eqref{res.tr.0}.

More precisely, the required estimates follow from \eqref{res.tr.0}, \eqref{stimaLR}, \eqref{eq:PhatR-est}, and \eqref{GR:stim}, together with the commutator identity \eqref{stime} and Cauchy estimates in the strong and weak topologies $q=p,p-4$. These are precisely the analytic estimates established in the proof of the normal-form theorem in \cite{BBG} for the corresponding Lie-expansion differences. The same argument applies here because the two-species quartic vector fields satisfy the required bounds by Lemma~\ref{sti.abs.2}, while their homological generators and flows satisfy the estimates recalled above. We thus obtain
\begin{equation*}
\sup_{B_{a,p}(\tR)} \|X_{P_{\rm BNF}^R}\|_{a,p-4} \le Ch\tR^5,
\end{equation*}
which completes \eqref{eq:PBNFR-est}.

\smallskip

Finally, all the generators preserve momentum and physical charge. Indeed, the homological equations divide the coefficients of the quartic monomials by scalar divisors without changing their supports. Moreover, $G^{\rm NLS}$ is locally gauge invariant. Since Poisson brackets and Hamiltonian flows preserve these symmetries, all the Hamiltonians in \eqref{NFS}--\eqref{NFS:NLS} are translation invariant and physically gauge invariant, and the NLS Hamiltonians are locally gauge invariant. This concludes the proof.

\end{proof}

\subsection{The normal form around the two-mode tori}
\label{sec:two-mode-normal-forms}

We now compute explicitly the resonant Hamiltonian entering the KAM analysis, which is the main difference from the real equation considered in \cite{BBG}. For later use we set
\begin{equation*}
  \mathsf P_j:=|\xi_j|^2+|\eta_j|^2 \qquad \mathsf C_j:=\xi_j\bar\eta_j, \qquad \overline{ \mathsf C_j}:=\bar\xi_j\eta_j \qquad j\in\mathbb Z.
\end{equation*}
where on the complexified phase space we use the algebraic notation $|a|^2:=a\bar a$.  Whenever $\Re M$ is used, it denotes the real Hamiltonian obtained from a complex monomial $M$ by adding its complex-conjugate monomial (notationally speaking) and dividing by two; on the real phase space this agrees with the usual notation.

\begin{lemma}

Let $P$ be as in \eqref{PKG} and $P^{\rm NLS}$ as in \eqref{P:NLS}. Then
\begin{equation}\label{ckg:Z4-closed}
\begin{aligned}
   \Pi_{\cIR}P=& \sum_{j\in\mathbb Z} \frac{\mathsf P_j^2+2 |\mathsf{C}_j|^2} {8\pi\widetilde w_j^2} +\sum_{j<k} \frac{\mathsf P_j\mathsf P_k +2\Re(\mathsf C_j\overline{\mathsf C_k})}
       {2\pi\widetilde w_j\widetilde w_k},& \\
\Pi_{\cIR}P^{\rm NLS}=&\sum_{j\in\mathbb Z} \frac{\mathsf P_j^2+2|\mathsf{C}_j|^2} {8\pi} +\sum_{j<k} \frac{\mathsf P_j\mathsf P_k +2\Re(\mathsf C_j\overline{\mathsf C_k})} {2\pi},&
\end{aligned}
\end{equation}
with $\Pi_{\cIR}$ defined as in Definition \ref{resonant}.
\end{lemma}

\begin{proof}
By Lemma  \ref{APP:lem:quartic-resonances}, every resonant monomial is obtained by pairing two equal Fourier indices with opposite Hamiltonian signs.  If all four indices coincide with $j$, the non-symmetrized expansion of $\frac14\int|\xi+\bar\eta|^4dx$ gives
\[
 \frac1{8\pi\widetilde w_j^2} \bigl(|\xi_j|^4+|\eta_j|^4+4|\xi_j|^2|\eta_j|^2\bigr).
\]
Since $\mathsf P_j^2+2|\mathsf C_j|^2 =|\xi_j|^4+|\eta_j|^4+4|\xi_j|^2|\eta_j|^2$, this is the first sum in \eqref{ckg:Z4-closed}.  If $j\neq k$, the action-type pairings give $\mathsf P_j\mathsf P_k/(2\pi\widetilde w_j\widetilde w_k)$, while the two exchange pairings give $\Re(\mathsf C_j\overline{\mathsf C_k})/ (\pi\widetilde w_j\widetilde w_k)$.  Summing over the unordered pair $\{j,k\}$ proves the formula.
\end{proof}

We separate the size of the actions from the radius of the complex neighbourhood. Consider $\tR$ in Theorem \ref{Teorema0} and $\varrho:=\tR^2$.  Following the scaling introduced in \cite{BBG} we set
\begin{equation}\label{ckg:varrho-r0}
 r_0:=\varrho^{3/4}, \qquad \Xi_0=\Xi_0(r_0):= \left[\frac12 r_0^{4/3},\frac32 r_0^{4/3} \right]^2.
\end{equation}

\subsubsection{Common KAM notation: parameter box, domains and jets}
\label{sec:common-KAM-notation}

Before distinguishing the two possible choices of tangential modes, we introduce notation common to both configurations.  The parameter $\rho=(\rho_1,\rho_2)\in\Xi_0$ parametrizes the two tangential actions. We set
\[
\cJ:=\{(m,+),(n,-)\}\subseteq\cZ,
\]
for some $m,n\in\mathbb Z$ (we are going to consider both $m=n$ and $m\ne n$). We consider a $\cF\subset\cZ$ with $\cF\cap\cJ=\varnothing$, which indexes the hyperbolic/elliptic directions, and we set
\begin{equation*}
\cE:=\cZ\setminus(\cF\cup\cJ),
\end{equation*}
which is going to be the indexes of the elliptic directions. 

We introduce then the KAM phase space
\begin{equation}\label{ckg:KAM-variables}
 \cP_{\textnormal{KAM}}^{a,p,\beta} :=\C^2/(2\pi\mathbb Z)^2\times\C^2 \times\cP_{\cE}^{a,p,\beta}\times(\C^{\cF})^2, \qquad (\theta,y,z,\bar z,W_{\cF})\in\cP_{\textnormal{KAM}}^{a,p,\beta}.
\end{equation}
Here $(z,\bar z)\in\cP_{\cE}^{a,p,\beta}$, the subspace of \eqref{Papb} indexed only on $\cE$ (we are going to omit the subscript), while $W_{\cF}$ denotes the possible finite-dimensional normal block indexed on $\cF$.  If $\cF=\varnothing$, all terms involving $W_{\cF}$ are omitted. 

The KAM coordinates are endowed with the symplectic form
\begin{equation}\label{ckg:KAM-symplectic-form}
 \Omega_{\rm KAM} = \sum_{i=1}^2d\theta_i\wedge dy_i +\iu\sum_{a\in\mathcal E}d\bar z_a\wedge dz_a +\iu\sum_{a\in\mathcal F}d\bar z_a\wedge dz_a,
\end{equation}
where the last sum is omitted if $\mathcal F=\varnothing$. Then, for any analytic Hamiltonian $F$ on $ \cP_{\textnormal{KAM}}^{a,p,\beta}$ we have
\begin{equation}\label{ckg:XF-KAM}
 X_F = \left( -\partial_yF,\, \partial_\theta F,\, -\iu\,\partial_{\bar z}F,\, \iu\,\partial_zF,\, X_F^{\cF} \right),
\end{equation}
where, if the finite block is present and $W_{\cF}=(z_{\cF},\bar z_{\cF})$,
\[
 X_F^{\cF} = \left( -\iu\,\partial_{\bar z_{\cF}}F,\, \iu\,\partial_{z_{\cF}}F \right).
\]

\medskip

For $q>1/2$, $\beta \in [0,1] $ and $s,r>0$, define the complex tubular domain
\begin{equation*}
 \mathfrak D_{q,\beta}(s,r):= \left\{ (\theta,y,z,\bar z,W_{\cF}):
 \begin{array}{l}
 |\operatorname{Im}\theta|_\infty<s,\\[0.2em]
 |y|<r^2,\\[0.2em]
 \|z\|_{a,q,\beta} +\|\bar z\|_{a,q,\beta} +|W_{\cF}|<r
 \end{array}
 \right\} \subset \cP_{\textnormal{KAM}}^{a,q,\beta}.
\end{equation*}
where we will denote $\mathfrak D_q(s,r):=\mathfrak D_{q,0}(s,r).$ All finite-dimensional norms without a subscript are Euclidean. 

For a vector $X=(X^\theta,X^y,X^z,X^{\bar z},X^{\cF})$ on $\cP_{\textnormal{KAM}}^{a,p,\beta}$ and $r>0$ we define the rescaled norm as
\begin{equation}\label{ckg:scaled-X-norm}
 \|X\|_{r,q,\beta} := |X^\theta| +\frac{|X^y|}{r^2} +\frac{ \|X^z\|_{a,q,\beta} +\|X^{\bar z}\|_{a,q,\beta} +|X^{\cF}|} {r}.
\end{equation}
Again, the $X^{\cF}$ term is omitted when the finite block is absent.

Consider $\Xi\subseteq\Xi_0$ and a vector field
\[
 X=(X^\theta,X^y,X^z,X^{\bar z},X^{\cF}) \, : \Xi \times \mathfrak{D}_{q_1,\beta_1}(s,r) \to \cP_{\textnormal{KAM}}^{a,q_2,\beta_2},
\]
with $q_1,q_2>1/2$ and $\beta_1,\beta_2 \in [0,1]$, we define its scaled norm as
\begin{align*}
 \|X\|_{r,q_2,\beta_2}^{\mathfrak D_{q_1,\beta_1}(s,r)} &:= \sup_{\rho\in\Xi} \sup_{w\in\mathfrak D_{q_1,\beta_1}(s,r)} \|X(\rho;w)\|_{r,q_2,\beta_2},
 \\
 \|X\|_{r,q_2,\beta_2}^{\lip,\mathfrak D_{q_1,\beta_1}(s,r)} &:= \sup_{\substack{\rho,\rho'\in\Xi\\ \rho\neq\rho'}} \frac{ \displaystyle \sup_{w\in\mathfrak D_{q_1,\beta_1}(s,r)} \|X(\rho;w)-X(\rho';w)\|_{r,q_2,\beta_2}} {|\rho-\rho'|}.
\end{align*}
For $\alpha>0$ we use the weighted Lipschitz norm
\begin{equation*}
 \|X\|_{r,q_2,\beta_2}^{\mathcal L(\alpha), \mathfrak D_{q_1,\beta_1}(s,r)} := \|X\|_{r,q_2,\beta_2}^{ \mathfrak D_{q_1,\beta_1}(s,r)} + \alpha \|X\|_{r,q_2,\beta_2}^{\lip, \mathfrak D_{q_1,\beta_1}(s,r)}.
\end{equation*}

For an analytic Hamiltonian $F$ on $\mathfrak D_{q,\beta}(s,r)$, depending Lipschitz-continuously on $\rho\in\Xi_0$, its Fourier-Taylor expansion has the form
\begin{equation}\label{ckg:Fourier-Taylor-expansion}
 F= \sum_{k,\tt{m}, \varpi,\varkappa,\gamma} F_{k,\tt{m},\varpi,\varkappa,\gamma}(\rho)\, e^{\iu k\cdot\theta} y^{\tt{m}} z^\varpi\bar z^{\varkappa}W_{\cF}^\gamma ,
\end{equation}
where $k\in\mathbb Z^2, \, \, \tt{m} \in\mathbb N^2$, the multi-indices $\varpi, \varkappa$ are finitely supported on the infinite-dimensional normal set, and $\gamma$ is a finite-dimensional multi-index associated with $W_{\cF}$.  We assign the weighted degrees
\begin{equation}\label{ckg:weighted-degree}
 \deg\theta=\deg\rho=0, \qquad \deg y=2, \qquad \deg z=\deg\bar z=\deg W_{\cF}=1.
\end{equation}
Thus the weighted degree of the monomial in \eqref{ckg:Fourier-Taylor-expansion} is
\begin{equation}\label{w:degree}
  2|\tt{m}|+|\varpi|+|\varkappa|+|\gamma|.
\end{equation}

\begin{definition}\label{ckg:def-KAM-jet-preliminary}
The KAM jet $F^T$ of $F$ is the Fourier-Taylor polynomial formed by all monomials of weighted degree at most two
\begin{equation*}
 F^T:=
 \sum_{\substack{k,\tt{m},\varpi,\varkappa,\gamma\\
2|\tt{m}|+|\varpi|+|\varkappa|+|\gamma|\le2}} F_{k,\tt{m},\varpi,\varkappa,\gamma}(\rho)\, e^{\iu k\cdot\theta} y^{\tt{m}} z^\varpi\bar z^{\varkappa}W_{\cF}^\gamma .
\end{equation*}
\end{definition}

\begin{remark} 
The KAM jet contains the constant and $y$-linear terms, the terms linear in the normal variables and all quadratic normal terms, whereas terms of the form $yz$, $yW_{\cF}$, $z^3$, $z^2W_{\cF}$ and all terms of weighted degree at least three have zero KAM jet. 

Moreover, let $q_1,q_2>1/2$, $\beta_1,\beta_2\in[0,1]$, $\alpha>0$, and $0<\eta<1/8$. For an analytic Hamiltonian $F$ with finite indicated vector-field norms, 
\begin{equation}\label{ckg:jet-analytic-estimates}
\begin{aligned}
\|X_{F^T}\|_{r,q_2,\beta_2}^
{\mathcal L(\alpha),\mathfrak D_{q_1,\beta_1}(s,r)}
&\le C\|X_F\|_{r,q_2,\beta_2}^
{\mathcal L(\alpha),\mathfrak D_{q_1,\beta_1}(s,r)},
\notag\\
\|X_{F-F^T}\|_{\eta r,q_2,\beta_2}^
{\mathcal L(\alpha),\mathfrak D_{q_1,\beta_1}(s,4\eta r)}
&\le C\eta
\|X_F\|_{r,q_2,\beta_2}^
{\mathcal L(\alpha),\mathfrak D_{q_1,\beta_1}(s,r)}.
\end{aligned}
\end{equation}
These estimates follow from Cauchy's formula and the weighted dilation (see \cite[Section~3]{poschel1996kam} and \cite[Remark~4.1]{BBG}).
\end{remark}

For scalar, vector or finite-dimensional matrix-valued maps $A=A(\rho)$ on $\Xi \subseteq \Xi_0$ we use the notation 
\[
\|A\|_\infty:=\sup_{\rho \in \Xi}\|A(\rho)\|, \qquad \|A\|_{\infty}^{\rm lip} := \sup_{\substack{\rho,\rho'\in\Xi\\ \rho\neq\rho'}} \frac{\|A(\rho)-A(\rho')\|}{|\rho-\rho'|}
\]
and then for $\alpha>0$
\begin{equation}\label{ckg:finite-Lip-convention}
 \|A\|^{\mathcal L(\alpha)} :=\|A\|_\infty+ \alpha \|A\|_\infty^{\rm lip},
\end{equation}
where $\|\cdot\|$ is the absolute value, Euclidean norm or operator norm, according to the object under consideration.

\subsubsection{The unperturbed invariant torus: identical mode excitation}

Fix $m\in\mathbb Z$ and take
\[
  \mathcal J_{m,m}:=\{(m,+1),(m,-1)\}.
\]
For $\rho=(\rho_1,\rho_2)\in\Xi_0$ we introduce action--angle coordinates on the two internal variables by
\begin{equation}\label{ckg:same-AA}
\begin{aligned}
 \xi_m&=\sqrt{\rho_1+x_1}\,e^{\iu\theta_1},
 &\bar\xi_m&=\sqrt{\rho_1+x_1}\,e^{-\iu\theta_1},\\
 \eta_m&=\sqrt{\rho_2+x_2}\,e^{\iu\theta_2}, &\bar\eta_m&=\sqrt{\rho_2+x_2}\,e^{-\iu\theta_2}.
\end{aligned}
\end{equation}
On the real phase space the corresponding linear torus is
\begin{equation*}
  \mathbb T^2_m(\rho) :=\left\{ |\xi_m|^2=\rho_1,\ |\eta_m|^2=\rho_2, \ \xi_k=\eta_k=0\quad(k\neq m) \right\}.
\end{equation*}
Restricting the partial Birkhoff normal form to $\mathcal J_{m,m}$, we set
\[
 \Lambda_+^{(m,m)} :=\Pi_{\mathcal{IR}\cap\mathcal L_{\mathcal J_{m,m}}}P, \qquad \Lambda_+^{{\rm NLS},(m,m)} :=\Pi_{\mathcal{IR}\cap\mathcal L_{\mathcal J_{m,m}}}P^{\rm NLS}.
\]
By \eqref{ckg:Z4-closed}, this is exactly the part of the resonant quartic Hamiltonian containing at least one internal variable, namely
\begin{equation}\label{ckg:Lambda-plus-same}
 \Lambda_+^{(m,m)}= \frac{\mathsf P_m^2+2|\mathsf C_m|^2}
 {8\pi\widetilde w_m^2}+\sum_{k\neq m} \frac{\mathsf P_m\mathsf P_k +2\Re(\mathsf C_m\overline{\mathsf C_k})} {2\pi\widetilde w_m\widetilde w_k}.
\end{equation}
In particular, no resonant monomial of odd normal degree occurs in this configuration.

To eliminate the tangential angles from the quadratic normal part we introduce
\begin{equation}\label{ckg:same-angle-change}
\begin{cases}
 a_k=\xi_k e^{-\iu\theta_1},\quad \bar a_k=\bar\xi_k e^{\iu\theta_1},\quad b_k=\eta_k e^{-\iu\theta_2},\quad
 \bar b_k=\bar\eta_k e^{\iu\theta_2},& k\neq m,\\[0.3em]
 y_1=x_1+\displaystyle\sum_{k\neq m}|a_k|^2,\qquad y_2=x_2+\displaystyle\sum_{k\neq m}|b_k|^2.&
\end{cases}
\end{equation}
We denote by $\Psi_m(\rho)$ the inverse change of coordinates
\begin{equation}\label{ckg:Psi-same}
\begin{aligned}
 x_1&=y_1-\sum_{k\neq m}|a_k|^2, &\qquad \xi_k&=a_ke^{\iu\theta_1},\quad
 \bar\xi_k=\bar a_ke^{-\iu\theta_1},\quad k\neq m,\\
 x_2&=y_2-\sum_{k\neq m}|b_k|^2, & \eta_k&=b_ke^{\iu\theta_2},\quad \bar\eta_k=\bar b_ke^{-\iu\theta_2},\quad k\neq m,
\end{aligned}
\end{equation}
together with \eqref{ckg:same-AA}.  We set, for $k \neq m$,
\[
z_{(k,+1)}=a_k, \qquad z_{(k,-1)}=b_k, \qquad \bar z_{(k,+1)}=\bar a_k, \qquad \bar z_{(k,-1)}=\bar b_k,
\]
We then consider $(\theta,y,z,\bar z)\in\cP_{\textnormal{KAM}}^{a,p,\beta}$ as in \eqref{ckg:KAM-variables}, with $\cF=\varnothing$. This is the analogue, in the coordinates \eqref{ckg:same-angle-change}, of the angle elimination used for the equal-mode torus in \cite[Section~4.2]{GVdR}.

The quadratic generators \eqref{eq:NQM-zeta}, expressed through \eqref{ckg:same-angle-change}, take the form
\begin{align}
 \mathcal N&=\rho_1+\rho_2+y_1+y_2,
 \label{ckg:sym-same-N}\\
 \mathcal Q&=\rho_1-\rho_2+y_1-y_2,
 \label{ckg:sym-same-Q}\\
 \mathcal M&=m(\rho_1+\rho_2+y_1+y_2) +\sum_{k\neq m}(k-m)(|a_k|^2+|b_k|^2).
 \label{ckg:sym-same-M}
\end{align}
Equivalently, the two local masses are $\rho_1+y_1$ and $\rho_2+y_2$. Hence every locally gauge invariant Hamiltonian is independent of $\theta_1,\theta_2$.

We have the following result.
\begin{proposition}[Identical-mode effective Hamiltonian]
\label{ckg:prop-same-effective}
There exist $h_*,r_*,s_*,C_*,C>0$, depending on $a,p$ and the fixed tangential sites, such that the following holds for $0<h<h_*$, $0<r_0<r_*$, $q=p,p-4$, $\beta\in[0,1]$ and $\rho\in\Xi_0$.  For any $0<r\leq r_0$ and $0<s\leq s_*$ the map $\Psi_m(\rho)$ defined in \eqref{ckg:Psi-same} is real analytic and symplectic, and
\begin{equation}\label{psi:m}
 \Psi_m(\rho):\mathfrak D_{q,\beta}(s,r) \longrightarrow B_{a,q,\beta}(C_*\tR).
\end{equation}
For the Klein-Gordon Hamiltonian we have
\begin{equation*}
 (\Lambda+\Lambda_+^{(m,m)})\circ\Psi_m =C_{m,m}(\rho;c)+H_{\rm eff}^{(m,m)}+\mathcal R_{m,m},
\end{equation*}
with
\begin{align}
\label{ckg:C-same}
&C_{m,m}(\rho;c) =\lambda_m(\rho_1+\rho_2) +\frac{\rho_1^2+\rho_2^2+4\rho_1\rho_2} {8\pi\widetilde w_m^2},& \\
\label{ckg:Heff-same}
&H_{\rm eff}^{(m,m)}=
 \langle\omega^{(m,m)}(\rho;c),y\rangle+\sum_{k\neq m} \Bigl( \alpha_k|a_k|^2+\beta_k|b_k|^2 +\gamma_k(a_k\bar b_k+\bar a_k b_k) \Bigr),& \\
\label{ckg:omega-same}
 &\omega_1^{(m,m)}=\lambda_m+
 \frac{\rho_1+2\rho_2}{4\pi\widetilde w_m^2},\qquad
 \omega_2^{(m,m)}=\lambda_m+ \frac{\rho_2+2\rho_1}{4\pi\widetilde w_m^2},&
\end{align}
where the Hermitian matrix
\begin{equation}\label{ckg:Mk-same}
 M_k(\rho;c):=
 \begin{pmatrix}\alpha_k&\gamma_k\\ \gamma_k&\beta_k\end{pmatrix}
\end{equation}
has coefficients
\begin{equation}\label{ckg:coeff-same}
\begin{aligned}
 &\alpha_k=\lambda_k-\lambda_m +\frac{\rho_1+\rho_2}{2\pi\widetilde w_m\widetilde w_k}
 -\frac{\rho_1+2\rho_2}{4\pi\widetilde w_m^2}, \qquad \gamma_k=\frac{\sqrt{\rho_1\rho_2}} {2\pi\widetilde w_m\widetilde w_k}& \\
 &\beta_k=\lambda_k-\lambda_m +\frac{\rho_1+\rho_2}{2\pi\widetilde w_m\widetilde w_k}
 -\frac{\rho_2+2\rho_1}{4\pi\widetilde w_m^2},&
\end{aligned}
\end{equation}
and the remainder $\mathcal{R}_{m,m}$ has zero KAM jet (see Definition \ref{ckg:def-KAM-jet-preliminary}), it is locally gauge invariant and satisfies
\begin{equation}
 \|X_{\mathcal R_{m,m}}\|_{r_0,q,0}^{\mathcal L(\varrho), \mathfrak D_q(s_*,r_0)} \le Cr_0^2.
 \label{ckg:R-same-size}
\end{equation} 
For the NLS Hamiltonian, with the same map $\Psi_m(\rho)$, it holds
\begin{equation}\label{ckg:same-NLS-decomposition}
 (\Lambda^{\rm NLS}+\Lambda_+^{{\rm NLS},(m,m)})\circ\Psi_m =C_{m,m}^{\rm NLS}(\rho)+H_{\rm eff}^{{\rm NLS},(m,m)} +\mathcal R_{m,m}^{\rm NLS},
\end{equation}
with
\begin{align}
 C_{m,m}^{\rm NLS}(\rho) ={}&\frac{m^2}{2}(\rho_1+\rho_2) +\frac{\rho_1^2+\rho_2^2+4\rho_1\rho_2}{8\pi} \nonumber \\
\label{ckg:Heff-same-NLS} 
 H_{\rm eff}^{{\rm NLS},(m,m)}={}&
 \langle\omega^{(m,m),\rm NLS}(\rho),y\rangle\\
 &+\sum_{k\neq m} \Bigl( \alpha_k^{\rm NLS}|a_k|^2+\beta_k^{\rm NLS}|b_k|^2 +\gamma_k^{\rm NLS}(a_k\bar b_k+\bar a_k b_k) \Bigr), \nonumber \\
 \omega_1^{(m,m),\rm NLS}
 ={}&\frac{m^2}{2}+\frac{\rho_1+2\rho_2}{4\pi}, \quad \, \, \,
 \omega_2^{(m,m),\rm NLS} =\frac{m^2}{2}+\frac{\rho_2+2\rho_1}{4\pi}, \nonumber
\end{align}
where the Hermitian matrix
\begin{equation*}
 M_k^{\rm NLS}(\rho):=
 \begin{pmatrix}\alpha_k^{\rm NLS}&\gamma_k^{\rm NLS}\\ \gamma_k^{\rm NLS}&\beta_k^{\rm NLS}\end{pmatrix}
\end{equation*}
has coefficients
\begin{equation}
\label{ckg:coeff-same-NLS}
 \alpha_k^{\rm NLS}=\frac{k^2-m^2}{2}+\frac{\rho_1}{4\pi},\qquad
 \beta_k^{\rm NLS}=\frac{k^2-m^2}{2}+\frac{\rho_2}{4\pi},\qquad
 \gamma_k^{\rm NLS}=\frac{\sqrt{\rho_1\rho_2}}{2\pi},
\end{equation}
and the remainder $\mathcal{R}^{\rm NLS}_{m,m}$ has zero KAM jet, it is locally gauge invariant and satisfies 
\begin{equation}
 \|X_{\mathcal R_{m,m}^{\rm NLS}}\|_{r_0,q}^{\mathcal L(\varrho), \mathfrak D_q(s_*,r_0)} \le Cr_0^2.
 \label{ckg:R-same-NLS-size}
\end{equation}
Finally it holds
\begin{align}
 \|X_{\mathcal R_{m,m}-\mathcal R_{m,m}^{\rm NLS}}\|_{r_0,p-4,0}^{\mathfrak D_p(s_*,r_0)} &\le Chr_0^2.
 \label{ckg:R-same-R-close} \\
 \bigl|\omega^{(m,m)}-c^2(1,1) -\omega^{(m,m),\rm NLS}\bigr|^{\mathcal L(\varrho)} &\le Ch,
 \label{ckg:omega-same-close}\\
 \|M_k-M_k^{\rm NLS}\|^{\mathcal L(\varrho)} &\le Ch\langle k\rangle^4, \qquad k\neq m.
 \label{ckg:M-same-close}
\end{align}
\end{proposition}

\begin{proof}
Put  $A=(A_1,A_2)$ with
\[
 A_1=\sum_{k\ne m}|a_k|^2,\qquad  A_2=\sum_{k\ne m}|b_k|^2,\qquad x=y-A,\qquad I=\rho+x.
\]
For $0<r\le r_0$ and $r_0$ sufficiently small, on $\mathfrak D_{q,\beta}(s,r)$ we have
\begin{equation}\label{ckg:Psi-same-domain}
 |x_i|\le 2r^2\le\frac12\rho_i,\qquad i=1,2,
\end{equation}
since $\rho_i\asymp\varrho=r_0^{4/3}$. Thus the square roots in \eqref{ckg:same-AA} are analytic. Moreover, the phase factors are bounded by $e^s$, and $r_0=o(\sqrt\varrho)$. Hence, uniformly for $q=p,p-4$ and $\beta\in[0,1]$,
\[
 \|\Psi_m(\rho)(\theta,y,z,\bar z)\|_{a,q,\beta}  \le C_*\sqrt\varrho=C_*\tR,
\]
making the compositions with the first Birkhoff map well defined upon applying Theorem~\ref{Teorema0} on a ball of radius $C_*\tR$, which only changes the constants in its estimates. The map is independent of $h$ and is symplectic: indeed, each normal phase rotation contributes $d\theta_i\wedge dA_i$ to the symplectic form, which is cancelled by $dx_i=dy_i-dA_i$, giving 
\begin{equation*}
 \sum_{i=1}^2d\theta_i\wedge dx_i +\iu\sum_{k\neq m} (d\bar\xi_k\wedge d\xi_k+d\bar\eta_k\wedge d\eta_k) = \sum_{i=1}^2d\theta_i\wedge dy_i +\iu\sum_{k\neq m} (d\bar a_k\wedge da_k+d\bar b_k\wedge db_k).
\end{equation*}

The quadratic Hamiltonian becomes
\[
 \Lambda\circ\Psi_m =\lambda_m(\rho_1+\rho_2+y_1+y_2) +\sum_{k\ne m}(\lambda_k-\lambda_m)(|a_k|^2+|b_k|^2).
\]
The purely internal quartic part in $I$ of \eqref{ckg:Lambda-plus-same} is
\[
 U(I)=\frac12 \langle I, T_mI\rangle,\qquad T_m=\frac1{4\pi\widetilde w_m^2}
 \begin{pmatrix}
1&2\\
2&1
\end{pmatrix}.
\]
Since
\[
 U(\rho+x)=U(\rho)+(T_m\rho)\cdot x +\frac12\langle x,T_mx \rangle,
\]
its constant and $y$-linear contributions give \eqref{ckg:C-same} and \eqref{ckg:omega-same}. The substitution $x=y-A$ also produces the
quadratic normal terms $-(T_m\rho)_1|a_k|^2-(T_m\rho)_2|b_k|^2$. On the other hand, the cross-index quartic terms, evaluated at $I=\rho$, contribute
\[
 \frac{(\rho_1+\rho_2)(|a_k|^2+|b_k|^2) +\sqrt{\rho_1\rho_2}(\bar a_kb_k+a_k\bar b_k)}{2\pi\widetilde w_m\widetilde w_k}.
\]
Together with the quadratic Hamiltonian, these give exactly \eqref{ckg:Heff-same}--\eqref{ckg:coeff-same}; in particular, the tangential angles cancel from the exchange term.

Subtracting the constant, tangential linear and normal quadratic terms just computed, we obtain
\[
\begin{aligned}
 \mathcal R_{m,m}={}&\frac12\langle x,T_mx \rangle\\
 &+\sum_{k\ne m} \frac{(x_1+x_2)(|a_k|^2+|b_k|^2)+(\sqrt{I_1I_2}-\sqrt{\rho_1\rho_2})(\bar a_kb_k+a_k\bar b_k)}{2\pi\widetilde w_m\widetilde w_k}.
\end{aligned}
\]
Here $x$ has weighted degree two, and $\sqrt{I_1I_2}-\sqrt{\rho_1\rho_2}=O(|x|)$. Thus every term has weighted degree \eqref{w:degree} at least four, so $(\mathcal R_{m,m})^T=0$. Moreover, \eqref{ckg:Psi-same-domain} implies
\[
 |x|+ |\sqrt{I_1I_2}-\sqrt{\rho_1\rho_2}| \lesssim r^2,\qquad |\partial_{I_i}\sqrt{I_1I_2}|\lesssim1.
\]
Differentiating the remainder, using $\widetilde w_k^{-1}\le1$, yields
\begin{equation}\label{loss:regula}
 |X_{\mathcal R_{m,m}}^\theta|\lesssim r^2, \qquad X_{\mathcal R_{m,m}}^y=0,\qquad \|X_{\mathcal R_{m,m}}^z\|_{a,q} +\|X_{\mathcal R_{m,m}}^{\bar z}\|_{a,q} \lesssim r^3.
\end{equation}
Since 
\[
\partial_{\rho}(\sqrt{I_1I_2}-\sqrt{\rho_1 \rho_2})=O(r^2/\varrho)
\]
we obtain also the weighted Lipschitz norm, proving \eqref{ckg:R-same-size}. Notice that $-\langle y,T_mA\rangle$ has a nondecaying normal derivative: the pullback therefore does not preserve the smoothing obtained for the first term in \eqref{res.0.2}.

For NLS, the same calculation applies with $\widetilde w_j=1$ and $\lambda_j=j^2/2$, giving \eqref{ckg:same-NLS-decomposition}--\eqref{ckg:coeff-same-NLS} and, by the same degree and vector-field estimates, $(\mathcal R_{m,m}^{\rm NLS})^T=0$ and \eqref{ckg:R-same-NLS-size}.
Both remainders are independent of the angles and hence locally gauge invariant.

Finally, with $T_0=(4\pi)^{-1} \left(\begin{smallmatrix}1&2\\2&1\end{smallmatrix}\right)$,
\[
 \|T_m-T_0\|\lesssim h,\qquad \left|\frac1{\widetilde w_m\widetilde w_k}-1\right| \lesssim h\langle k\rangle^2,\qquad \left|\nu_k(h)-\frac{k^2}{2}\right| \lesssim h\langle k\rangle^4,
\]
since $m$ is fixed. Comparing the two explicit remainder formulas, the first two bounds and the Sobolev loss from $p$ to $p-4$ give \eqref{ckg:R-same-R-close}. The explicit frequency and matrix coefficients similarly yield \eqref{ckg:omega-same-close} and \eqref{ckg:M-same-close}. These comparisons also hold in the stated weighted parameter Lipschitz norms.
\end{proof}

\begin{remark}[Block structure in the identical-mode case]

The matrices $M_k$ are Hermitian and could be diagonalized for each fixed parameter.  We do not need such a parameter-dependent diagonalization in the KAM iteration: retaining the $2\times2$ species matrix, and later the $\{k,-k\}$ multiplicity, as a single uniformly finite Hermitian block keeps the parameter dependence explicit and avoids introducing a separate choice of eigenvectors.
\end{remark}

\subsubsection{The unperturbed invariant torus: distinct mode excitation}

Let $m\neq n$ and take
\[
 \mathcal J_{m,n}:=\{(m,+1),(n,-1)\}.
\]
The internal action--angle variables are for $ \rho\in\Xi_0$
\begin{equation*}
\begin{aligned}
 \xi_m&=\sqrt{\rho_1+x_1}\,e^{\iu\theta_1},
 &\bar\xi_m&=\sqrt{\rho_1+x_1}\,e^{-\iu\theta_1},\\
 \eta_n&=\sqrt{\rho_2+x_2}\,e^{\iu\theta_2}, &\bar\eta_n&=\sqrt{\rho_2+x_2}\,e^{-\iu\theta_2}.
\end{aligned}
\end{equation*}
The corresponding linear torus is
\begin{equation*}
 \mathbb T^2_{m,n}(\rho) :=\left\{ |\xi_m|^2=\rho_1,\ |\eta_n|^2=\rho_2, \ \xi_k=0 \, \, \forall k\neq m, \, \, \,  \eta_k=0 \, \, \forall k\neq n \right\}.
\end{equation*}
The two normal variables $\eta_m$ and $\xi_n$ will be called \emph{semi-external modes}.  As in the unstable-torus construction of \cite[Section~4.1]{GVdR}, they are the only normal modes which have to be rotated in order to remove the tangential angles from the quadratic normal form.  We set
\begin{equation}\label{ckg:distinct-angle-change}
\begin{cases}
 a_n=\xi_n e^{-\iu\theta_2},\quad \bar a_n=\bar\xi_n e^{\iu\theta_2},\quad b_m=\eta_m e^{-\iu\theta_1},\quad
 \bar b_m=\bar\eta_m e^{\iu\theta_1},\\[0.2em]
 a_k=\xi_k,\quad \bar a_k=\bar\xi_k,\quad
 b_k=\eta_k,\quad \bar b_k=\bar\eta_k,\qquad k\neq m,n,\\[0.2em]
 y_1=x_1+|b_m|^2,\qquad y_2=x_2+|a_n|^2.
\end{cases}
\end{equation}
Let $\Psi_{m,n}(\rho)$ denote the inverse map, then
\begin{equation}\label{ckg:Psi-dist}
\begin{aligned}
 x_1&=y_1-|b_m|^2, & x_2&=y_2-|a_n|^2,\\
 \eta_m&=b_me^{\iu\theta_1}, &
 \bar\eta_m&=\bar b_me^{-\iu\theta_1},\\
 \xi_n&=a_ne^{\iu\theta_2}, & \bar\xi_n&=\bar a_ne^{-\iu\theta_2}.
\end{aligned}
\end{equation}
The generic normal variables are unchanged.  We set $\cF:=\{(n,+1),(m,-1)\}$ and
\[
 W_{\cF}:=(a_n,b_m,\bar a_n,\bar b_m), \qquad z_{(j,+1)}:=\xi_j,\quad z_{(j,-1)}:=\eta_j, \qquad j\ne m,n,
\]
and analogously for $\bar z$.  We then consider $(\theta,y,z,\bar z,W_{\cF})\in\cP_{\textnormal{KAM}}^{a,p,\beta}$ as in \eqref{ckg:KAM-variables}.

Restricting the partial Birkhoff normal form to $\mathcal J_{m,n}$, we put
\begin{equation*}
 \Lambda_+^{(m,n)} :=\Pi_{\mathcal{IR}\cap\mathcal L_{\mathcal J_{m,n}}}P, \qquad \Lambda_+^{{\rm NLS},(m,n)} :=\Pi_{\mathcal{IR}\cap\mathcal L_{\mathcal J_{m,n}}}P^{\rm NLS}.
\end{equation*}
The quadratic generators \eqref{eq:NQM-zeta} take the following form in the distinct-mode coordinates:
\begin{align}
 \mathcal N={}&\rho_1+\rho_2+y_1+y_2 +\sum_{k\neq m,n}(|a_k|^2+|b_k|^2),
 \label{ckg:sym-dist-N}\\
 \mathcal Q={}&\rho_1-\rho_2+y_1-y_2 +2|a_n|^2-2|b_m|^2 +\sum_{k\neq m,n}(|a_k|^2-|b_k|^2),
  \notag\\
 \mathcal M={}&m(\rho_1+y_1)+n(\rho_2+y_2) +\sum_{k\neq m,n}k(|a_k|^2+|b_k|^2).
 \label{ckg:sym-dist-M}
\end{align}
In particular, the semi-external variables do not occur in $\mathcal N$ or $\mathcal M$.  This cancellation is the same geometric mechanism as in the angle elimination of \cite[Section~4.1]{GVdR}.

We have the following result.

\begin{proposition}[Distinct-mode effective Hamiltonian]
\label{ckg:prop-dist-effective}
There exist $h_*,r_*,s_*,C_*,C>0$, depending only on $a,p$ and the fixed tangential sites, such that the following holds for $0<h<h_*$, $0<r_0<r_*$, $q=p,p-4$, $\beta\in[0,1]$ and $\rho\in\Xi_0$.  For every $0<r\le r_0$ and $0<s\le s_*$, the map $\Psi_{m,n}(\rho)$ defined in \eqref{ckg:Psi-dist} is real analytic and symplectic, and
\begin{equation}\label{psi_mn}
 \Psi_{m,n}(\rho):\mathfrak D_{q,\beta}(s,r) \longrightarrow B_{a,q,\beta}(C_*\tR).
\end{equation}
For the Klein-Gordon Hamiltonian we have
\begin{equation}\label{ckg:dist-decomposition}
 (\Lambda+\Lambda_+^{(m,n)})\circ\Psi_{m,n} =C_{m,n}(\rho;c)+H_{\rm eff}^{(m,n)}+Z_{4,3}+\mathcal R_{m,n},
\end{equation}
with
\begin{align}\label{ckg:C-dist}
 C_{m,n}(\rho;c) &=\lambda_m\rho_1+\lambda_n\rho_2 +\frac{\rho_1^2}{8\pi\widetilde w_m^2} +\frac{\rho_2^2}{8\pi\widetilde w_n^2} +\frac{\rho_1\rho_2}{2\pi\widetilde w_m\widetilde w_n},& \\
\label{ckg:Heff-dist}
 H_{\rm eff}^{(m,n)}&= \langle\omega^{(m,n)}(\rho;c),y\rangle
 +\sum_{k\neq m,n}\Lambda_k(\rho;c)(|a_k|^2+|b_k|^2)+H_\cF(a_n,b_m;\rho,c),& \\
\label{ckg:omega-dist}
 \omega_1^{(m,n)}&=\lambda_m+ \frac{\rho_1}{4\pi\widetilde w_m^2}
 +\frac{\rho_2}{2\pi\widetilde w_m\widetilde w_n}, \qquad  \omega_2^{(m,n)}=\lambda_n+ \frac{\rho_2}{4\pi\widetilde w_n^2} +\frac{\rho_1}{2\pi\widetilde w_m\widetilde w_n},& \\
\label{ckg:Lambda-dist}
 \Lambda_k(\rho;c)&=\lambda_k+ \frac1{2\pi\widetilde w_k} \left(\frac{\rho_1}{\widetilde w_m} +\frac{\rho_2}{\widetilde w_n}\right), \qquad k\neq m,n, \\
\label{ckg:HF-dist}
 H_\cF&=\alpha|a_n|^2+\beta|b_m|^2 +\gamma(a_nb_m+\bar a_n\bar b_m),& \\
\label{ckg:abg-dist}
 \alpha&=\frac{\rho_2}{4\pi\widetilde w_n^2},\qquad \beta=\frac{\rho_1}{4\pi\widetilde w_m^2},\qquad \gamma=\frac{\sqrt{\rho_1\rho_2}} {2\pi\widetilde w_m\widetilde w_n},& \\
\label{ckg:cubic-res-dist-newcoords}
 Z_{4,3}&=\sum_{k\neq m,n}\Re\!\left[ \left( \frac{\sqrt{\rho_2}}{\pi\widetilde w_n\widetilde w_k}\,a_n +\frac{\sqrt{\rho_1}}{\pi\widetilde w_m\widetilde w_k}\,\bar b_m \right)\bar a_k b_k\right],&
\end{align}
where the Hamiltonians $Z_{4,3}$ and $\mathcal R_{m,n}$ have zero KAM jet  (see Definition \ref{ckg:def-KAM-jet-preliminary}), are locally gauge invariant and satisfy
\begin{equation} \label{ckg:Z43-size}
 \|X_{Z_{4,3}}\|_{r_0,q,1}^{\mathcal L(\varrho), \mathfrak D_q(s_*,r_0)} \le Cr_0^{5/3}, \qquad \|X_{\mathcal R_{m,n}}\|_{r_0,q,1}^{\mathcal L(\varrho), \mathfrak D_q(s_*,r_0)} \le Cr_0^2.
\end{equation}
For the NLS Hamiltonian the same map $\Psi_{m,n}(\rho)$ gives
\begin{equation}\label{ckg:dist-NLS-decomposition}
 (\Lambda^{\rm NLS}+\Lambda_+^{{\rm NLS},(m,n)})\circ\Psi_{m,n} =C_{m,n}^{\rm NLS}(\rho)+H_{\rm eff}^{{\rm NLS},(m,n)} +Z_{4,3}^{\rm NLS}+\mathcal R_{m,n}^{\rm NLS},
\end{equation}
where
\begin{align}
 C_{m,n}^{\rm NLS}(\rho) &=\frac{m^2}{2}\rho_1+\frac{n^2}{2}\rho_2 +\frac{\rho_1^2+\rho_2^2}{8\pi} +\frac{\rho_1\rho_2}{2\pi},& \nonumber \\
\label{ckg:Heff-dist-NLS}
 H_{\rm eff}^{{\rm NLS},(m,n)}&= \langle\omega^{(m,n),\rm NLS}(\rho),y\rangle +\sum_{k\neq m,n}\Lambda_k^{\rm NLS}(\rho)
 (|a_k|^2+|b_k|^2)+H_\cF^{\rm NLS}(a_n,b_m),& \\
\label{ckg:dist-NLS-coeff}
 \omega^{(m,n),\rm NLS} &=\left(\frac{m^2}{2}+\frac{\rho_1}{4\pi}+\frac{\rho_2}{2\pi},\frac{n^2}{2}+\frac{\rho_2}{4\pi}+\frac{\rho_1}{2\pi} \right), \quad
 \Lambda_k^{\rm NLS}
 =\frac{k^2}{2}+\frac{\rho_1+\rho_2}{2\pi}, &\\
\label{ckg:HFNLS-dist}
 H_\cF^{\rm NLS} &=\alpha^{\rm NLS}|a_n|^2+\beta^{\rm NLS}|b_m|^2 +\gamma^{\rm NLS}(a_nb_m+\bar a_n\bar b_m),& \\
\alpha^{\rm NLS}&=\frac{\rho_2}{4\pi},\qquad \beta^{\rm NLS}=\frac{\rho_1}{4\pi},\qquad \gamma^{\rm NLS}=\frac{\sqrt{\rho_1\rho_2}}{2\pi},& \nonumber \\
\label{ckg:cubic-res-dist-NLS}
 Z_{4,3}^{\rm NLS} &=\sum_{k\neq m,n}\Re\!\left[ \left( \frac{\sqrt{\rho_2}}{\pi}\,a_n +\frac{\sqrt{\rho_1}}{\pi}\,\bar b_m \right)\bar a_k b_k\right],&
\end{align}
where the Hamiltonians $Z_{4,3}^{\rm NLS}$ and $\mathcal R_{m,n}^{\rm NLS}$ have zero KAM jet, are locally gauge invariant and satisfy
\begin{equation} \label{ckg:Z43-NLS-size}
 \|X_{Z_{4,3}^{\rm NLS}}\|_{r_0,q}^{\mathcal L(\varrho), \mathfrak D_q(s_*,r_0)} \le Cr_0^{5/3}, \quad \|X_{\mathcal R_{m,n}^{\rm NLS}}\|_{r_0,q}^{\mathcal L(\varrho), \mathfrak D_q(s_*,r_0)} \le Cr_0^2.
\end{equation}
Finally, it holds
\begin{align}
 &\bigl|\omega^{(m,n)}-c^2(1,1) -\omega^{(m,n),\rm NLS}\bigr|^{\mathcal L(\varrho)} \le Ch, \qquad |\Lambda_k-c^2-\Lambda_k^{\rm NLS}|^{\mathcal L(\varrho)} \le Ch\langle k\rangle^4,&
 \label{ckg:omega-dist-close}\\
 &|\alpha-\alpha^{\rm NLS}|^{\mathcal L(\varrho)} +|\beta-\beta^{\rm NLS}|^{\mathcal L(\varrho)} +|\gamma-\gamma^{\rm NLS}|^{\mathcal L(\varrho)} \le Ch\varrho,&
 \label{ckg:HF-dist-close}\\
& \|X_{Z_{4,3}-Z_{4,3}^{\rm NLS}}\|_{r_0,p-4}^{ \mathfrak D_p(s_*,r_0)} \le Ch r_0^{5/3}, \qquad \|X_{\mathcal R_{m,n}-\mathcal R_{m,n}^{\rm NLS}}\|_{r_0,p-4,0}^{\mathfrak D_p(s_*,r_0)} \le Chr_0^2.&
 \label{ckg:Z43-diff}
\end{align}
\end{proposition}

\begin{proof}
Set $A=|a_n|^2$, $B=|b_m|^2$, $x=(y_1-B,y_2-A)$, and $I=\rho+x$. As in Proposition~\ref{ckg:prop-same-effective}, $|x_i|\le 2r^2\le\rho_i/2$ on $\mathfrak D_{q,\beta}(s,r)$, for $r\le r_0$ sufficiently small. Hence the square roots are analytic. Since only the two semi-external modes are rotated, their phase factors are bounded by $e^s$, and
\begin{equation*}
\|\Psi_{m,n}(\rho)(\theta,y,z,\bar z,W_{\cF})\|_{a,q,\beta}\le C_*\sqrt\varrho=C_*\tR,\qquad q=p,p-4,\quad\beta\in[0,1].
\end{equation*}
Thus the compositions with the first Birkhoff normal form are well defined on the corresponding enlarged input ball. The map is symplectic: the phase rotations contribute $d\theta_1\wedge dB+d\theta_2\wedge dA$, cancelled by $dx_1=dy_1-dB$ and $dx_2=dy_2-dA$. In particular,
\begin{equation*}
\Lambda\circ\Psi_{m,n}=\lambda_m(\rho_1+y_1)+\lambda_n(\rho_2+y_2)+\sum_{k\ne m,n}\lambda_k(|a_k|^2+|b_k|^2).
\end{equation*}

If we write $d_{jk}:=(2\pi\widetilde w_j\widetilde w_k)^{-1}$, the quartic formula \eqref{ckg:Z4-closed} can be evaluated using
\begin{equation*}
\mathsf P_m=I_1+B,\quad \mathsf P_n=I_2+A,\quad \mathsf C_m=\sqrt{I_1}\,\bar b_m,\quad \mathsf C_n=\sqrt{I_2}\,a_n,\quad \mathsf C_k=a_k\bar b_k\ (k\ne m,n).
\end{equation*}
Its purely internal part in $I$ reads as
\begin{equation}\label{ckg:Zint-dist}
U(I)=\tfrac14d_{mm}I_1^2+\tfrac14d_{nn}I_2^2+d_{mn}I_1I_2.
\end{equation}
By denoting $D=(D_1,D_2)=\nabla U(\rho)$, the expansion of $U(\rho+x)$ gives \eqref{ckg:C-dist} and \eqref{ckg:omega-dist}. For the generic modes the mixed action terms contribute $(d_{mk}\rho_1+d_{nk}\rho_2)(|a_k|^2+|b_k|^2)$; these modes are not rotated, proving \eqref{ckg:Lambda-dist}.

The free semi-external quadratic frequencies cancel, while the remaining quadratic terms in $a_n,b_m$ are
\begin{equation*}
\begin{aligned}
&(d_{nn}\rho_2+d_{mn}\rho_1-D_2)A+(d_{mm}\rho_1+d_{mn}\rho_2-D_1)B\\
&\hspace{3em}+2d_{mn}\sqrt{\rho_1\rho_2}\,\Re(a_nb_m),
\end{aligned}
\end{equation*}
which gives \eqref{ckg:HF-dist}--\eqref{ckg:abg-dist}. The exchanges $\mathsf C_n\overline{\mathsf C_k}$ and $\mathsf C_m\overline{\mathsf C_k}$ produce \eqref{ckg:cubic-res-dist-newcoords} upon replacing $\sqrt{I_i}$ by $\sqrt{\rho_i}$.

For completeness, write $J_k=|a_k|^2+|b_k|^2$ and $Q_k=\bar a_kb_k$. Subtracting the terms computed above yields
\begin{equation*}
\begin{aligned}
\mathcal R_{m,n}={}&\tfrac12\langle x,D^2Ux\rangle+d_{mm}x_1B+d_{nn}x_2A+\tfrac14(d_{mm}B^2+d_{nn}A^2)\\
&+d_{mn}(x_1A+x_2B+AB)+2d_{mn}(\sqrt{I_1I_2}-\sqrt{\rho_1\rho_2})\Re(a_nb_m)\\
&+\sum_{k\ne m,n}\Bigl[(d_{mk}y_1+d_{nk}y_2)J_k\\
&\hspace{5em}+2\Re\!\bigl((d_{mk}(\sqrt{I_1}-\sqrt{\rho_1})\bar b_m+d_{nk}(\sqrt{I_2}-\sqrt{\rho_2})a_n)Q_k\bigr)\Bigr].
\end{aligned}
\end{equation*}
Since $x$ has weighted degree two and the square-root differences vanish at $x=0$, both $\mathcal R_{m,n}$ and $Z_{4,3}$ have zero KAM jet. On the above domain,
\begin{equation*}
|x|+|\sqrt{I_1I_2}-\sqrt{\rho_1\rho_2}|\lesssim r^2,\qquad |\sqrt{I_i}-\sqrt{\rho_i}|\lesssim r^2/\sqrt\varrho.
\end{equation*}
Every generic-mode coefficient carries $\widetilde w_k^{-1}$, which compensates the output smoothing weight. Differentiating the formula gives, for $q=p,p-4$,
\begin{equation*}
\begin{aligned}
|X_{\mathcal R_{m,n}}^\theta|&\lesssim r^2+r^3/\sqrt\varrho,\qquad X_{\mathcal R_{m,n}}^y=0,\\
\|X_{\mathcal R_{m,n}}^z\|_{a,q,1}+\|X_{\mathcal R_{m,n}}^{\bar z}\|_{a,q,1}+|X_{\mathcal R_{m,n}}^{\cF}|&\lesssim r^3+r^4/\sqrt\varrho.
\end{aligned}
\end{equation*}
The same orders hold for the $\varrho$-weighted parameter Lipschitz seminorms by computing the derivatives. As $r\le r_0=o(\sqrt\varrho)$, this proves the remainder bound in \eqref{ckg:Z43-size}. Moreover, the cubic coefficients are $O(\sqrt\varrho\,\widetilde w_k^{-1})$, including their weighted Lipschitz bounds, and hence
\begin{equation*}
\|X_{Z_{4,3}}\|_{r_0,q,1}^{\mathcal L(\varrho),\mathfrak D_q(s_*,r_0)}\lesssim\sqrt\varrho\,r_0=r_0^{5/3}.
\end{equation*}

The NLS identities follow from the same calculation with $\widetilde w_j=1$ and $\lambda_j=j^2/2$, yielding \eqref{ckg:dist-NLS-decomposition}--\eqref{ckg:cubic-res-dist-NLS}, the zero jets and \eqref{ckg:Z43-NLS-size}. Finally, for fixed $m,n$,
\begin{equation*}
\left|\nu_k(h)-\frac{k^2}{2}\right|\lesssim h\langle k\rangle^4,\qquad \left|\frac1{\widetilde w_j\widetilde w_k}-1\right|\lesssim h(\langle j\rangle^2+\langle k\rangle^2).
\end{equation*}
These bounds applied to the explicit coefficients give \eqref{ckg:omega-dist-close}--\eqref{ckg:HF-dist-close}, including their weighted parameter Lipschitz estimates. Applied to the cubic and remainder formulas, the same bounds, together with the loss from $p$ to $p-4$, yield both comparison estimates in \eqref{ckg:Z43-diff}. Indeed, the quartic Hamiltonian \eqref{ckg:Z4-closed} and all terms extracted above are invariant under the transported local-mass flows, so $\mathcal R_{m,n}$ is locally gauge invariant.
\end{proof}

We now isolate the finite block responsible for instability.  Put
\begin{equation}\label{ckg:UF-def}
 U:=\binom{a_n}{\bar b_m}, \qquad W_{\cF}:=(a_n,b_m,\bar a_n,\bar b_m).
\end{equation}
The Hamilton equations generated by $H_\cF$  and $H_\cF^{\rm NLS}$ in \eqref{ckg:HF-dist} and \eqref{ckg:HFNLS-dist} on the $U$-subspace are $\dot U=L_\cF^{(2)}U$ and $\dot U=L_{\cF, \rm NLS}^{(2)}U$, with
\begin{equation}\label{ckg:LF-def}
 L_{\cF}^{(2)}:=
 \begin{pmatrix}
 -\iu\alpha&-\iu\gamma\\
 \iu\gamma&\iu\beta
 \end{pmatrix},
 \qquad L_{\cF,\rm NLS}^{(2)}:=
 \begin{pmatrix}
 -\iu\alpha^{\rm NLS}&-\iu\gamma^{\rm NLS}\\
 \iu\gamma^{\rm NLS}&\iu\beta^{\rm NLS}
 \end{pmatrix}.
\end{equation}
Let 
\begin{equation}\label{ckg:sigma-NLS-def}
 \delta^{\sharp}:=\frac{\alpha^{\sharp}-\beta^{\sharp}}{2}, \qquad (\sigma^{\sharp})^2 :=(\gamma^{\sharp})^2 -\left(\frac{\alpha^{\sharp}+\beta^{\sharp}}{2}\right)^2,
\end{equation}
where from now on when we use $\sharp$ it stands for either the KG or the NLS considered quantities. We have the following 
\begin{lemma}
\label{ckg:lem-hyperbolic}
There exist \(h_*,\kappa_{\cF}>0\) such that, for every \(0<r_0<r_*\), where \(r_*\) is as in Proposition \ref{ckg:prop-dist-effective}, for any $0<h<h_*$ and for any $\rho\in\Xi_0$ it holds $(\sigma^{\sharp})^2>0$, for both the equations. By setting $\sigma^\sharp$ as the positive square root, the full finite-dimensional Hamiltonian blocks of both NLS and KG have the four eigenvalues
\begin{equation}\label{ckg:hyper-spec-NLS}
 -\iu\delta^{\sharp}\pm\sigma^{\sharp}, \qquad \iu\delta^{\sharp}\pm\sigma^{\sharp},
\end{equation}
and it holds
\begin{equation}\label{ckg:sigma-def}
\sigma^{\rm NLS}\ge\frac{\sqrt3}{8\pi}\,\varrho \qquad \sigma(\rho;c)\ge\kappa_{\cF}\varrho.
\end{equation}
Thus both the NLS and Klein--Gordon real normal linearizations have nontrivial stable and unstable directions, with exponential rates uniformly bounded from below by a constant times $\varrho$.
\end{lemma}

\begin{proof}
The characteristic polynomial of $L_{\cF}^{(2)}$ gives the two eigenvalues $-\iu\delta\pm\sigma$; the conjugate block gives $\iu\delta\pm\sigma$. The same computation for \(L_{\mathcal F,\rm NLS}^{(2)}\) gives the corresponding NLS eigenvalues, thereby proving \eqref{ckg:hyper-spec-NLS} in both cases. Using the explicit NLS coefficients in Proposition~\ref{ckg:prop-dist-effective} we have
\[
 (\sigma^{\rm NLS})^2 =\frac{16\rho_1\rho_2-(\rho_1+\rho_2)^2}{64\pi^2}.
\]
Writing $\rho_i=\varrho q_i$ with $q_i\in[1/2,3/2]$, we have
\[
 16q_1q_2-(q_1+q_2)^2\ge3,
\]
and hence we obtain the first estimate in \eqref{ckg:sigma-def} . By \eqref{ckg:HF-dist-close},
\[
 |\sigma^2-(\sigma^{\rm NLS})^2|\le Ch\varrho^2
\]
uniformly on $\Xi_0$.  After reducing $h_*$, this implies the second estimate in \eqref{ckg:sigma-def} and completes the proof.
\end{proof}

We do not diagonalize this block in the KAM variables, since the subsequent construction of the homological equations requires a uniform resolvent estimate for the Hamiltonian matrix, not a parameter-dependent choice of eigenvectors.

\subsection{Elimination of the resonant cubic tail}

The polynomials \eqref{ckg:cubic-res-dist-newcoords} and \eqref{ckg:cubic-res-dist-NLS} have zero KAM jet, but \eqref{ckg:Z43-size} and \eqref{ckg:Z43-NLS-size} show that their full rescaled vector-field norm is $O(r_0^{5/3})$, which is larger than the perturbative size $O(r_0^2)$ of $\mathcal{R}_{m,n}$ and $\mathcal{R}_{m,n}^{\rm NLS}$. We therefore remove it by a preliminary homological equation whose only nontrivial inversion is the finite hyperbolic block.

For $k\neq m,n$ we denote
\begin{equation}\label{ckg:Qkc-def}
 Q_k:=\bar a_k b_k, \qquad c_k:=\left( \frac{\sqrt{\rho_2}}{\pi\widetilde w_n\widetilde w_k}, \frac{\sqrt{\rho_1}}{\pi\widetilde w_m\widetilde w_k} \right), \qquad  c_k^{\rm NLS}:= \left(\frac{\sqrt{\rho_2}}{\pi}, \frac{\sqrt{\rho_1}}{\pi}\right).
\end{equation}
Then recalling \eqref{ckg:cubic-res-dist-newcoords} and \eqref{ckg:cubic-res-dist-NLS} it holds
\begin{equation}\label{Z43:newdef}
 Z^{\sharp}_{4,3}=\sum_{k\neq m,n}\Re\bigl[(c_k^{\sharp} U)Q_k\bigr],
\end{equation}
where $U$ is given in \eqref{ckg:UF-def}. Moreover, a direct computation gives
\begin{equation}\label{ckg:Q-commutes}
 \{\Lambda_k^{\sharp}(|a_k|^2+|b_k|^2),Q_k\}=0.
\end{equation}

\begin{lemma}
\label{ckg:lem-exact-finite-resolvent}
There exist $h_*,C_{\cF}>0$, depending only on $m,n$, such that, for every $\rho\in\Xi_0$, $t\in\mathbb R$ and $h \in (0,h_*)$ it holds
\begin{align*}
 \|(L_{\cF,\rm NLS}^{(2)}+\iu tI)^{-1}\| \le\frac{C_{\cF}}{\varrho}, \qquad \|(L_{\cF}^{(2)}+\iu tI)^{-1}\| \le\frac{C_{\cF}}{\varrho}.
\end{align*}
In particular,
\begin{equation}\label{ckg:LF-inv-Lip}
 \|(L_{\cF}^{(2)})^{-1}\|^{\mathcal L(\varrho)} +\|(L_{\cF,\rm NLS}^{(2)})^{-1}\|^{\mathcal L(\varrho)} \le\frac{C_{\cF}}{\varrho}, \quad \|(L_{\cF}^{(2)})^{-1}-(L_{\cF,\rm NLS}^{(2)})^{-1}\|^{\mathcal L(\varrho)} \le\frac{C_{\cF}h}{\varrho}.
\end{equation}
Moreover, for $E$ any $2\times2$ matrix satisfying $\|E\|\le(2C_{\cF})^{-1}\varrho$, then
\begin{equation}\label{ckg:pert-finite-resolvent}
 \|(L_\sharp^{(2)}+E+\iu tI)^{-1}\| \le\frac{2C_{\cF}}{\varrho}, \qquad t\in\mathbb R,
\end{equation}
with the Klein--Gordon estimate understood for $0<h<h_*$.
\end{lemma}

\begin{proof}
For both the matrices, we have
\[
 L_\sharp^{(2)}+\iu tI
 =\iu\begin{pmatrix}t-\alpha^\sharp&-\gamma^\sharp\\[0.2em]
                       \gamma^\sharp&t+\beta^\sharp\end{pmatrix}.
\]
Recalling the definition of $\delta^\sharp$ and $\sigma^{\sharp}$ in \eqref{ckg:sigma-NLS-def}, up to a constant the determinant is
\begin{align*}
 (t-\alpha^\sharp)(t+\beta^\sharp)+(\gamma^\sharp)^2 &=(t-\delta^\sharp)^2
   +(\gamma^\sharp)^2-\left(\frac{\alpha^\sharp+\beta^\sharp}{2}\right)^2\\
 &=(t-\delta^\sharp)^2+(\sigma^\sharp)^2.
\end{align*}
For the KG matrix, Lemma~\ref{ckg:lem-hyperbolic} gives $\sigma\ge\kappa_{\cF}\varrho$, while $|\alpha|+|\beta|+|\gamma|\le C\varrho$.  The adjugate formula yields
\[
 \|(L_{\cF}^{(2)}+\iu tI)^{-1}\| \le C\frac{|t-\delta|+\varrho} {(t-\delta)^2+\sigma^2} \le\frac{C_{\cF}}{\varrho}.
\]
The NLS estimate is obtained similarly by using $\sigma^{\rm NLS}\ge\sqrt3\,\varrho/(8\pi)$.

Since $\operatorname{Lip}_\rho L_{\cF}^{(2)}+\operatorname{Lip}_\rho L_{\cF,\rm NLS}^{(2)} \le C$, the resolvent identity gives
\[
 \operatorname{Lip}_\rho (L_{\cF}^{(2)})^{-1} +\operatorname{Lip}_\rho (L_{\cF,\rm NLS}^{(2)})^{-1} \le C\varrho^{-2},
\]
which proves the first in \eqref{ckg:LF-inv-Lip}.  Moreover,
\[
 (L_{\cF}^{(2)})^{-1}-(L_{\cF,\rm NLS}^{(2)})^{-1} =(L_{\cF}^{(2)})^{-1}(L_{\cF,\rm NLS}^{(2)}-L_{\cF}^{(2)})(L_{\cF,\rm NLS}^{(2)})^{-1}.
\]
Together with \eqref{ckg:HF-dist-close} and the corresponding parameter Lipschitz bounds, this gives the second in \eqref{ckg:LF-inv-Lip}. Finally,
\[
 L_\sharp^{(2)}+E+\iu tI =\bigl(I+E(L_\sharp^{(2)}+\iu tI)^{-1}\bigr)(L_\sharp^{(2)}+\iu tI),
\]
and the assumption on $E$ allows one to invert the first factor by a Neumann series, proving \eqref{ckg:pert-finite-resolvent}.
\end{proof}

\begin{remark}
Recalling \eqref{localmasses}, we use $\Psi$ to denote either of the symplectic inverse coordinate maps $\Psi_m$ and $\Psi_{m,n}$ defined in \eqref{ckg:Psi-same} and \eqref{ckg:Psi-dist}, respectively, and set
\begin{equation*}
 \mathcal N_\pm^\Psi:=\mathcal N_\pm\circ\Psi.
\end{equation*}
The two Hamiltonians $\mathcal N_+^\Psi$ and $\mathcal N_-^\Psi$ Poisson commute. We denote
\[
 \mathcal R_{\phi_+,\phi_-}(\xi,\eta) := \bigl(e^{\iu\phi_+}\xi,e^{\iu\phi_-}\eta\bigr), \qquad (\phi_+,\phi_-)\in\mathbb T^2,
\]
as the independent phase rotations of the two species.  We define the corresponding local gauge action in the transformed coordinates by
\begin{equation*}
 \mathcal G_{\phi_+,\phi_-} := \Psi^{-1}\circ \mathcal R_{\phi_+,\phi_-}\circ \Psi .
\end{equation*}
Since $\Psi$ is symplectic it holds
\begin{equation*}
 \mathcal G_{\phi_+,\phi_-} = \Phi_{\mathcal N_+^\Psi}^{-\phi_+} \circ \Phi_{\mathcal N_-^\Psi}^{-\phi_-},
\end{equation*}
where $\Phi_H^t$ denotes the Hamiltonian flow generated by $H$. 

For an analytic Hamiltonian $F$ we then set
\begin{equation}\label{ckg:PiLG-averaging}
 \Pi_{\rm LG}F := \frac{1}{(2\pi)^2} \int_{\mathbb T^2} F\circ\mathcal G_{\phi_+,\phi_-} \,d\phi_+\,d\phi_-.
\end{equation}
Thus $\Pi_{\rm LG}$ is the coefficientwise projection onto the Fourier--Taylor monomials which commute with both local masses,
\[
 \{\Pi_{\rm LG}F,\mathcal N_+^\Psi\} = \{\Pi_{\rm LG}F,\mathcal N_-^\Psi\} =0.
\]
\end{remark}

\medskip

To perform the preliminary step, for both the cases $m \neq n$ and $m=n$ we denote 
\begin{equation}\label{initial-H2}
 H^{(2)}:=(H\circ\mathcal T_0)\circ\Psi, \qquad H^{(2),\rm NLS}:=(H^{\rm NLS}\circ\mathcal T_0^{\rm NLS}) \circ\Psi.
\end{equation}
We recall that $H, H^{\rm NLS}$ are the Hamiltonians in \eqref{H:xi-eta} and \eqref{HNLS} , $\mathcal{T}_0, \mathcal{T}_0^{\rm NLS}$ are introduced in Theorem \ref{Teorema0}, while $\Psi$ denotes both $\Psi_m$ and $\Psi_{m,n}$ in \eqref{ckg:Psi-same}, \eqref{ckg:Psi-dist}. Moreover we denote
\begin{equation*}
 C(\rho):=\begin{cases}C_{m,m}(\rho;c),&m=n,\\C_{m,n}(\rho;c),&m\ne n,\end{cases}
 \quad C^{\rm NLS}(\rho):=\begin{cases}C_{m,m}^{\rm NLS}(\rho),&m=n,\\C_{m,n}^{\rm NLS}(\rho),&m\ne n,\end{cases}
\end{equation*}
and finally 
\begin{equation}\label{ckg:case-beta-definition}
 \ckgcasebeta:=\begin{cases}0,&m=n,\\1,&m\ne n.\end{cases}
\end{equation}
We have the following
\begin{proposition}
\label{ckg:prop-cubic-elimination}
There exist $h_*,r_*,s_*,C>0$, depending only on $a,p,m,n$, such that, for either choice of tangential sites, for $0<h<h_*$ and $0<r_0<r_*$ the following holds. There exist two families of real analytic symplectic maps, Lipschitz in $\rho \in \Xi_0$,
\begin{align}
 \mathcal C_3(\rho)&:\mathfrak D_q(s_*,3r_0/4) \longrightarrow\mathfrak D_q(s_*,r_0),
 \label{ckg:C3-map}\\
 \mathcal C_3^{\rm NLS}(\rho)&:\mathfrak D_q(s_*,3r_0/4) \longrightarrow\mathfrak D_q(s_*,r_0),
 \label{ckg:C3-NLS-map}
\end{align}
for $q=p,p-4$ such that
\begin{align}
 H^{(2)}\circ\mathcal C_3 &=C(\rho)+N_0+P_0,
  \label{ckg:cubiKGnew}\\
 H^{(2),\rm NLS}\circ\mathcal C_3^{\rm NLS} &=C^{\rm NLS}(\rho)+N_0^{\rm NLS}+P_0^{\rm NLS},
 \label{ckg:cubic-eliminated-H-NLS}
\end{align}
where $N_0$ and $N_0^{\rm NLS}$ are the effective Hamiltonians of both Propositions~\ref{ckg:prop-same-effective} and~\ref{ckg:prop-dist-effective}, the perturbations satisfy
\begin{align}
 \|X_{P_0}\|_{r_0,q,\ckgcasebeta}^{\mathcal L(\varrho), \mathfrak D_q(s_*,3r_0/4)} \le Cr_0^2, \qquad \|X_{P_0^{\rm NLS}}\|_{r_0,q}^{\mathcal L(\varrho), \mathfrak D_q(s_*,3r_0/4)} \le Cr_0^2,
 \label{ckg:cubic-eliminated-size}
\end{align}
and the changes of coordinates satisfy
\begin{align}
 \|\mathcal C_3-\Id\|_{r_0,q}^{\mathcal L(\varrho), \mathfrak D_q(s_*,3r_0/4)} \le Cr_0^{1/3}, \qquad  \|\mathcal C_3^{\rm NLS}-\Id\|_{r_0,q}^{\mathcal L(\varrho), \mathfrak D_q(s_*,3r_0/4)}\le Cr_0^{1/3}.
 \label{ckg:C3-close}
\end{align}
Moreover both maps commute with the momentum, physical-charge and local-gauge actions and 
\begin{align}
 \|\mathcal C_3-\mathcal C_3^{\rm NLS}\|_{r_0,p-4}^{ \mathfrak D_p(s_*,3r_0/4)} \le Chr_0^{1/3}, \qquad  \|X_{\Pi_{\rm LG}P_0-P_0^{\rm NLS}}\|_{r_0,p-4}^{ \mathfrak D_p(s_*,3r_0/4)} \le Chr_0^2.
 \label{ckg:C3-diff}
\end{align}
In the identical-mode case the two maps are the identity, as explained in the proof.
\end{proposition}

\begin{proof}
If $m=n$, we make the choice $\mathcal C_3=\mathcal C_3^{\rm NLS}=\Id$, with $Z_{4,3}=Z_{4,3}^{\rm NLS}=0$. By Theorem~\ref{Teorema0} and Proposition~\ref{ckg:prop-same-effective}, the perturbations are exactly
\begin{equation*}
P_0=\mathcal R_{m,m}+(\widehat P+P_{\rm BNF})\circ\Psi_m,\qquad P_0^{\rm NLS}=\mathcal R_{m,m}^{\rm NLS}+(\widehat P^{\rm NLS}+P_{\rm BNF}^{\rm NLS})\circ\Psi_m.
\end{equation*}
For any Hamiltonian $F$, symplecticity gives
\begin{equation*}
X_{F\circ\Psi_m}(w)=d\Psi_m(w)^{-1}X_F(\Psi_m(w)).
\end{equation*}
The quartic Hamiltonians $\widehat P$ and $\widehat P^{\rm NLS}$ contain only external variables, so their pullbacks have scaled vector-field norm $O(r_0^2)$. The Birkhoff remainders have degree at least six, and the bounds \eqref{res.0.0}--\eqref{res.NLS}, together with $\|\Psi_m(w)\|_{a,q}\lesssim\sqrt\varrho$ and $\varrho=r_0^{4/3}$, give the same $O(r_0^2)$ bound after pullback. Combining these estimates with \eqref{ckg:R-same-size} and \eqref{ckg:R-same-NLS-size} proves \eqref{ckg:cubic-eliminated-size} for $q=p,p-4$. Here $\ckgcasebeta=0$ because the infinite-mode corotation in $\Psi_m$ introduces action corrections which do not preserve the one-weight smoothing from the Birkhoff Normal Form Step, as observed after \eqref{loss:regula}.

For the KG-NLS comparison, recall the decompositions \eqref{comp:PP}. The terms $\widehat P^{\mathcal{NG}}$ and $P_{\rm BNF}^{\mathcal{NG}}$ have vanishing local-gauge projection, whereas the NLS remainders are locally gauge invariant. Since $\Psi_m$ transports the local-gauge action and $\mathcal R_{m,m}$ is locally gauge invariant, we obtain the exact identity
\begin{equation*}
\Pi_{\rm LG}P_0-P_0^{\rm NLS}=\mathcal R_{m,m}-\mathcal R_{m,m}^{\rm NLS}+\Pi_{\rm LG}\bigl[(\widehat P^R+P_{\rm BNF}^R)\circ\Psi_m\bigr],
\end{equation*}
where $P_{\rm BNF}^{\mathcal{NG}}$ is the particular non-LG contribution defined in \eqref{P:0NG}. The purely external quartic structure combined with \eqref{eq:PhatR-est} and the sixth-order Birkhoff estimate combined with the pullback formula yield
\begin{equation*}
\|X_{\widehat P^R\circ\Psi_m}\|_{r_0,p-4}+\|X_{P_{\rm BNF}^R\circ\Psi_m}\|_{r_0,p-4}\lesssim hr_0^2.
\end{equation*}
Together with \eqref{ckg:R-same-R-close}, these estimates prove the second bound in \eqref{ckg:C3-diff}. All bounds involving the identity maps are immediate.

\medskip

For the rest of the proof we study the case $m\ne n$. By Theorem~\ref{Teorema0} and Proposition~\ref{ckg:prop-dist-effective}, the two Hamiltonians \eqref{initial-H2} can be written as
\begin{align}
 H^{(2)}&=C(\rho)+N_0+Z_{4,3}+R,
 \label{ckg:H2-before-C3}\\
 H^{(2),\rm NLS}&=C^{\rm NLS}(\rho)+N_0^{\rm NLS} +Z_{4,3}^{\rm NLS}+R^{\rm NLS},
 \label{ckg:H2-before-C3-NLS}
\end{align}
where
\begin{align}
 \|X_R\|_{r_0,q,1}^{\mathcal L(\varrho), \mathfrak D_q(s_*,r_0)}\le Cr_0^2, \qquad  \|X_{R^{\rm NLS}}\|_{r_0,q}^{\mathcal L(\varrho), \mathfrak D_q(s_*,r_0)}\le Cr_0^2
 \label{ckg:R-before-C3-size}
\end{align}
for both $q=p,p-4$. The estimates \eqref{ckg:R-before-C3-size} follow from the same weighted-degree counting as in the identical-mode case, using Theorem~\ref{Teorema0} and Proposition~\ref{ckg:prop-dist-effective}. Unlike $\Psi_m$, the map $\Psi_{m,n}$ rotates only the finite-dimensional semi-external block and leaves all generic normal variables unchanged. Consequently, the pullback preserves the one-weight smoothing of the Klein-Gordon remainders, yielding the stronger estimate with $\beta=1$.

We construct now the change of coordinates. Recalling \eqref{ckg:Qkc-def} and \eqref{Z43:newdef}, we introduce 
\begin{align}\label{ckg:S3-def}
 s_k&:=c_k(L_{\cF}^{(2)})^{-1},& \qquad &S_3:=\sum_{k\neq m,n}\Re\bigl[(s_kU)Q_k\bigr],& \\
\label{c:kdef}
 s_k^{\rm NLS}&:=c_k^{\rm NLS}(L_{\cF,\rm NLS}^{(2)})^{-1},& \quad &S_3^{\rm NLS}:=\sum_{k\neq m,n} \Re\bigl[(s_k^{\rm NLS}U)Q_k\bigr].&
\end{align}
Since $\dot U=L_{\cF}^{(2)}U$, antisymmetry of the Poisson bracket gives $\{H_\cF,U\}=-L_{\cF}^{(2)}U$ and
\[
 \{H_\cF,s_kU\}=-(s_kL_{\cF}^{(2)})U=-c_kU.
\]
The same holds for the NLS. Recalling that by \eqref{ckg:H2-before-C3}, \eqref{ckg:H2-before-C3-NLS}, \eqref{ckg:dist-NLS-decomposition} and \eqref{ckg:dist-decomposition}, combined with \eqref{ckg:Q-commutes} and \eqref{Z43:newdef}, we have the exact homological equation
\begin{equation}\label{ckg:S3-homological}
 \{N_0,S_3\}+Z_{4,3}=0, \qquad \{N_0^{\rm NLS},S_3^{\rm NLS}\}+Z_{4,3}^{\rm NLS}=0,
\end{equation}
The monomials in $S_3$ satisfy
\begin{equation}\label{ckg:S3-symmetries}
 \{S_3,\mathcal M\}=\{S_3,\mathcal Q\} =\{S_3,\mathcal N\}=0,
\end{equation}
and the same for $S_3^{\rm NLS}$. We next estimate the generator, including its parameter dependence.  From \eqref{ckg:Qkc-def}, \eqref{ckg:LF-inv-Lip} and \eqref{c:kdef}
\begin{align}\label{ckg:sk-bound}
 \widetilde w_k|c_k|^{\mathcal L(\varrho)} + |c_k^{\rm NLS}|^{\mathcal L(\varrho)} \le C\sqrt\varrho, \qquad 
 \widetilde w_k|s_k|^{\mathcal L(\varrho)}+ |s_k^{\rm NLS}|^{\mathcal L(\varrho)} \le\frac{C}{\sqrt\varrho}
\end{align}
A direct differentiation of \eqref{ckg:S3-def} and \eqref{c:kdef} gives
\begin{align}
 \|X_{S_3}\|_{r_0,q,1}^{\mathcal L(\varrho), \mathfrak D_q(s_*,r_0)} \le C\frac{r_0}{\sqrt\varrho} =Cr_0^{1/3}, \qquad  \|X_{S_3^{\rm NLS}}\|_{r_0,q}^{\mathcal L(\varrho), \mathfrak D_q(s_*,r_0)} \le C\frac{r_0}{\sqrt\varrho} =Cr_0^{1/3},
 \label{ckg:S3-size}
\end{align}
for both $q=p,p-4$. Indeed, for the Klein--Gordon infinite-dimensional components the factor $\widetilde w_k^{-1}$ in \eqref{ckg:sk-bound} is exactly absorbed by the output $\beta=1$ weight, while the finite components are controlled by Cauchy--Schwarz in the $k$-sum.  The NLS estimate is the same computation without the smoothing weight.  Since both generators are independent of $(\theta,y)$, their flows fix $\theta$ and $y$. The estimates \eqref{ckg:S3-size} are taken on the larger domain $\mathfrak D_q(s_*,r_0)$. We introduce $\mathcal C_3:=\Phi_{S_3}^1$, $\mathcal C_3^{\rm NLS}:=\Phi_{S_3^{\rm NLS}}^1$ that by standard flow estimates are therefore well defined from $\mathfrak D_q(s_*,3r_0/4)$ to $\mathfrak D_q(s_*,r_0)$ and satisfy the first in \eqref{ckg:C3-close}.  They are symplectic and, by \eqref{ckg:S3-symmetries}, preserve the stated symmetries.

Using the second-order Lie formula and \eqref{ckg:S3-homological},
\begin{align*}
 H^{(2)}\circ\Phi_{S_3}^1-N_0-C(\rho)
 ={}&R+\{Z_{4,3},S_3\}+\{R,S_3\}\\
 &+\int_0^1(1-t) \{\{H^{(2)},S_3\},S_3\}\circ\Phi_{S_3}^t\,dt.
\end{align*}
By \eqref{ckg:R-before-C3-size} $X_R$ is of order $r_0^2$, while the scaled Poisson-bracket estimate, using the vector-field formula \eqref{ckg:XF-KAM}, gives
\[
 \|X_{\{Z_{4,3},S_3\}}\|_{r_0,q,1}^{\mathfrak D_q(s_*,3r_0/4)} \lesssim r_0^{5/3}r_0^{1/3}=r_0^2.
\]
The above estimate is obtained by combining \eqref{ckg:S3-size}, \eqref{ckg:Z43-size}, the second identity in \eqref{stime}, applied to the Hamiltonians present here, and standard Cauchy estimates (see \cite[Lemmas~A.3-A.4]{poschel1996kam}). On the other hand, the third term is $O(r_0^{7/3})$. In the integral, the leading contribution is $\{Z_{4,3},S_3\}$, of order $r_0^2$, while all other contributions are of higher order. Hence the entire integral is $O(r_0^2)$, proving the first estimate in \eqref{ckg:cubic-eliminated-size}. The identical argument, using the NLS homological equation in \eqref{ckg:S3-homological}, proves \eqref{ckg:cubic-eliminated-H-NLS} and the second in \eqref{ckg:cubic-eliminated-size}.

For the Klein--Gordon--NLS comparison, since $m,n$ are fixed,
\begin{align*}
 \|L_{\cF}^{(2)}-L_{\cF,\rm NLS}^{(2)}\|^{\mathcal L(\varrho)}\le Ch\varrho, \qquad  |c_k-c_k^{\rm NLS}|^{\mathcal L(\varrho)} \le Ch\sqrt\varrho\,\langle k\rangle^2.
\end{align*}
By the second in \eqref{ckg:LF-inv-Lip}, we therefore obtain
\begin{equation}\label{ckg:sk-diff}
 |s_k-s_k^{\rm NLS}|^{\mathcal L(\varrho)} \le\frac{Ch\langle k\rangle^2}{\sqrt\varrho} \implies  \|X_{S_3-S_3^{\rm NLS}}\|_{r_0,p-4}^{\mathfrak D_p(s_*,r_0)} \le Chr_0^{1/3}.
\end{equation}
The flow comparison in \eqref{ckg:C3-diff} then follows from the standard comparison estimates for analytic flows in a weaker Banach topology (see \cite[Appendix~B]{BBG}), applied with the strong and weak normal spaces $\ell^{a,p}$ and $\ell^{a,p-4}$, respectively.

For the comparison of the transformed remainders, we first note that the estimates above, together with Theorem~\ref{Teorema0} and the explicit coefficient comparison in Proposition~\ref{ckg:prop-dist-effective}, give
\begin{equation}\label{ckg:R-preC3-diff}
 \|X_{\Pi_{\rm LG}R-R^{\rm NLS}}\|_{r_0,p-4}^{\mathfrak D_p(s_*,3r_0/4)} \le Ch r_0^2, \qquad \|X_{Z_{4,3}-Z_{4,3}^{\rm NLS}}\|_{r_0,p-4}^{\mathfrak D_p(s_*,3r_0/4)} \le Ch r_0^{5/3}.
\end{equation}
Here the same angle-elimination map $\Psi_{m,n}$ is used in the Klein--Gordon and NLS systems, so it produces no additional Klein--Gordon--NLS difference at fixed parameter.

Since $S_3$ is locally gauge invariant, its flow commutes with the local gauge action. Consequently
\[
 \Pi_{\rm LG}\{F,S_3\} = \{\Pi_{\rm LG}F,S_3\}, \qquad \Pi_{\rm LG}(F\circ\Phi_{S_3}^t) = (\Pi_{\rm LG}F)\circ\Phi_{S_3}^t.
\]
Using $\{N_0,S_3\}=-Z_{4,3}$, the transformed perturbation can be written as
\[
\begin{aligned}
 P_0={}&R+\{Z_{4,3},S_3\}+\{R,S_3\}\\
 &+\int_0^1(1-t) \bigl\{ -Z_{4,3}+\{Z_{4,3},S_3\}+\{R,S_3\}, S_3 \bigr\}\circ\Phi_{S_3}^t\,dt ,
\end{aligned}
\]
and the same identity holds for the NLS quantities.

To prove \eqref{ckg:C3-diff} we apply $\Pi_{\rm LG}$ to the expression of $P_0$ written above. Since both cubic generators are locally gauge invariant, the projection commutes with their Poisson brackets and flows. Moreover, recalling that $Z_{4,3}$ and $S_3$ are locally gauge invariant, the resulting differences contain at least one of
\[
\Pi_{\rm LG}R-R^{\rm NLS},\qquad
Z_{4,3}-Z_{4,3}^{\rm NLS},\qquad
S_3-S_3^{\rm NLS},
\qquad
\Phi_{S_3}^t-\Phi_{S_3^{\rm NLS}}^t.
\]
The bounds \eqref{ckg:R-preC3-diff}, \eqref{ckg:S3-size}, \eqref{ckg:sk-diff} and the first in \eqref{ckg:C3-diff}, combined with the standard commutator and flow-comparison estimates, therefore yield
\[
\|X_{\Pi_{\rm LG}P_0-P_0^{\rm NLS}}\|_{r_0,p-4}^{\mathfrak D_p(s_*,3r_0/4)}
\le Chr_0^2,
\]
applying the same analytic estimates used in Theorem \ref{Teorema0} and in \cite[Section~3]{BBG}. This proves the second in \eqref{ckg:C3-diff}.
\end{proof}

\section{Block normal forms and the KAM phase space}
\label{sec:block-KAM-setting}

The transformations of Theorem~\ref{Teorema0}, Propositions~\ref{ckg:prop-same-effective}--\ref{ckg:prop-dist-effective} and Proposition~\ref{ckg:prop-cubic-elimination} reduce the Hamiltonian to a perturbation of a quadratic normal form which is no longer scalar on the normal variables.  The multiplicity $\lambda_j=\lambda_{-j}$ and the two species are retained in uniformly finite Hermitian blocks.  In the distinct-mode case one also has the finite hyperbolic block produced by the semi-external modes. The following KAM block approach is based on standard known results, see for instance \cite{EGK} and \cite{GVdR}.

\subsection{Elliptic clusters, block norms and normal projection}

Let $\mathcal E$ be the elliptic normal index set introduced in \eqref{ckg:KAM-variables}, with the sequence weights of \eqref{eq:norm-apb}.  For $j\in\mathbb Z$ define
\begin{equation}\label{ckg:block-partition}
 [j]:=\{(\ell,s)\in\mathcal E:\ |\ell|=|j|\},
\end{equation}
and the resulting family of non-empty blocks is denoted by $\mathfrak B$.  Since there are at most two signs of the Fourier index and two species,
\begin{equation}\label{ckg:block-size}
 \forall \,  B\in\mathfrak B, \qquad d_B:=\#B\le4.
\end{equation}
If $B=[j]$, we set
\begin{equation*}
 \langle B\rangle:=\langle j\rangle, \qquad \widetilde w_B:=(1+h j^2)^{1/2}.
\end{equation*}

Let $B, C\in\mathfrak B$, for a matrix $M:\mathbb C^C\to\mathbb C^B$ we use the operator norm $\|M\|$, and, when diagonalizing a finite Sylvester operator, the Hilbert--Schmidt norm $\|M\|_{\rm HS}$.  The associated inner product is
$\langle X,Y\rangle_{\rm HS}:=\operatorname{Tr}(X^*Y)$.  By \eqref{ckg:block-size},
\begin{equation}\label{ckg:finite-norm-equivalence}
 \|M\|\le \|M\|_{\rm HS} \le 2\,\|M\|,
\end{equation}
uniformly in the blocks.  For linear operators $\mathscr L$ acting on matrix spaces, $\|\mathscr L\|_{\mathrm{op}\to\mathrm{op}}$ and $\|\mathscr L\|_{\mathrm{HS}\to\mathrm{HS}}$ denote the induced norms with respect to the operator and Hilbert--Schmidt matrix norms, respectively.When no subscript is displayed, $\|\mathscr L\|$ denotes the operator-to-operator norm

For a parameter-dependent block family $D=(D_B)_{B\in\mathfrak B}$ on $\Xi \subseteq \Xi_0$ and $\beta\ge0$ we put
\begin{align*}
 |D|_{\infty,\beta} &:={\sup}_{\rho\in\Xi}\,{\sup}_{B\in\mathfrak B} \widetilde w_B^\beta\|D_B(\rho)\|, \quad 
 |D|_{\infty,\beta}^{\lip} :={\sup}_{\substack{\rho,\rho'\in\Xi\\ \rho\ne\rho'}} \frac{{\sup}_{B\in\mathfrak B} \widetilde w_B^\beta \|D_B(\rho)-D_B(\rho')\|} {|\rho-\rho'|}, \\
|D|_{\infty,\beta}^{\mathcal L(\alpha)} &:=|D|_{\infty,\beta} +\alpha |D|_{\infty,\beta}^{\lip}, \qquad \qquad  \quad  \, \,\forall \,  \alpha>0.
\end{align*}
For finite-dimensional vectors and matrices we use the convention \eqref{ckg:finite-Lip-convention}. When $\beta=0$ the index is dropped.

In the distinct-mode case let $d_\cF:=\#\mathcal F=2$ and write, as in \eqref{ckg:UF-def},
\begin{equation*}
W_{\cF}=(z_\cF,\bar z_\cF)\in\mathbb C^{2d_\cF},\qquad \mathcal C_{\cF}:=\begin{pmatrix}0&I_{d_\cF}\\ I_{d_\cF}&0\end{pmatrix},\qquad \mathcal J_\cF:=\begin{pmatrix}0&I_{d_\cF}\\-I_{d_\cF}&0\end{pmatrix}.
\end{equation*}
The real subspace is characterized by $W_{\cF}=\mathcal C_{\cF}\overline{W_{\cF}}$. In the identical-mode case $\mathcal F=\varnothing$, and all expressions involving the finite block are omitted. Here $z_B$ and $W_\cF$ are regarded as column vectors.

\begin{definition}\label{def:blocknormalform}
We say that a real Hamiltonian is in \emph{block normal form} if
\begin{equation}\label{ckg:block-N}
N=\langle\omega(\rho),y\rangle+\sum_{B\in\mathfrak B} \langle \bar z_B, A_B(\rho)z_B\rangle+\frac12\langle W_{\cF},K_\cF(\rho)W_{\cF}\rangle,
\end{equation}
where $\omega(\rho)\in\mathbb R^2$, each elliptic matrix $A_B(\rho)$ is Hermitian, $A_B=A_B^*$, and, when $\mathcal F\ne\varnothing$, the finite matrix $K_\cF(\rho)\in\mathbb C^{2d_\cF\times2d_\cF}$ satisfies
\begin{equation*}
K_\cF^{t}=K_\cF,\qquad K_\cF=\mathcal C_{\cF}\overline{K_\cF}\mathcal C_{\cF}.
\end{equation*}
The coefficients are required to satisfy the exact symmetry relations
\begin{equation}\label{ckg:block-exact-symmetries}
\{N,\mathcal M\}=0,\qquad \{N,\mathcal Q\}=0.
\end{equation}
\end{definition}

\medskip

For $\mathcal F\ne\varnothing$, the finite quadratic part of $N$ generates the linear Hamiltonian system $\dot W_{\cF}=\mathbb L_{\cF}W_{\cF}$, where, by \eqref{ckg:KAM-symplectic-form} and \eqref{ckg:XF-KAM},
\begin{equation}\label{ckg:LF-from-KF}
\mathbb L_{\cF}:=-\iu\mathcal J_\cF K_\cF.
\end{equation}
This $2d_\cF\times2d_\cF$ matrix is distinct from the $2\times2$ matrix $L_{\cF}^{(2)}$ in \eqref{ckg:LF-def}, which acts only on $U$ from \eqref{ckg:UF-def}.
For the Klein--Gordon normal form, we separate the rest energy (the term of order $c^2$) as
\begin{equation}\label{ckg:block-fast-split}
\omega=c^2(1,1)+\widehat\omega,\qquad A_B=q_Bc^2I_B+\widetilde A_B,\qquad q_B\in\{0,1\}.
\end{equation}
By \eqref{ckg:Heff-same} and \eqref{ckg:Heff-dist}, for the identical-mode torus $q_B=0$ for every $B$. For the distinct-mode torus $q_B=1$ on every generic elliptic block. The semi-external modes are contained in $\mathcal F$ and, after the angle elimination \eqref{ckg:distinct-angle-change}, carry no uncancelled rest energy. The corresponding initial NLS normal forms are obtained from \eqref{ckg:same-NLS-decomposition} and \eqref{ckg:dist-NLS-decomposition} by removing the explicit $c^2$ terms and setting $h=0$ in the shifted coefficients.

It is convenient to introduce the corresponding fast generator
\begin{equation*}
 \mathcal F_*:=y_1+y_2+ \sum_{B\in\mathfrak B}q_B \sum_{a\in B}z_a\bar z_a.
\end{equation*}
By \eqref{ckg:sym-same-N} and \eqref{ckg:sym-dist-N}, in the two coordinate systems considered above we have
\begin{equation}\label{ckg:Fstar-mass}
 \mathcal F_*=\mathcal N-(\rho_1+\rho_2).
\end{equation}
Moreover every block normal form satisfies
\begin{equation}\label{ckg:block-fast-commute}
 \{N,\mathcal F_*\}=0, \qquad N=c^2\mathcal F_*+N_{\rm slow},
\end{equation}
where $N_{\rm slow}$ is obtained from \eqref{ckg:block-N} by replacing $\omega$ by $\widehat\omega$ and $A_B$ by $\widetilde A_B$. Indeed, each retained elliptic quadratic coefficient is of $z_B\bar z_B$ type and $q_B$ is constant on $B$, while the finite block does not occur in $\mathcal F_*$. 

Let
\begin{equation}\label{mony}
M_{k,\tt{m}, \varpi, \varkappa, \gamma}:=e^{\iu k\cdot\theta}y^{\tt{m}} z^\varpi\bar z^\varkappa W_{\cF}^\gamma
\end{equation}
be a Fourier--Taylor monomial, and put $\ell:=\varpi-\varkappa$ on the elliptic index set.  Its \emph{fast defect} is
\begin{equation}\label{ckg:fast-defect}
 \mathfrak L(k,\ell):= k_1+k_2+ \sum_{B\in\mathfrak B}q_B\sum_{a\in B}\ell_a.
\end{equation}
For a physically gauge invariant monomial, local gauge invariance is equivalent to
\begin{equation*}
 \mathfrak L(k,\ell)=0.
\end{equation*}
Indeed, from \eqref{ckg:Fstar-mass}, under the physical gauge constraint, local gauge invariance is equivalent to invariance under the total mass. For a physically gauge invariant Fourier--Taylor series, the projection \eqref{ckg:PiLG-averaging} is equivalently
\begin{equation*}
 \Pi_{\rm LG}P:= \sum_{\mathfrak L(k,\ell)=0}P_{k,\tt{m},\varpi,\varkappa,\gamma} e^{\iu k\cdot\theta}y^{\tt{m}} z^\varpi\bar z^\varkappa W_{\cF}^\gamma.
\end{equation*}
Since \eqref{ckg:block-fast-commute} holds, $\operatorname{ad}_N$ preserves every fast-defect eigenspace and hence commutes with $\Pi_{\rm LG}$ on the class of Hamiltonians under consideration.

\begin{definition}\label{ckg:def-normal-projection}
Let $P^T$ be a real, translation invariant and physically gauge invariant KAM jet (see Definition \ref{ckg:def-KAM-jet-preliminary}). Its projection $[P^T]$ is the sum of the following coefficients:
\begin{enumerate}[label=(\roman*)]
 \item the zero Fourier modes of the scalar and $y$-linear parts;
  \item for each $B\in\mathfrak B$, the complete zero-Fourier Hermitian quadratic form
\[
\langle \bar z_B,D_Bz_B\rangle, \qquad D_B=D_B^*,
\]
which is compatible with momentum and physical charge.
\item in the distinct-mode case, the complete zero-Fourier real quadratic form
\[
\frac12 \langle W_{\cF},E_{\cF}W_{\cF}\rangle, \qquad E_{\cF}^{t}=E_{\cF}, \qquad E_{\cF}=\mathcal C_{\cF}\overline{E_{\cF}}\mathcal C_{\cF},
\]
which is compatible with momentum and physical charge.
\end{enumerate}
Here \(D_B\) and \(E_{\mathcal F}\) denote the matrices of the corresponding zero-Fourier quadratic coefficients of \(P^T\).
The scalar average is dynamically irrelevant and is absorbed into the parameter-dependent constant of the Hamiltonian.  The normalization $[S]=0$ means that the same coefficients of the generating function vanish.
\end{definition}

\medskip

The projection is bounded uniformly in all the norms used below.  More precisely, for $q>1/2$, $\beta \in [0,1]$, $\alpha>0$ and every parameter set $\Xi$,
\begin{equation}\label{ckg:normal-projection-bounded}
 \|X_{[P^T]}\|_{r,q,\beta}^{\mathcal L(\alpha), \mathfrak D_q(s,r)} \le C \|X_{P^T}\|_{r,q,\beta}^{\mathcal L(\alpha), \mathfrak D_q(s,r)}.
\end{equation}
Indeed, $[\cdot]$ is the zero-Fourier projection followed by orthogonal projections on spaces of dimension bounded by \eqref{ckg:block-size}; the finite block has fixed dimension.  Moreover every coefficient retained by $[\cdot]$ has zero fast defect.  Hence
\begin{equation}\label{ckg:normal-projection-LG}
 \Pi_{\rm LG}[P^T]=[P^T], \qquad [(I-\Pi_{\rm LG})P^T]=0.
\end{equation}

\subsection{Exact selection rules}

The exact symmetries substantially reduce the Fourier support of a KAM jet.
\begin{lemma}[Identical-mode selection]\label{ckg:lem-selection-same}
Let a monomial \eqref{mony} satisfy momentum conservation and physical gauge invariance in the identical-mode coordinates.  If the monomial is locally gauge invariant, then $k=0$.  If it is not locally gauge invariant, then
\begin{equation}\label{ckg:same-k-diag}
 k_1=k_2, \qquad \mathfrak L(k,\ell)=2k_1\ne0.
\end{equation}
\end{lemma}

\begin{proof}
By \eqref{ckg:sym-same-N}, local gauge invariance gives $k_1+k_2=0$, because the normal variables do not occur in the total mass in these coordinates.  Physical gauge invariance and \eqref{ckg:sym-same-Q} give $k_1-k_2=0$.  Thus $k=0$ in the local sector. Without local gauge invariance the physical-charge relation still gives $k_1=k_2$.  Since $q_B=0$ for every identical-mode elliptic block, \eqref{ckg:fast-defect} becomes $\mathfrak L(k,\ell)=k_1+k_2=2k_1$, which is nonzero precisely in the non-local sector.  This proves \eqref{ckg:same-k-diag}.
\end{proof}

In the distinct-mode coordinates, the invariance conditions can be expressed directly in terms of the Fourier--Taylor exponents of a monomial \eqref{mony}. By \eqref{ckg:sym-dist-N}--\eqref{ckg:sym-dist-M}, invariance under the total mass $\mathcal N$, the physical charge $\mathcal Q$, and the momentum $\mathcal M$ is equivalent, respectively, to
\begin{align}
k_1+k_2+\sum_{j\ne m,n}(\ell_{j,+}+\ell_{j,-})&=0, \label{ckg:sel-dist-N}\\
k_1-k_2+\sum_{j\ne m,n}(\ell_{j,+}-\ell_{j,-})+2\ell_{a_n}-2\ell_{b_m}&=0, \notag\\
mk_1+nk_2+\sum_{j\ne m,n}j(\ell_{j,+}+\ell_{j,-})&=0, \label{ckg:sel-dist-M}
\end{align}
where $\ell=\varpi-\varkappa$ and $\ell_{a_n}$,$\ell_{b_m}$ the signed exponents of the two semi-external variables. In particular, for a physically gauge invariant monomial, local gauge invariance is equivalent to the first of these conditions.

\begin{lemma}[Distinct-mode selection]\label{ckg:lem-selection-dist}
Let a monomial \eqref{mony} be translation invariant and locally gauge invariant in the distinct-mode coordinates, and let $d$ denote its total normal degree, then $|k|_1\le 2d.$ In particular, for a KAM-jet monomial, whose total normal degree is at most two,
\begin{equation}\label{ckg:k-local-bound}
|k|_1\le4.
\end{equation}
Moreover, if no elliptic normal variables occur then $k=0$.
\end{lemma}

\begin{proof}
Set
\begin{equation*}
S:=\sum_{j\ne m,n}(\ell_{j,+}+\ell_{j,-}),\qquad T:=\sum_{j\ne m,n}(\ell_{j,+}-\ell_{j,-})+2\ell_{a_n}-2\ell_{b_m}.
\end{equation*}
Let $g_{\cE}:=|\varpi|+|\varkappa|$ and $g_{\cF}:=|\gamma|$, so that $d=g_{\cE}+g_{\cF}$. Then $|S|\le g_{\cE}$ and $|T|\le g_{\cE}+2g_{\cF}$. Since local gauge invariance gives both the total-mass and physical-charge relations, we have $k_1+k_2=-S$ and $k_1-k_2=-T$. Consequently,
\begin{equation*}
|k|_1=\frac12\bigl(|S+T|+|S-T|\bigr)=\max\{|S|,|T|\}\le g_{\cE}+2g_{\cF}\le2d.
\end{equation*}
This proves \eqref{ckg:k-local-bound} when $d\le2$. Finally, if $\varpi=\varkappa=0$, then $S=0$, and the total-mass and momentum relations reduce to $k_1+k_2=0$ and $mk_1+nk_2=0$. Since $m\ne n$, it follows that $k=0$.
\end{proof}

\subsection{The initial KAM data}

We recall that the maps $\mathcal C_3$ \eqref{ckg:C3-map} introduced in Proposition \ref{ckg:prop-cubic-elimination} to cancel the cubic tail, are defined on $\mathfrak D_q(s_*,3r_0/4)$.  We therefore start the KAM scheme on the slightly smaller radius
\begin{equation}\label{ckg:initial-KAM-radius}
 s_{\rm K}:=s_*, \qquad r_{\rm K}:=\frac34r_0,
\end{equation}
and we notice that $r_{\rm K}\asymp r_0$ and $\varrho=r_0^{4/3}\asymp r_{\rm K}^{4/3}$. Define the fully transformed Hamiltonians by
\begin{equation}\label{ckg:H0-def}
 H_0:=((H\circ\mathcal T_0)\circ\Psi)\circ\mathcal C_3, \qquad H_0^{\rm NLS}:= ((H^{\rm NLS}\circ\mathcal T_0^{\rm NLS})\circ\Psi) \circ\mathcal C_3^{\rm NLS},
\end{equation}
where we recall that in the identical-mode case the convention $\mathcal C_3=\mathcal C_3^{\rm NLS}=\mathrm{Id}$, while $\Psi$ is either $\Psi_m$ or $\Psi_{m,n}$, defined in \eqref{ckg:Psi-same} and \eqref{ckg:Psi-dist}, respectively.
For the effective quadratic Hamiltonians \eqref{ckg:Heff-same} and \eqref{ckg:Heff-dist}, denote by $\omega_0$ and $A_{0,B}$ their tangential frequencies and elliptic block matrices in the representation \eqref{ckg:block-N}. We use $\omega_0^{\rm NLS}$ and $A_{0,B}^{\rm NLS}$ for the corresponding coefficients of \eqref{ckg:Heff-same-NLS} and \eqref{ckg:Heff-dist-NLS}. With the splitting \eqref{ckg:block-fast-split} and the finite generator convention \eqref{ckg:LF-from-KF}, set
\begin{equation}\label{ckg:initial-shifted-coefficients}
\begin{aligned}
 \widehat\omega_0^{\rm KG}&:=\omega_0-c^2(1,1),
 &\widetilde A_{0,B}^{\rm KG}&:=A_{0,B}-q_Bc^2I_B,\\
 K_{\cF,0}^{\rm KG}
 &:=\begin{pmatrix}
 0&\gamma&\alpha&0\\
 \gamma&0&0&\beta\\
 \alpha&0&0&\gamma\\
 0&\beta&\gamma&0
 \end{pmatrix},
 &K_{\cF,0}^{\rm NLS}
 &:=\begin{pmatrix}
 0&\gamma^{\rm NLS}&\alpha^{\rm NLS}&0\\
 \gamma^{\rm NLS}&0&0&\beta^{\rm NLS}\\
 \alpha^{\rm NLS}&0&0&\gamma^{\rm NLS}\\
 0&\beta^{\rm NLS}&\gamma^{\rm NLS}&0
 \end{pmatrix},\\
 \mathbb L_{\cF,0}^{\rm KG}
 &:=-\iu\mathcal J_\cF K_{\cF,0}^{\rm KG},
 &\mathbb L_{\cF,0}^{\rm NLS}
 &:=-\iu\mathcal J_\cF K_{\cF,0}^{\rm NLS},
\end{aligned}
\end{equation}
Here the finite-block coefficients are those of \eqref{ckg:abg-dist} and \eqref{ckg:HFNLS-dist}, respectively. When $\cF=\varnothing$ the quantities in the last two lines of \eqref{ckg:initial-shifted-coefficients} are not present.
\begin{theorem}[Initial block normal form]
\label{ckg:thm-initial-KAM-data}
There exist $h_*,r_*,C>0$, depending only on $a,p$ and the fixed tangential sites and independent of $c,r_0,\rho$, such that, for $0<h<h_*$ and $0<r_0<r_*$, the Hamiltonians \eqref{ckg:H0-def} are real analytic on $\mathfrak D_q(s_{\rm K},r_{\rm K})$, $q=p,p-4$, and can be written as
\begin{align}
 H_0&=C_0(\rho)+N_0+P_0,
 \label{ckg:H0-block}\\
 H_0^{\rm NLS}&=C_0^{\rm NLS}(\rho)+N_0^{\rm NLS}+P_0^{\rm NLS}.
 \label{ckg:H0-block-NLS}
\end{align}
Here, recalling the constants in \eqref{ckg:C-same} and \eqref{ckg:C-dist} and their NLS counterparts in \eqref{ckg:same-NLS-decomposition} and \eqref{ckg:dist-NLS-decomposition},
\begin{equation}\label{ckg:C0-piecewise}
 C_0=
 \begin{cases}
 C_{m,m}(\rho;c),&m=n,\\
 C_{m,n}(\rho;c),&m\ne n,
 \end{cases}
 \qquad C_0^{\rm NLS}=
 \begin{cases}
 C_{m,m}^{\rm NLS}(\rho),&m=n,\\
 C_{m,n}^{\rm NLS}(\rho),&m\ne n,
 \end{cases}
\end{equation}
and $N_0$, $N_0^{\rm NLS}$ are in block normal form as in Definition \ref{def:blocknormalform}, being respectively the quadratic Hamiltonians \eqref{ckg:Heff-same} and \eqref{ckg:Heff-same-NLS} for $m=n$, and \eqref{ckg:Heff-dist} and \eqref{ckg:Heff-dist-NLS} when $m \neq n$.

If we define the KG--NLS coefficient differences by
\begin{equation}\label{ckg:N0-decomp}
\begin{aligned}
 \widehat\omega_0^{\rm KG}-\omega_0^{\rm NLS}&=:\omega_0^R,\quad \widetilde A_{0,B}^{\rm KG}-A_{0,B}^{\rm NLS}&=:A_{0,B}^R,\quad K_{\cF,0}^{\rm KG}-K_{\cF,0}^{\rm NLS}&=:K_{\cF,0}^R.
\end{aligned}
\end{equation}
then it holds
\begin{equation}\label{ckg:N0R-est}
 |\omega_0^R|^{\mathcal L(\varrho)} +\sup_{B\in\mathfrak B}\langle B\rangle^{-4} \|A_{0,B}^R\|^{\mathcal L(\varrho)} +\varrho^{-1}\|K_{\cF,0}^R\|^{\mathcal L(\varrho)} \le Ch.
\end{equation}

The perturbations $P_0$ and $P_0^{\rm NLS}$ are real analytic, translation invariant and physically gauge invariant.  The following decomposition hold
\begin{equation}\label{ckg:P0-decomp}
 P_0=P_0^{\rm NLS}+P_0^{\mathcal{NG}}+P_0^R, \qquad P_0^{\mathcal{NG}}:=(I-\Pi_{\rm LG})P_0, \qquad P_0^R:=\Pi_{\rm LG}P_0-P_0^{\rm NLS}.
\end{equation}
In particular $P_0^{\rm NLS}$ and $P_0^R$ are locally gauge invariant and $\Pi_{\rm LG}P_0^{\mathcal{NG}}=0$.  
Moreover there exists a constant $C_{\rm in}>0$, independent of $h$ and $r_0$, such that by setting $\epsilon_0:=C_{\rm in}r_{\rm K}^2$ we have
\begin{align}
\max_{q\in\{p,p-4\}} \|X_{P_0}\|_{r_{\rm K},q,\ckgcasebeta}^{\mathcal L(\varrho), \mathfrak D_q(s_{\rm K},r_{\rm K})} &\le \epsilon_0, \label{ckg:E0-full}\\
\epsilon_0^{\rm NLS}:= \max_{q\in\{p,p-4\}} \|X_{P_0^{\rm NLS}}\|_{r_{\rm K},q}^{\mathcal L(\varrho),\mathfrak D_q(s_{\rm K},r_{\rm K})} &\le \epsilon_0,
\label{ckg:E0-NLS}\\
\|X_{P_0^R}\|_{r_{\rm K},p-4}^{\mathfrak D_p(s_{\rm K},r_{\rm K})}&\le h\epsilon_0.
\label{ckg:P0R-jet}
\end{align}

Finally, the shifted tangential frequency maps are lipeomorphisms on $\Xi_0$, with constants uniform in $c$.
\end{theorem}

The theorem collects the preceding normal-form results in the standard KAM format (see for instance Theorem 3.10 in \cite{BBG} or Section 7 in \cite{kuksin1996invariant}). The additional points specific to the present setting are the finite hyperbolic block in the distinct-mode case and the uniform KG–NLS comparison of the initial normal forms and perturbations.

\begin{proof}
The decompositions and perturbation bounds follow from Proposition~\ref{ckg:prop-cubic-elimination}. Since $r_{\rm K}=3r_0/4$, the estimates \eqref{ckg:cubic-eliminated-size} and the second in \eqref{ckg:C3-diff} imply
\eqref{ckg:E0-full}--\eqref{ckg:P0R-jet}, after enlarging $C_{\rm in}$. The stated symmetry decomposition follows from the gauge invariance of the NLS Hamiltonian and the symmetry preservation of the preceding transformations.

The coefficient estimates \eqref{ckg:N0R-est} follow from \eqref{ckg:omega-same-close}--\eqref{ckg:M-same-close} when $m=n$, and from
\eqref{ckg:omega-dist-close}--\eqref{ckg:HF-dist-close} when $m\ne n$. In the former case the rest energy cancels in every elliptic block; in the latter, one copy of $c^2$ remains in the generic elliptic frequencies.

Finally, the tangential twist matrix differs by $O(h)$ from
\begin{equation}\label{ckg:twist-matrix}
 A_*:=\frac1{4\pi}
 \begin{pmatrix}1&2\\2&1\end{pmatrix}.
\end{equation}
Then its singular values are $1/(4\pi)$ and $3/(4\pi)$ and from \eqref{ckg:omega-same}, \eqref{ckg:omega-dist} and the fixedness of the internal sites,
\begin{equation}\label{ckg:twist-close}
 \sup_{\rho\in\Xi_0} \|D_\rho\widehat\omega_0^{\rm KG}(\rho)-A_*\| \le Ch, \qquad D_\rho\omega_0^{\rm NLS}=A_*.
\end{equation}
Since $\Xi_0$ is convex, integrating the derivative along the segment joining $\rho$ and $\rho'$ gives
\[
 |\widehat\omega_0^{\rm KG}(\rho)-\widehat\omega_0^{\rm KG}(\rho')| \ge \left(\frac1{4\pi}-Ch\right)|\rho-\rho'|,
\]
and the corresponding upper bound follows from $\|D_\rho\widehat\omega_0^{\rm KG}\|\le\|A_*\|+Ch$.  After decreasing $h_*$ this proves the lipeomorphism nature of the frequency map; the same property for the NLS estimate follows directly from the affine map with matrix $A_*$.  The proof is complete.
\end{proof}

\section{The block homological equation}

We now solve the linear homological equation for a block normal form \eqref{ckg:block-N}.  The analytic part of the argument is the standard one of \cite[Lemma~1]{poschel1996kam}, the singular KG--NLS comparison follows the mechanism of \cite[Lemma~4.2]{BBG}, while for a reference using a block approach see for instance \cite[Sections 6-8]{EGK}.

\subsection{Perturbed normal forms and finite-dimensional homological operators}

Fix $0\le\theta_0\le1$ and set
\begin{equation*}
 \beta_0:=1-\theta_0.
\end{equation*}
With the initial coefficients defined in \eqref{ckg:initial-shifted-coefficients}, at a generic KAM step we write
\begin{align*}
 \widehat\omega &=\widehat\omega_0^{\rm KG}+\delta_\omega, \qquad  &\omega^{\rm NLS} =\omega_0^{\rm NLS}+\delta_\omega^{\rm NLS},
 \\
 \widetilde A_B &=\widetilde A_{0,B}^{\rm KG}+D_B, \qquad  &A_B^{\rm NLS} =A_{0,B}^{\rm NLS}+D_B^{\rm NLS},
 \\
 K_\cF &=K_{\cF,0}^{\rm KG}+E_{\cF}, \qquad  &K_\cF^{\rm NLS} =K_{\cF,0}^{\rm NLS}+E_{\cF}^{\rm NLS}.
\end{align*}
The last line is omitted if $\mathcal F=\varnothing$.
On the current parameter set $\Xi\subseteq\Xi_0$ we assume
\begin{equation}\label{ckg:freq-corrections}
 \begin{aligned}
 &|\delta_\omega|^{\mathcal L(\alpha)} +|D|_{\infty,\beta_0}^{\mathcal L(\alpha)} +\|E_{\cF}\|^{\mathcal L(\alpha)}
 \\
 &\qquad +|\delta_\omega^{\rm NLS}|^{\mathcal L(\alpha)} +|D^{\rm NLS}|_{\infty}^{\mathcal L(\alpha)} +\|E_{\cF}^{\rm NLS}\|^{\mathcal L(\alpha)} \le C_*\varepsilon_{\rm NF},
 \\
 &|\delta_\omega|_\infty^{\rm lip} +|D|_{\infty,\beta_0}^{\rm lip} +\|E_{\cF}\|_\infty^{\rm lip} +|\delta_\omega^{\rm NLS}|_\infty^{\rm lip} +|D^{\rm NLS}|_{\infty}^{\rm lip} +\|E_{\cF}^{\rm NLS}\|_\infty^{\rm lip} \le C_*\kappa_{\rm NF}.
 \end{aligned}
\end{equation}
for some $\varepsilon_{\rm NF} \in (0,1)$ and $\kappa_{\rm NF}\ge0$. The finite term is again omitted when $\mathcal F=\varnothing$.  

For the local KG--NLS comparison below we assume, on the common parameter set,
\begin{equation}\label{ckg:freq-KG-NLS-close}
 {\sup}_{\rho\in\Xi}|\widehat\omega(\rho)-\omega^{\rm NLS}(\rho)| +{\sup}_{\substack{\rho\in\Xi\\B\in\mathfrak B}} \langle B\rangle^{-4} \|\widetilde A_B(\rho)-A_B^{\rm NLS}(\rho)\| \le Ch,
\end{equation}
and, in the distinct-mode case,
\begin{equation}\label{ckg:freq-KG-NLS-close-F}
 {\sup}_{\rho\in\Xi}\|K_\cF(\rho)-K_\cF^{\rm NLS}(\rho)\| \le hC_{\cF}\Delta_{\cF},
\end{equation}
where $\Delta_{\cF}\ge0$ is a finite-block comparison size, which is going to remain uniformly bounded in the iterations and it is omitted in the identical-mode case.

Let $\mathbb L_{\cF}=-\iu\mathcal J_\cF K_\cF$ as in \eqref{ckg:LF-from-KF}, and denote by $\mathbb L_{\cF,0}^{\rm KG}$ and $\mathbb L_{\cF,0}^{\rm NLS}$ the corresponding initial matrices in \eqref{ckg:initial-shifted-coefficients}.

\begin{lemma}\label{ckg:lem-hyp-resolvent}
There exist $C_{0,\cF},C>0$, independent of $c,\rho,t$, such that, if
\begin{equation}\label{ckg:finite-resolvent-neighbourhood}
 \|\mathbb L_{\cF}-\mathbb L_{\cF,0}^{\rm KG}\|\le C_{0,\cF}\varrho,
\end{equation}
then
\begin{equation}\label{ckg:hyp-resolvent-allt}
 \|(\mathbb L_{\cF}+\iu tI)^{-1}\|\le\frac{C}{\varrho} \qquad\forall t\in\mathbb R.
\end{equation}
The same statement holds for the NLS finite block in a $C_{0,\cF}\varrho$-neighbourhood of $\mathbb L_{\cF,0}^{\rm NLS}$.
\end{lemma}

\begin{proof}
Let $\mathscr P_{\cF}$ be the fixed permutation sending $W_{\cF}=(a_n,b_m,\bar a_n,\bar b_m)$ to $(a_n,\bar b_m,\bar a_n,b_m)$.  By the definitions of ${L}_\cF^{(2)}$ and ${L}_{\cF,{\rm NLS}}^{(2)}$ in \eqref{ckg:LF-def} we have 
\begin{align*}
 \mathscr P_{\cF}\mathbb L_{\cF,0}^{\rm KG}\mathscr P_{\cF}^{-1}
 &=\operatorname{diag}\bigl(L_{\cF}^{(2)},\overline{L_{\cF}^{(2)}}\bigr),\\
 \mathscr P_{\cF}\mathbb L_{\cF,0}^{\rm NLS}\mathscr P_{\cF}^{-1} &=\operatorname{diag}\bigl(L_{\cF,\rm NLS}^{(2)}, \overline{L_{\cF,\rm NLS}^{(2)}}\bigr).
\end{align*}
Since $\|(\overline{L_{\cF}^{(2)}}+\iu tI)^{-1}\| =\|(L_{\cF}^{(2)}-\iu tI)^{-1}\|$, and similarly for NLS, Lemma~\ref{ckg:lem-exact-finite-resolvent} yields
\[
 \|(\mathbb L_{\cF,0}^{\rm KG}+\iu tI)^{-1}\| +\|(\mathbb L_{\cF,0}^{\rm NLS}+\iu tI)^{-1}\| \le C_{\cF}\varrho^{-1}, \qquad t\in\mathbb R.
\]
Choose $C_{0,\cF}<(2C_{\cF})^{-1}$.  Then
\[
 \mathbb L_{\cF}+\iu tI =\bigl[I+(\mathbb L_{\cF}-\mathbb L_{\cF,0}^{\rm KG}) (\mathbb L_{\cF,0}^{\rm KG}+\iu tI)^{-1}\bigr] (\mathbb L_{\cF,0}^{\rm KG}+\iu tI),
\]
and the first factor is invertible by a Neumann series, proving \eqref{ckg:hyp-resolvent-allt}. The NLS proof follows the same idea.
\end{proof}

We now define all finite-dimensional homological operators. For $B,C\in\mathfrak B$, set
\[
 \operatorname{Mat}(B,C):= \mathcal L(\mathbb C^C,\mathbb C^B).
\]
For $k\in\mathbb Z^2$, the elliptic homological operators are
\begin{align}
 \mathscr L_k^0&:\mathbb C\longrightarrow\mathbb C, &\mathscr L_k^0u&:=\iu\langle\omega,k\rangle u,
 \label{ckg:L0}\\
 \mathscr L_{k,B}^{1,\pm}&:\mathbb C^B\longrightarrow\mathbb C^B, &\mathscr L_{k,B}^{1,\pm}u &:=\iu\langle\omega,k\rangle u\pm\iu A_Bu,
  \notag\\
 \mathscr L_{k,B,C}^{-}&:\operatorname{Mat}(B,C) \longrightarrow\operatorname{Mat}(B,C), &\mathscr L_{k,B,C}^{-}X &:=\iu\langle\omega,k\rangle X+\iu A_BX-\iu XA_C,
 \label{ckg:Lminus}\\
 \mathscr L_{k,B,C}^{+}&:\operatorname{Mat}(B,C) \longrightarrow\operatorname{Mat}(B,C), &\mathscr L_{k,B,C}^{+}X &:=\iu\langle\omega,k\rangle X+\iu A_BX+\iu XA_C^{t}.
 \label{ckg:Lplus}
\end{align}
For a $y$-linear coefficient, $\mathscr L_k^0$ acts componentwise on $\mathbb C^2$.

When $\mathcal F\ne\varnothing$, let $V_{\cF}:=\mathbb C^{2d_\cF}$.  For row covectors $s\in V_{\cF}^*$ define
\begin{equation}\label{ckg:Lfinite-linear}
 \mathscr L_{k,\cF }^{1}:V_{\cF}^*\longrightarrow V_{\cF}^*, \qquad \mathscr L_{k,\cF }^{1}s :=s\bigl(\iu\langle\omega,k\rangle I-\mathbb L_{\cF}\bigr).
\end{equation}
For a mixed finite--elliptic coefficient $X\in\operatorname{Mat}(V_{\cF},B):=\mathcal L(\mathbb C^B,V_{\cF})$ we use the model operators
\begin{equation}\label{ckg:Lfinite-mixed}
 \mathscr L_{k,\cF ,B}^{1,\pm}X :=\iu\langle\omega,k\rangle X-\mathbb L_{\cF}^{t}X \pm\iu XA_B,
\end{equation}
acting from $\operatorname{Mat}(V_{\cF},B)$ to itself.  Depending on whether the finite or elliptic variable is conjugated, the actual coefficient equation is one of \eqref{ckg:Lfinite-mixed}, its transpose or its complex conjugate; all these operators have the same inverse norm.  Finally, on the space $\operatorname{Sym}(V_{\cF})$ of complex symmetric matrices we set
\begin{equation}\label{ckg:Lfinite-quadratic}
 \mathscr L_{k,\cF }^{2}X :=\iu\langle\omega,k\rangle X -\mathbb L_{\cF}^{t}X-X\mathbb L_{\cF}.
\end{equation}
This is the coefficient operator for a quadratic finite Hamiltonian $\frac12e^{\iu k\cdot\theta}\langle W_{\cF},XW_{\cF}\rangle$, with $X^{t}=X$.

For each $k$ and each choice of blocks, we denote by $\mathfrak X_{k,B}^{1,\pm}$, $\mathfrak X_{k,B,C}^{\pm}$ and the analogous finite spaces the subspaces of coefficients whose associated monomials satisfy the exact momentum and physical-charge selection rules.  By \eqref{ckg:block-exact-symmetries}, all operators above preserve the corresponding symmetry subspaces.  Every inverse below is understood on these invariant subspaces.

By Lemma~\ref{ckg:lem-hyp-resolvent}, unitary diagonalization of the Hermitian blocks $A_B$ and the uniform norm equivalence \eqref{ckg:finite-norm-equivalence}, we obtain the following.

\begin{corollary}\label{ckg:cor-one-finite-factor}
Under the assumptions of Lemma~\ref{ckg:lem-hyp-resolvent}, the operators \eqref{ckg:Lfinite-linear} and \eqref{ckg:Lfinite-mixed}, restricted to any symmetry-allowed coefficient space, are invertible and
\begin{equation}\label{ckg:one-finite-factor-inverse}
 \|(\mathscr L_{k,\cF }^{1})^{-1}\| +\|(\mathscr L_{k,\cF ,B}^{1,\pm})^{-1}\| \le\frac{C}{\varrho}.
\end{equation}
The same estimate holds for their transposes and complex conjugates, and for the corresponding NLS operators.
\end{corollary}

\begin{lemma}
\label{ckg:lem-finite-finite-nonlocal}
Consider the distinct-mode case.  A symmetry-allowed quadratic coefficient with two finite variables is either locally gauge invariant, in which case Lemma~\ref{ckg:lem-selection-dist} forces $k=0$ and the coefficient belongs to the normal projection, or it is non-locally-gauge-invariant.  In the latter case the restriction of \eqref{ckg:Lfinite-quadratic} to its symmetry subspace is invertible for $c$ sufficiently large and
\begin{equation}\label{ckg:finite-finite-c2-inverse}
 \|(\mathscr L_{k,\cF }^{2})^{-1}\|\le Cc^{-2}=Ch.
\end{equation}
\end{lemma}

\begin{proof}
If the coefficient is locally gauge invariant and contains only finite variables, Lemma~\ref{ckg:lem-selection-dist} gives $k=0$; by Definition~\ref{ckg:def-normal-projection} it is retained.

Assume now that the coefficient is not locally gauge invariant.  Since no generic elliptic variable occurs, its fast defect is $d=\mathfrak L(k,\ell)=k_1+k_2\ne0$.  Physical charge gives $|k_1-k_2|\le4$, whereas momentum gives $mk_1+nk_2=0$.  If $m+n\ne0$, these two relations leave only finitely many possible $k$, depending on the fixed sites $m,n$. Hence the shifted tangential part is $O(1)$, and \eqref{ckg:abg-dist}, \eqref{ckg:LF-from-KF}, and \eqref{ckg:freq-corrections} give
$\|\mathbb L_{\cF}\|\le C\varrho+C_*\varepsilon_{\rm NF}=O(1)$.
Thus, by \eqref{ckg:Lfinite-quadratic}, the operator has the form
\[
 \iu c^2d\,I+\mathcal A_{k,\cF }, \qquad \|\mathcal A_{k,\cF }\|\le C_{m,n}.
\]
For $c$ large, a Neumann series gives \eqref{ckg:finite-finite-c2-inverse}.

If $m+n=0$, momentum and $m\ne n$ imply $k_1=k_2$.  Thus $d=2k_1\ne0$ and
\[
 \langle\omega,k\rangle =k_1(\omega_1+\omega_2) =2k_1c^2+O(|k_1|).
\]
The two $\mathbb L_{\cF}$ contributions in \eqref{ckg:Lfinite-quadratic} are uniformly bounded, so the same Neumann argument is uniform for every $k_1\ne0$.  
\end{proof}

\subsection{Elliptic Sylvester operators}

\begin{lemma}\label{ckg:lem-spectral-reduction}
Let $\mu_{B,a}$, $a=1,\dots,d_B$, and $\mu_{C,b}$, $b=1,\dots,d_C$, be the ordered eigenvalues of the Hermitian matrices $A_B$ and $A_C$.  On the full coefficient spaces and with the Euclidean norm for vectors and the Hilbert--Schmidt norm for matrices,
\begin{align}
 \|(\mathscr L_{k,B}^{1,\pm})^{-1}\|^{-1} &=\min_a|\langle\omega,k\rangle\pm\mu_{B,a}|,
 \label{ckg:spectral-L1}\\
 \|(\mathscr L_{k,B,C}^{-})^{-1}\|_{\rm HS\to HS}^{-1} &=\min_{a,b}|\langle\omega,k\rangle+ \mu_{B,a}-\mu_{C,b}|,
 \label{ckg:spectral-Lminus}\\
 \|(\mathscr L_{k,B,C}^{+})^{-1}\|_{\rm HS\to HS}^{-1} &=\min_{a,b}|\langle\omega,k\rangle+ \mu_{B,a}+\mu_{C,b}|. \notag\end{align} More precisely, if $\delta$ is the minimum of the allowed absolute scalar divisors, then, for blocks of size at most four,
\[
 \|\mathscr L^{-1}\|_{\mathrm{op}\to\mathrm{op}}\le 2/\delta,
 \qquad \|X\|_{\mathrm{HS}}\le 2\|X\|.
\]
On the invariant symmetry spaces $\mathfrak X_{k,B}^{1,\pm}$ and $\mathfrak X_{k,B,C}^{\pm}$ the minima are restricted to the eigenvalue entries which belong to those spaces.  The symmetry-allowed coefficient $k=0$, $B=C$ of $z_B\bar z_B$ type is retained by $[\cdot]$ and is not inverted.
\end{lemma}

\begin{proof}
Because \eqref{ckg:block-exact-symmetries} holds, $A_B$ and $A_C$ commute with the restrictions of the diagonal momentum and physical-charge generators to their blocks.  We may therefore diagonalize them by unitary matrices which preserve the corresponding joint symmetry eigenspaces.  Let
\[
 U_B^*A_BU_B=\mathscr D_B, \qquad U_C^*A_CU_C=\mathscr D_C,
\]
where $\mathscr D_B$, $\mathscr D_C$ are real diagonal. For \eqref{ckg:Lminus} we use the Hilbert--Schmidt isometry $Y=U_B^*XU_C$.  Then
\[
 \widetilde{\mathscr L} Y:= U_B^*\mathscr L_{k,B,C}^{-}(X)U_C =\iu\bigl(\langle\omega,k\rangle Y +\mathscr D_BY-Y\mathscr D_C\bigr),
\]
and then 
\[
(\widetilde{\mathscr L} Y)_{a,b}= \iu ( \langle \omega,k \rangle +\mu_{B,a}-\mu_{C,b}) Y_{a,b},
\]
and since $\|Y\|_{\rm HS}^2= \sum_{a,b} |Y_{a,b}|^2$, it follows 
\[
\| \widetilde{\mathscr L}^{-1}\|_{\rm HS \to \rm HS} = \max |\langle \omega,k \rangle +\mu_{B,a}-\mu_{C,b}|^{-1},
\]
and by the isometry property we have \eqref{ckg:spectral-Lminus}.

For \eqref{ckg:Lplus} the estimate is the same, by using $Y=U_B^*X\overline{U_C}$.  The vector identity \eqref{ckg:spectral-L1} is the same computation with one block.  Restricting to a symmetry subspace simply removes the forbidden matrix entries. 
\end{proof}

The next gap lemma is needed only for the local KG--NLS comparison. 

\begin{lemma}
\label{ckg:lem-local-large-gap}
Assume \eqref{ckg:freq-corrections}.  There exist $J_*>0$, $c_*>0$ and $C>0$, depending only on the fixed tangential sites, such that every symmetry-allowed, locally gauge invariant elliptic KAM-jet coefficient which is not retained by $[\cdot]$ satisfies one of the following alternatives:
\begin{enumerate}[label=(\roman*)]
 \item every elliptic block appearing in the coefficient satisfies
 $\langle B\rangle\le J_*$;
 \item the corresponding KG and NLS homological operators are both
 uniformly invertible,
 \begin{equation}\label{ckg:local-large-gap-inverse}
  \|\mathscr L_{KG}^{-1}\|+ \|\mathscr L_{NLS}^{-1}\|\le C,
 \end{equation}
 and, for a quadratic coefficient involving $B$ and $C$,
 \begin{equation}\label{ckg:local-large-weight-comparable}
 C^{-1}\langle B\rangle\le\langle C\rangle \le C\langle B\rangle.
 \end{equation}
\end{enumerate}
The constants are uniform for $c\ge c_*$, $\rho\in\Xi_0$, and for $\varrho$ and the KAM corrections sufficiently small.
\end{lemma}

\begin{proof}
In the \emph{identical-mode case}, Lemma~\ref{ckg:lem-selection-same} gives $k=0$ in the local sector.  A linear elliptic monomial is excluded by momentum. For a quadratic difference $z_i\bar z_j$, the momentum identity \eqref{ckg:sym-same-M} gives $i=j$; the coefficient therefore belongs to a retained block.  An unretained same-sign quadratic pair, by \eqref{ckg:sym-same-M} has signed Fourier indices $m+d$ and $m-d$ with $d\ne0$.  At $h=0$ the sum of the two leading shifted normal frequencies is
\[
 \frac{(m+d)^2-m^2}{2}+ \frac{(m-d)^2-m^2}{2}=d^2.
\]
For KG the corresponding leading quantity is
\begin{equation}\label{ckg:same-convex-gap}
 \lambda_{m+d}+\lambda_{m-d}-2\lambda_m.
\end{equation}
For fixed $d\ne0$, \eqref{ckg:same-convex-gap} converges to $d^2$ as $c\to\infty$.  Moreover, for $d>0$ it is increasing in $d$, because $j\mapsto\lambda_j$ is strictly convex.  Hence, after fixing $c_*$, its minimum over $c\ge c_*$ and $d\in\mathbb Z\setminus\{0\}$ is strictly positive. For $\varrho+\varepsilon_{\rm NF}$ sufficiently small, the explicit $O(\varrho)$ shifts in \eqref{ckg:coeff-same} and \eqref{ckg:coeff-same-NLS}, together with the KAM corrections controlled by \eqref{ckg:freq-corrections}, do not close this gap.  Finally $||m+d|-|m-d||\le2|m|$, which proves \eqref{ckg:local-large-weight-comparable}.  Thus every unretained local coefficient in the identical-mode case satisfies alternative (ii).

Consider now the \emph{distinct-mode case}.  By Lemma~\ref{ckg:lem-selection-dist}, $|k|_1\le4$.  For a linear elliptic coefficient, momentum fixes its Fourier index in a finite set depending only on $m,n$, so alternative (i) applies.  For a quadratic difference, momentum fixes $i-j=d(k)$, with $|d(k)|\le C_{m,n}.$ If $d(k)=0$, local gauge invariance gives $k_1+k_2=0$, and momentum then gives $k=0$; the coefficient is retained.  Hence an unretained difference has $d(k)\ne0$.  Since
\begin{equation}\label{ckg:nu-difference-exact}
 |\nu_i(h)-\nu_j(h)| =\frac{c\,|i-j|\,|i+j|} {\sqrt{c^2+i^2}+\sqrt{c^2+j^2}},
\end{equation}
this difference is bounded away from zero when $|i|+|j|$ is large, uniformly for $c\ge c_*$.  Indeed, if $\max\{|i|,|j|\}\le c$, the right-hand side of \eqref{ckg:nu-difference-exact} is bounded below by a fixed multiple of $|i-j||i+j|$; if $\max\{|i|,|j|\}\ge c$, then $i$ and $j$ have the same sign for large indices and the right-hand side is bounded below by a fixed multiple of $c_*|i-j|$.  The tangential and KAM corrections are uniformly bounded because $|k|_1\le4$.  The NLS difference $|i^2-j^2|/2$ has the same property.

For a quadratic sum, momentum fixes $i+j=s(k)$ with $|s(k)|\le C_{m,n}$.  The elementary estimate $\nu_j(h) \ge\frac{|j|}{3}$ for $c\ge1$ and $|j|\ge1$, shows that $\nu_i(h)+\nu_j(h)$ tends uniformly to infinity with $|i|+|j|$; the NLS statement is immediate.  Thus all sufficiently large unretained local coefficients satisfy \eqref{ckg:local-large-gap-inverse}. In both the difference and sum cases, either $i-j$ or $i+j$ is bounded, which gives \eqref{ckg:local-large-weight-comparable}.  Increasing $J_*$ to cover the remaining finitely many coefficients completes the proof.
\end{proof}

\subsection{Melnikov conditions and solution of the homological equation}

Fix
\begin{equation*}
 0\le\theta_0\le1, \qquad 0\le\theta_1\le\frac{1-\theta_0}{2}, \qquad \tau\ge1, \qquad 0<\alpha<1,
\end{equation*}
and set
\begin{equation*}
 \beta_0:=1-\theta_0, \qquad \beta_1:=1-\theta_0-2\theta_1\ge0, \qquad \mu:=2\tau+6.
\end{equation*}
For an elliptic homological operator $\mathscr L$, let $\widetilde w(\mathscr L)$ be the minimum of the $\widetilde w_B$ of the elliptic blocks occurring in the coefficient, with the convention $\widetilde w(\mathscr L)=1$ when no elliptic block occurs.  On every symmetry-allowed elliptic coefficient which is not retained by $[\cdot]$, we impose
\begin{equation}\label{ckg:block-Melnikov}
 \|\mathscr L^{-1}\| \le \frac{\langle k\rangle^{\tau+1} \widetilde w(\mathscr L)^{\theta_1}}{\alpha}.
\end{equation}
A convenient sufficient formulation is
\begin{align}
 |\langle\omega,k\rangle|^{-1} &\le\frac{\langle k\rangle^{\tau+1}}{\alpha}, \qquad k\ne0, \qquad \|(\mathscr L_{k,B}^{1,\pm})^{-1}\| \le\frac{\langle k\rangle^{\tau+1} \widetilde w_B^{\theta_1}}{\alpha},
 \label{ckg:Melnikov-L1}\\
 \|(\mathscr L_{k,B,C}^{\pm})^{-1}\| &\le\frac{\langle k\rangle^{\tau+1} \min\{\widetilde w_B,\widetilde w_C\}^{\theta_1}}{\alpha}.
 \label{ckg:Melnikov-L2}
\end{align}
The inverse norms in \eqref{ckg:Melnikov-L1}--\eqref{ckg:Melnikov-L2} are understood on the exact symmetry subspaces defined above.  Coefficients with one finite factor are controlled by Corollary~\ref{ckg:cor-one-finite-factor}; finite--finite coefficients are controlled by Lemma~\ref{ckg:lem-finite-finite-nonlocal}.  Hence no Melnikov exclusion is imposed on finite-containing operators. For NLS we use the same conditions with $h=0$ and $\theta_1=0$; for the locally gauge invariant NLS coefficients of this paper, they are automatic by Lemma~\ref{ckg:lem-NLS-no-bad-strips}. The endpoint $\theta_0=1$, $\theta_1=0$ is the ordinary unweighted homological estimate; its proof is the same coefficient argument with no normal-weight transfer.

\begin{theorem}[Block homological equation]\label{ckg:thm-block-homological}
Fix $p>9/2$ and $m,n \in \mathbb{Z}$. There exist $h_*,\varrho_*,\varepsilon_*,c_\alpha,C>0$ such that the following holds.  Let $0<h<h_*$, $0<\varrho<\varrho_*$, $s,r>0$ and $0<\sigma<\min\{s/4,1\}$.  Let $N$ be a KG block normal form \eqref{ckg:block-N} satisfying \eqref{ckg:freq-corrections} and
\begin{equation}\label{ckg:homological-smallness}
 \frac{\varepsilon_{\rm NF}}{\varrho}+\kappa_{\rm NF} \le\varepsilon_*, \qquad 0<\alpha\le c_\alpha\varrho.
\end{equation}
Suppose that \eqref{ckg:block-Melnikov} holds on a current parameter set $\Xi\subseteq\Xi_0$.

Let $P^T$ be a real analytic, translation invariant and physically gauge invariant KG KAM jet, with finite weighted Lipschitz vector-field norms on $\mathfrak D_q(s,r)$ for $q=p,p-4$. Then there exist a normalized real KAM jet $S$ and a normal correction $\widehat N:=[P^T]$ such that
\begin{equation}\label{ckg:block-homological-eq}
 \{N,S\}+P^T=\widehat N, \qquad [S]=0.
\end{equation}
The constant component of $\widehat N$ is absorbed into the scalar part of the Hamiltonian, while the remaining part is a block normal form of the type \eqref{ckg:block-N}.  The generator $S$ preserves momentum and physical charge.  We have then, for $q=p,p-4$,
\begin{align}\label{ckg:block-homological-est}
 \|X_S\|_{r,q,\beta_1}^{\mathcal L(\alpha), \mathfrak D_q(s-\sigma,r)} \le \frac{C}{\alpha\sigma^\mu} \|X_{P^T}\|_{r,q,\beta_0}^{\mathcal L(\alpha), \mathfrak D_q(s,r)}, \\
\label{ckg:Nhat-est}
 \|X_{\widehat N}\|_{r,q,\beta_0}^{\mathcal L(\alpha), \mathfrak D_q(s,r)} \le C\|X_{P^T}\|_{r,q,\beta_0}^{\mathcal L(\alpha), \mathfrak D_q(s,r)}.
\end{align}

Similarly, let $N^{\rm NLS}$ be the corresponding NLS block normal form, satisfying the analogous smallness and finite-resolvent assumptions, and suppose that the NLS Melnikov conditions hold. If $P^{T,\rm NLS}$ is a real analytic, translation invariant and locally gauge invariant KAM jet, with finite vector-field Lipschitz norms on $\mathfrak D_q(s,r)$ for $q=p,p-4$, then there exist $S^{\rm NLS}$ and $\widehat N^{\rm NLS}=[P^{T,\rm NLS}]$ such that
\begin{equation}\label{ckg:block-homological-NLS-eq}
 \{N^{\rm NLS},S^{\rm NLS}\}+P^{T,\rm NLS} =\widehat N^{\rm NLS}, \qquad [S^{\rm NLS}]=0,
\end{equation}
and, for $q=p,p-4$,
\begin{align}
 \|X_{S^{\rm NLS}}\|_{r,q}^{\mathcal L(\alpha), \mathfrak D_q(s-\sigma,r)} &\le \frac{C}{\alpha\sigma^\mu} \|X_{P^{T,\rm NLS}}\|_{r,q}^{\mathcal L(\alpha), \mathfrak D_q(s,r)},
  \notag\\
 \|X_{\widehat N^{\rm NLS}}\|_{r,q}^{\mathcal L(\alpha), \mathfrak D_q(s,r)} &\le C\|X_{P^{T,\rm NLS}}\|_{r,q}^{\mathcal L(\alpha), \mathfrak D_q(s,r)}.
 \label{ckg:Nhat-NLS-est}
\end{align}

Finally assume \eqref{ckg:freq-KG-NLS-close} and, when present, \eqref{ckg:freq-KG-NLS-close-F}, work on a common parameter set on which the
KG and NLS Melnikov conditions above hold. Define
\begin{equation*}
 P^{T,\mathcal{NG}}:=(I-\Pi_{\rm LG})P^T, \qquad P^{T,R}:=\Pi_{\rm LG}P^T-P^{T,\rm NLS}.
\end{equation*}
Then $P^T=P^{T,\rm NLS}+P^{T,\mathcal{NG}}+P^{T,R}$ and $P^{T,R}$ is locally gauge invariant.  Assume that, for some $\varepsilon_R\ge0$,
\begin{equation}\label{ckg:Rjet-diff}
 \|X_{P^{T,R}}\|_{r,p-4}^{\mathfrak D_p(s,r)} \le h\varepsilon_R,
\end{equation}
and put
\begin{align*}
 \varepsilon_{\mathcal{NG}} &:={} \|X_{P^{T,\mathcal{NG}}}\|_{r,p,\beta_0}^{\mathcal L(\alpha), \mathfrak D_p(s,r)} +\|X_{P^{T,\mathcal{NG}}}\|_{r,p-4,\beta_0}^{\mathcal L(\alpha), \mathfrak D_{p-4}(s,r)},
 \\
 \varepsilon_{\rm NLS} &:={} \|X_{P^{T,\rm NLS}}\|_{r,p}^{\mathcal L(\alpha), \mathfrak D_p(s,r)} +\|X_{P^{T,\rm NLS}}\|_{r,p-4}^{\mathcal L(\alpha), \mathfrak D_{p-4}(s,r)}.
\end{align*}
Then
\begin{equation}\label{ckg:block-generator-comparison}
 \|X_{S-S^{\rm NLS}}\|_{r,p-4}^{\mathfrak D_p(s-\sigma,r)} \le Ch\left[ \frac{\varepsilon_R+\varepsilon_{\mathcal{NG}}} {\alpha\sigma^\mu} +\frac{(1+\Delta_{\cF})\varepsilon_{\rm NLS}} {\alpha^2\sigma^\mu} \right].
\end{equation}
The factor $\Delta_{\cF}$ is omitted when $\mathcal F=\varnothing$.  Moreover
\begin{equation}\label{ckg:Nhat-comparison}
 \|X_{\widehat N-\widehat N^{\rm NLS}}\|_{r,p-4}^{ \mathfrak D_p(s,r)} \le Ch\varepsilon_R.
\end{equation}
All constants are independent of $c,h,r,\alpha,\sigma$ and of the block labels; they may depend on $p,a,\tau$ and on the fixed internal sites.
\end{theorem}

\begin{proof}
Write the KAM jet as $P^T=P^{\rm tan}+P^{(1)}+P^{(2)}$, where $P^{\rm tan}=P^{00}(\theta)+\langle P^y(\theta),y\rangle$ and $P^{(j)}$ is homogeneous of degree $j$ in the normal variables.  On each Fourier mode $k$, the coefficient equations associated with \eqref{ckg:block-homological-eq} are
\begin{equation}\label{ckg:proof-coefficient-equations}
 \mathscr L_{k,\bullet}S_{k,\bullet}=-P^T_{k,\bullet}
 \quad\text{outside the range of }[\cdot],\qquad [S]=0.
\end{equation}
Here $\mathscr L_{k,\bullet}$ is one of \eqref{ckg:L0}--\eqref{ckg:Lplus} or, if $\mathcal F\ne\varnothing$, one of \eqref{ckg:Lfinite-linear}--\eqref{ckg:Lfinite-quadratic}. The projection $[\cdot]$ retains the zero Fourier scalar and $y$-linear coefficients, the entire symmetry-allowed Hermitian $k=0$, $B=C$ block of $z_B\bar z_B$ type, and the admissible finite quadratic block.
Every other elliptic equation is invertible by \eqref{ckg:block-Melnikov}, the equations with one finite factor by Corollary~\ref{ckg:cor-one-finite-factor}, and the remaining finite--finite equations by Lemma~\ref{ckg:lem-finite-finite-nonlocal}.  By \eqref{ckg:homological-smallness}, the finite-resolvent neighbourhood \eqref{ckg:finite-resolvent-neighbourhood} holds after reducing the smallness constants.  These inversions are performed on the exact symmetry subspaces;
by \eqref{ckg:block-exact-symmetries} and the reality involution they preserve momentum, physical charge and reality.  They determine a unique normalized solution, and the retained coefficients give $\widehat N=[P^T]$.

To estimate the solution, consider for instance $\mathscr L_{k,B,C}^{\pm}$ with $B=[i]$, $C=[j]$, and endow a quadratic coefficient with the target-weighted block norm
\[
 \|M_{BC}\|_{q,\beta;B\leftarrow C}:= \widetilde w_B^\beta\frac{e^{a|i|}\langle i\rangle^q}{e^{a|j|}\langle j\rangle^q} \|M_{BC}\|,\qquad q\in\{p,p-4\}.
\]
This is the operator norm after conjugation by the normal weights in \eqref{eq:norm-apb}, including the one-sided target factor in \eqref{ckg:scaled-X-norm}.  The vector-field component with the opposite orientation is estimated with $B$ and $C$ interchanged.  With $w_{BC}:=\min\{\widetilde w_B,\widetilde w_C\}$, \eqref{ckg:block-Melnikov} gives, in both quadratic orientations,
\begin{equation}\label{ckg:proof-weighted-coefficient}
 \|S_{k,B,C}^{\pm}\|_{q,\beta_1;B\leftarrow C}
 \le\frac{\langle k\rangle^{\tau+1}}{\alpha}
 \|P_{k,B,C}^{\pm}\|_{q,\beta_0;B\leftarrow C}.
\end{equation}
Indeed, $\beta_0-\beta_1=2\theta_1$ and, crucially,
\begin{equation}\label{ckg:proof-weight-absorption}
 \widetilde w_B^{\beta_1}w_{BC}^{j_\sharp \theta_1}
 \le\widetilde w_B^{\beta_0},\qquad j_\sharp=1,2.
\end{equation}
For linear elliptic coefficients the same argument uses only the one-block weight $\widetilde w_B$.  In particular, the case $j_\sharp=2$ absorbs the two Melnikov inverses arising in parameter differences.  The resolvent identity
\[
 \mathscr L(\rho)^{-1}-\mathscr L(\rho')^{-1}
 =\mathscr L(\rho)^{-1}
 [\mathscr L(\rho')-\mathscr L(\rho)]\mathscr L(\rho')^{-1}
\]
and the Lipschitz bounds on the reference frequencies and on \eqref{ckg:freq-corrections} therefore give
\begin{equation}\label{ckg:proof-parameter-coefficient}
 \alpha\sup_{\substack{\rho,\rho'\in\Xi\\\rho\ne\rho'}}
 \frac{\|S_{k,B,C}^{\pm}(\rho)-S_{k,B,C}^{\pm}(\rho')\|_{q,\beta_1;B\leftarrow C}}
 { |\rho-\rho'|}
 \le\frac{C\langle k\rangle^{2\tau+3}}{\alpha}
 \|P_{k,B,C}^{\pm}\|_{q,\beta_0;B\leftarrow C}^{\mathcal L(\alpha)}.
\end{equation}
Here the weighted Lipschitz norm of a block coefficient is defined by the convention \eqref{ckg:finite-Lip-convention}, using the displayed block norm.  The same argument gives the corresponding one-block estimate for linear elliptic coefficients.

For finite-linear and mixed finite--elliptic coefficients recall Corollary \ref{ckg:cor-one-finite-factor}.  In the mixed case,
\[
\begin{aligned}
 &\bigl(\mathscr L_{k,\cF ,B}^{1,\pm}(\rho)
       -\mathscr L_{k,\cF ,B}^{1,\pm}(\rho')\bigr)X\\
 &\quad=\iu\langle\omega(\rho)-\omega(\rho'),k\rangle X
 -(\mathbb L_{\cF}(\rho)-\mathbb L_{\cF}(\rho'))^{t}X
 \pm\iu X(A_B(\rho)-A_B(\rho')).
\end{aligned}
\]
The uniformly bounded parameter-Lipschitz seminorms of the reference coefficients, together with \eqref{ckg:freq-corrections} and \eqref{ckg:homological-smallness}, bound this difference by $C\langle k\rangle|\rho-\rho'|$. The same bound holds for the finite-linear operator.  The resolvent identity thus gives
\[
 \|\mathscr L^{-1}\|^{\mathcal L(\alpha)}
 \le C\bigl(\varrho^{-1}
         +\alpha\langle k\rangle\varrho^{-2}\bigr)
 \le C\langle k\rangle\alpha^{-1},
\]
by $\alpha\le c_\alpha\varrho$.  The product estimate for $S_k=-\mathscr L^{-1}P_k$ now yields the finite-containing analogues of the preceding coefficient bounds.  The nonlocal finite--finite inverses satisfy stronger estimates by Lemma~\ref{ckg:lem-finite-finite-nonlocal}.  For mixed coefficients, the elliptic vector-field component has the same target block on both sides of the coefficient equation; its target weight is therefore absorbed by $\beta_1\le\beta_0$.
The finite component carries no target weight.

At fixed $k$, momentum restricts every row and column of a quadratic normal coefficient to a uniformly bounded number of elliptic blocks. The coefficient inversion respects this support; hence Lemma~\ref{ckg:lem-sparse-block-Schur} reconstructs the full normal operator from the coefficient estimates with no loss in the block labels.  The other jet terms obey the same scaled estimates.
Fourier Cauchy estimates and the standard analytic argument of \cite[Lemma~1]{poschel1996kam} then apply.  The angle component $X_S^y=\partial_\theta S$ costs at most one further factor $\langle k\rangle$.  Thus the coefficientwise Lipschitz loss $\langle k\rangle^{2\tau+3}$, this angle derivative, and the two-dimensional Fourier summation give the strip exponent $\mu=2\tau+6$, proving \eqref{ckg:block-homological-est}.
The projection bound \eqref{ckg:normal-projection-bounded} gives \eqref{ckg:Nhat-est}.  Repeating the coefficient argument for the NLS operators, with $\theta_1=0$ and no target weights, proves \eqref{ckg:block-homological-NLS-eq}--\eqref{ckg:Nhat-NLS-est}; in particular, the NLS hyperbolic block remains covered by Lemma~\ref{ckg:lem-hyp-resolvent}.

We turn to the estimate of the KG-NLS difference, following the cut-off approach of \cite[Lemma~4.2]{BBG}.  Since $\operatorname{ad}_N$ commutes with $\Pi_{\rm LG}$ by \eqref{ckg:block-fast-commute}, the normalized solution $S^{\mathcal{NG}}$ driven by $P^{T,\mathcal{NG}}$ has only coefficients with a nonzero integer fast defect $d$.  Let $\Pi_{\rm int}$ project onto elliptic normal indices $|j|<c/2$, and $\Pi_{\rm ext}:=I-\Pi_{\rm int}$.  When all elliptic indices are low and $|k|_1\le\kappa_*\sqrt c$, the coefficient operator has the form $\mathscr L=\iu c^2dI+\mathscr A$, with $d\ne0$ as in \eqref{ckg:fast-defect}. The shifted tangential contribution is $O(\langle k\rangle)$, there are at most two elliptic factors with $\nu_j(h)<c^2/8$ in this range, and the finite matrix is $O(\varrho)$.  By choosing $\kappa_*$ and the smallness parameters appropriately, $\|\mathscr A\|\le\tfrac12c^2|d|$ and thus
\begin{equation}\label{ckg:NG-low-inverse}
 \|\mathscr L^{-1}\|\le\frac{2}{c^2|d|}\le 2h.
\end{equation}
This also treats coefficients with one finite factor.  Purely finite nonlocal coefficients satisfy the corresponding bound by Lemma~\ref{ckg:lem-finite-finite-nonlocal}.  For low elliptic indices but $|k|_1>\kappa_*\sqrt c$, use the ordinary homological estimate. Since $h=c^{-2}$, on this range one has $1\le Ch\langle k\rangle^4$ and the Fourier tail is therefore bounded by $Ch\sigma^{-(\tau+7)}$.

For high elliptic indices the loss of four derivatives supplies a stronger factor.  For every $\beta\ge0$,
\begin{equation}\label{ckg:four-derivative-cutoff}
 \|\Pi_{\rm ext}v\|_{a,p-4,\beta}
 \le Cc^{-4}\|v\|_{a,p,\beta}
 =Ch^2\|v\|_{a,p,\beta}.
\end{equation}
To use this estimate in both orientations of a quadratic normal coefficient, write $S_k=\Pi_{\rm int}S_k\Pi_{\rm int}+\Pi_{\rm ext}S_k\Pi_{\rm int}+S_k\Pi_{\rm ext}$. For the second term solve at regularity $p$ and apply the cutoff on the output; for the third solve at regularity $p-4$ and apply it to the input.  More precisely,
\begin{equation}\label{ckg:proof-oriented-high}
 \begin{aligned}
 \|\Pi_{\rm ext}S_k\Pi_{\rm int}\|_{\ell^{a,p}\to\ell^{a,p-4}}
 &\le Ch^2\|S_k\|_{\ell^{a,p}\to\ell^{a,p}},\\
 \|S_k\Pi_{\rm ext}\|_{\ell^{a,p}\to\ell^{a,p-4}}
 &\le Ch^2\|S_k\|_{\ell^{a,p-4}\to\ell^{a,p-4}}.
 \end{aligned}
\end{equation}
These inequalities hold as well with the relevant target weights. The analogous one-sided argument applies to linear normal and finite--elliptic coefficients.  It is essential here that the ordinary homological estimate has been established at both $p$ and $p-4$.

These estimates also control the \emph{entire} Hamiltonian vector field.  In particular, Lemma~\ref{ckg:lem-sparse-block-Schur}, together with Cauchy--Schwarz, controls the bilinear expressions defining $X_{S^{\mathcal{NG}}}^y$ as well as the normal gradients.  For low indices, $\mathscr L_k$ commutes with $\partial_\theta$ on the $k$-th Fourier mode. Applying the inverse directly to the $\theta$-derivative of the datum avoids an additional angular factor beyond that already in the norm of $X_P^y$.
In the high-index range, \eqref{ckg:proof-oriented-high} also controls the bilinear contractions; for $X_S^y$ they cost at most an additional $\langle k\rangle$.  The resulting low-index, high-angular and high-normal bounds are respectively $Ch\sigma^{-2}$, $Ch/(\alpha\sigma^{\tau+7})$, and $Ch^2/(\alpha\sigma^{\tau+4})$. Since $\tau\ge1$ and $\mu=2\tau+6$, all are bounded by $Ch/(\alpha\sigma^\mu)$.  We have proved
\begin{equation}\label{ckg:SNG-est}
 \|X_{S^{\mathcal{NG}}}\|_{r,p-4}^{\mathfrak D_p(s-\sigma,r)}
 \le\frac{Ch}{\alpha\sigma^\mu}\varepsilon_{\mathcal{NG}}.
\end{equation}

Let $S^R$ and $\widetilde S^{\rm NLS}$ denote the normalized KG solutions driven by $P^{T,R}$ and $P^{T,\rm NLS}$, respectively. Linearity gives
\begin{equation}\label{ckg:proof-comparison-decomposition}
 S-S^{\rm NLS}=S^{\mathcal{NG}}+S^R+
                 (\widetilde S^{\rm NLS}-S^{\rm NLS}).
\end{equation}
Both right-hand sides $P^{T,R}$ and $P^{T,\rm NLS}$ are locally gauge invariant.  By Lemma~\ref{ckg:lem-local-large-gap}, each unretained elliptic coefficient is either supported on a fixed bounded family of blocks, where its inverse is $O(\alpha^{-1})$, or belongs to the large-index alternative, where the KG and NLS inverses are $O(1)$. In the latter case the momentum relations in the proof of that lemma give $|i-j|\le C_{m,n}$ or $|i+j|\le C_{m,n}$, hence
$\bigl||i|-|j|\bigr|\le C_{m,n}$.  In particular, the source--target Sobolev and exponential weights are uniformly comparable.  This removes the weighted Melnikov loss for large local indices.  Locally gauge invariant mixed finite--elliptic coefficients involve only bounded elliptic indices: the selection rule gives $|k|_1\le4$ and momentum then fixes the normal index.  Their inverses are $O(\varrho^{-1})=O(\alpha^{-1})$ by Corollary~\ref{ckg:cor-one-finite-factor}.
Hence the assumption \eqref{ckg:Rjet-diff}, with the coefficient estimates at the required regularities, implies
\begin{equation}\label{ckg:SR-est}
 \|X_{S^R}\|_{r,p-4}^{\mathfrak D_p(s-\sigma,r)}
 \le\frac{Ch}{\alpha\sigma^\mu}\varepsilon_R.
\end{equation}

To estimate the third term in the r.h.s. of \eqref{ckg:proof-comparison-decomposition}, it remains to compare the KG and NLS inversions on the same local right-hand side.  Denote their operators on an exact coefficient space by $L_h$ and $L_0$.  The basic identity is
\begin{equation}\label{ckg:KG-NLS-resolvent-identity}
 L_h^{-1}-L_0^{-1}=L_h^{-1}(L_0-L_h)L_0^{-1}.
\end{equation}
On the bounded family of Lemma~\ref{ckg:lem-local-large-gap}, $|k|_1$ and the normal weights are uniformly bounded, and \eqref{ckg:freq-KG-NLS-close} gives $\|L_h-L_0\|\le Ch$. Two inverse bounds therefore yield $Ch/\alpha^2$. On the large-index alternative, both inverses are $O(1)$, and the fourth-order bound in \eqref{ckg:freq-KG-NLS-close} is absorbed by the loss of four derivatives.  In particular,
\begin{equation}\label{ckg:proof-four-loss}
 \langle B\rangle^{p-4}
 \|\widetilde A_B-A_B^{\rm NLS}\|
 \langle B\rangle^{-p}\le Ch.
\end{equation}
For a quadratic Sylvester operator, left multiplication is handled by \eqref{ckg:proof-four-loss}.  For right multiplication, the momentum relations recalled above give $\bigl||i|-|j|\bigr|\le C_{m,n}$, so the exponential ratio $e^{a(|i|-|j|)}$ is bounded; moreover $\langle i\rangle\asymp\langle j\rangle$, and the four powers of the input weight absorb the $O(h\langle j\rangle^4)$ multiplier difference.  The conjugate orientations are identical.
Consequently, using the two uniform inverse bounds in \eqref{ckg:local-large-gap-inverse}, the resolvent identity \eqref{ckg:KG-NLS-resolvent-identity} gives an $O(h)$ difference of the large-index local inverses from regularity $p$ to $p-4$.

If there is one finite factor, the local symmetry and momentum constraints again force the elliptic index to be bounded. Both mixed inverses are $O(\varrho^{-1})$.  In their difference, the elliptic and tangential multipliers contribute $O(h)$ by \eqref{ckg:freq-KG-NLS-close}, whereas the finite multiplier contributes $O(h\Delta_{\cF})$ by \eqref{ckg:freq-KG-NLS-close-F}.  Hence
\begin{equation}\label{ckg:finite-resolvent-comparison}
 \|L_h^{-1}-L_0^{-1}\|
 \le\frac{Ch(1+\Delta_{\cF})}{\varrho^2}
 \le\frac{Ch(1+\Delta_{\cF})}{\alpha^2}.
\end{equation}
The locally gauge invariant finite--finite block is retained by $[\cdot]$.  The three estimates and finite local Fourier support give
\[
 \|X_{\widetilde S^{\rm NLS}-S^{\rm NLS}}\|_{r,p-4}^{
 \mathfrak D_p(s-\sigma,r)}
 \le\frac{Ch(1+\Delta_{\cF})}{\alpha^2\sigma^\mu}
 \varepsilon_{\rm NLS},
\]
with the factor involving $\Delta_{\cF}$ omitted if $\mathcal F=\varnothing$.  Combining this with \eqref{ckg:SNG-est}, \eqref{ckg:SR-est} and \eqref{ckg:proof-comparison-decomposition} proves \eqref{ckg:block-generator-comparison}.  Finally, \eqref{ckg:normal-projection-LG} implies $\widehat N-\widehat N^{\rm NLS}=[P^{T,R}]$. The same blockwise projection argument as for \eqref{ckg:normal-projection-bounded}, applied to the $p\to p-4$ norm, together with \eqref{ckg:Rjet-diff}, yields \eqref{ckg:Nhat-comparison}.
\end{proof}

\section{Measure estimates}
\label{sec:CKG-measure-estimates}

We now verify the block Melnikov conditions on a set of asymptotically full measure.  

For the coupled NLS reference system, our two mode choice yields no Melnikov parameter exclusions for the condition \eqref{ckg:block-Melnikov}: every symmetry-allowed elliptic coefficient not retained by the normal projection has a uniform spectral gap, and finite-containing operators, when present, are invertible for all parameters (Lemma~\ref{ckg:lem-NLS-no-bad-strips}). Thus all exclusions below come from KG. We nevertheless construct the two KAM transformations on the same retained parameter set in order to compare their embeddings and slow frequencies.

Three points require care.  First, the tangential parameter space has dimension two, whereas the scalar measure argument in \cite{BBG} invokes at one stage a non-degeneracy lemma stated for at least three tangential sites.  We therefore prove the required strip transversality directly from the explicit two-dimensional twist.  Second, the identical-mode normal frequencies are eigenvalues of $2\times2$ Hermitian matrices and may cross.  Accordingly, all weighted large-index Lipschitz asymptotics are proved for the matrices themselves; only pointwise spectral tails are deduced from Weyl's inequality.  Third, in the distinct-mode case, after the preliminary weighted KAM step the resulting $O(\epsilon_0)$ normal-form correction is frozen into the new reference frequencies.  This is the precise analogue of the separation between the first weighted step and the subsequent iteration in \cite[Section~6]{BBG}.

\subsection{Branch divisors and reference tails}

By Lemma~\ref{ckg:lem-spectral-reduction}, every elliptic homological operator is reduced to finitely many scalar functions of the form
\begin{equation}\label{ckg:branch-divisors}
 \Phi(\rho)=\langle\omega(\rho),k\rangle, \qquad \langle\omega(\rho),k\rangle\pm\mu_{B,a}(\rho), \qquad \langle\omega(\rho),k\rangle+ \mu_{B,a}(\rho)\pm\mu_{C,b}(\rho).
\end{equation}
For a branch involving normal blocks we put $\widetilde w(\Phi)=\min\{\widetilde w_B,\widetilde w_C\}$, with the obvious modification for one block, and $\widetilde w(\Phi)=1$ for a tangential branch.  The corresponding resonant strip is
\begin{equation}\label{ckg:branch-resonant-strip}
 \mathcal R_\Phi^{\theta_1}(\alpha):= \left\{\rho\in\Xi_0: |\Phi(\rho)|< \frac{2\alpha}{\langle k\rangle^{\tau} \widetilde w(\Phi)^{\theta_1}} \right\}.
\end{equation}
The exponent $\tau$ and the harmless factor $2$ are intentional.  The Hilbert--Schmidt reduction of Lemma~\ref{ckg:lem-spectral-reduction}, followed by the operator-norm inequality there, shows that failure of \eqref{ckg:block-Melnikov} implies, on a symmetry-allowed branch,
\[
 |\Phi(\rho)|< \frac{2\alpha}{\langle k\rangle^{\tau+1} \widetilde w(\Phi)^{\theta_1}},
\]
which is contained in \eqref{ckg:branch-resonant-strip}.  Using the factor $2$ in the \emph{excluded} strips as well as in their measure bounds is necessary for this implication; it changes only absolute constants in the latter.  The extra factor $\langle k\rangle^{-1}$ is recovered from transversality.

For each branch associated with an NLS symmetry-allowed coefficient, the corresponding NLS branch $\Phi^{\rm NLS}$ is obtained by replacing the KG tangential and elliptic normal frequencies in $\Phi$ with their NLS counterparts. We set $\mathcal R_\Phi^{\rm NLS}(\alpha) :=\mathcal R_{\Phi^{\rm NLS}}^0(\alpha)$.

First, we observe that combining Lemma~\ref{ckg:lem-spectral-reduction}, Corollary~\ref{ckg:cor-one-finite-factor}, and Lemma~\ref{ckg:lem-finite-finite-nonlocal}, we obtain the following.
\begin{lemma}[Finite multiplicity reduction]\label{ckg:lem-finite-multiplicity}
For fixed $k$ and fixed elliptic blocks, the bad set of every elliptic block operator is contained in the union of at most sixteen branch strips of the form \eqref{ckg:branch-resonant-strip}.  Operators containing exactly one finite hyperbolic factor contribute no resonant parameters provided the accumulated finite correction is $o(\varrho)$ and $0<\alpha\le c_\alpha\varrho$, with $c_\alpha$ as in \eqref{ckg:homological-smallness}.  Finite-finite coefficients outside the normal form are non-locally-gauge-invariant and are uniformly invertible by Lemma~\ref{ckg:lem-finite-finite-nonlocal}.
\end{lemma}

The momentum relation is generated by the conserved Hamiltonian \eqref{eq:NQM-zeta} and does not see the species label.  In the distinct-mode coordinates it is
\begin{equation}\label{ckg:momentum-kell}
 mk_1+nk_2+\sum_jj\ell_j=0,
\end{equation}
for the generic elliptic variables.  Hence, with $J_*:=\max\{|m|,|n|,1\}$,
\begin{align*}
 |\ell|_1=1 \Rightarrow |j|\le J_*|k|_1, \, \,   \ell=\pm(\delta_i+\delta_j)\Rightarrow |i+j|\le J_*|k|_1, \, \, \ell=\pm(\delta_i-\delta_j)\Rightarrow |i-j|\le J_*|k|_1.
\end{align*}
In the identical-mode coordinates the corresponding identity is
\begin{equation}\label{ckg:momentum-kell-same}
 m(k_1+k_2)+\sum_{j\ne m}(j-m)\ell_j=0.
\end{equation}
In the locally gauge invariant identical-mode sector, Lemma~\ref{ckg:lem-selection-same} gives $k=0$, and hence \eqref{ckg:momentum-kell-same} reduces simply to $\sum_{j\ne m}(j-m)\ell_j=0$.  In the non-locally-gauge-invariant sector we use \eqref{ckg:momentum-kell-same} together with $k=(r,r)$, $r\ne0$, again by Lemma~\ref{ckg:lem-selection-same}.

The normal forms \eqref{ckg:Heff-same} and \eqref{ckg:Heff-dist} also provide a useful large-index description of the reference branches.

\begin{lemma}\label{ckg:lem-intrinsic-tail}
Let $B=[j]\in\mathfrak B$ be an elliptic block with $|j|$ sufficiently large and write $A^0_B:=q_Bc^2I+\widetilde A_{0,B}^{\rm KG}$ with $\widetilde A_{0,B}^{\rm KG}$ as in \eqref{ckg:initial-shifted-coefficients}, for the corresponding full bare reference block.  After a fixed permutation of the coordinates inside $B$, independent of $\rho,j,c$, it has the form
\begin{equation}\label{ckg:intrinsic-tail}
 A^0_B =q_Bc^2 I+\nu_j(h)I+B_B^{\infty}(\rho;c)+R^0_B(\rho;c), \qquad \|R^0_B\|^{\mathcal L(\varrho)} \le C\frac{\varrho}{\widetilde w_j}.
\end{equation}
Here $B_B^{\infty}$ is Hermitian, independent of the large Fourier index on each fixed asymptotic block type.  In the distinct-mode case one may take $B_B^{\infty}=0$ on every generic elliptic block.  In the identical-mode case, on each $2\times2$ species component,
\begin{equation}\label{ckg:same-B-infty}
 B^{\infty}(\rho;c) =-\nu_m(h)I- \frac1{4\pi\widetilde w_m^2}
 \begin{pmatrix}
  \rho_1+2\rho_2&0\\
  0&\rho_2+2\rho_1
 \end{pmatrix},
\end{equation}
and
\begin{equation}\label{ckg:same-Rj-matrix}
 R_j^0(\rho;c)= \frac{\rho_1+\rho_2}{2\pi\widetilde w_m\widetilde w_j}I +\frac{\sqrt{\rho_1\rho_2}} {2\pi\widetilde w_m\widetilde w_j}
 \begin{pmatrix}0&1\\1&0\end{pmatrix}.
\end{equation}
The ordered eigenvalues $\beta_{B,a}$ of $B_B^{\infty}$ are Lipschitz on $\Xi_0$.  If $\mu^0_{B,a}$ denotes the corresponding ordered eigenvalue of $A^0_B$, then, pointwise in $\rho$,
\begin{equation}\label{ckg:intrinsic-tail-eigen-pointwise}
 \left|\mu^0_{B,a}- \bigl(q_Bc^2+\nu_j(h)+\beta_{B,a}(\rho;c)\bigr)\right| \le C\frac{\varrho}{\widetilde w_j}.
\end{equation}
No weighted Lipschitz estimate is asserted for the difference on the left-hand side of \eqref{ckg:intrinsic-tail-eigen-pointwise}.  In the identical-mode case the two limiting eigenvalues are
\begin{equation}\label{ckg:same-beta-infty}
 \beta_\pm(\rho;c)= -\nu_m(h)-\frac{3(\rho_1+\rho_2)}{8\pi\widetilde w_m^2} \pm\frac{|\rho_1-\rho_2|}{8\pi\widetilde w_m^2}.
\end{equation}
\end{lemma}

\begin{proof}
For the distinct-mode generic block, \eqref{ckg:Lambda-dist} is scalar in the species variables and gives
\[
 A^0_B=(c^2+\nu_j(h))I+ \frac1{2\pi\widetilde w_j} \left(\frac{\rho_1}{\widetilde w_m} +\frac{\rho_2}{\widetilde w_n}\right)I.
\]
Thus $B_B^{\infty}=0$, and both the supremum and the scaled Lipschitz seminorm of the remainder are bounded by $C\varrho/\widetilde w_j$.

In the identical-mode case, using $\lambda_j-\lambda_m=\nu_j(h)-\nu_m(h)$ in \eqref{ckg:coeff-same}, the matrix \eqref{ckg:Mk-same} is exactly
\[
 M_j=\nu_j(h)I+B^{\infty}+R_j^0,
\]
with \eqref{ckg:same-B-infty}--\eqref{ckg:same-Rj-matrix}.  Since $\rho_i\in[\varrho/2,3\varrho/2]$, one has
\[
 \sup_{\rho\in\Xi_0}\|R_j^0(\rho)\| \le C\frac{\varrho}{\widetilde w_j}, \qquad \varrho\sup_{\rho\ne\rho'} \frac{\|R_j^0(\rho)-R_j^0(\rho')\|} {|\rho-\rho'|} \le C\frac{\varrho}{\widetilde w_j}.
\]
The same estimate holds for the direct sums produced by the two Fourier signs in a full block $B=[j]$.  This proves \eqref{ckg:intrinsic-tail}.

The formula \eqref{ckg:same-beta-infty} is the spectrum of the diagonal matrix \eqref{ckg:same-B-infty}.  Finally, \eqref{ckg:intrinsic-tail-eigen-pointwise} follows from Lemma~\ref{ckg:lem-Weyl-branches} at each fixed parameter.
\end{proof}

We next separate the first-step tangential and elliptic normal-form corrections from those of the subsequent Newton iteration. We call $(\widehat\omega_0^{\rm KG},\widetilde A_0^{\rm KG})$ the \emph{bare reference} frequencies and incorporate the first correction into a fixed \emph{reference}:

\begin{equation}\label{ckg:frozen-reference-def}
 \widehat\omega^{\rm ref}:=\widehat\omega_0^{\rm KG}+\delta^{\rm ref}, \qquad \widetilde A_B^{\rm ref}:=\widetilde A_{0,B}^{\rm KG}+D_B^{\rm ref}.
\end{equation}
At every subsequent Newton step we write
\begin{equation}\label{ckg:frozen-decomposition}
 \widehat\omega=\widehat\omega^{\rm ref}+\delta, \qquad \widetilde A_B=\widetilde A_B^{\rm ref}+D_B.
\end{equation}

The first correction has size $O(\epsilon_0)$ and raw Lipschitz seminorm $O(\epsilon_0/\alpha_0)=o(1)$, where $\alpha_0$ is the preliminary Melnikov threshold; for $m\ne n$, it must be kept separate from the smaller restarted threshold $\alpha_1$, introduced in Section~\ref{sec:CKG-KAM-theorem}.

Let $\beta_{\rm tail}\in[0,1]$ be the available normal decay exponent, and let $\varepsilon_{\rm ref},\kappa_{\rm ref}\ge0$ be bounds for the frozen corrections \eqref{ckg:frozen-reference-def} as follows
\begin{equation}\label{ckg:frozen-tail-assumption}
 |\delta^{\rm ref}|_\infty +|D^{\rm ref}|_{\infty,\beta_{\rm tail}} \le C\varepsilon_{\rm ref}, \qquad |\delta^{\rm ref}|_\infty^{\rm lip} +|D^{\rm ref}|_{\infty,\beta_{\rm tail}}^{\rm lip} \le \kappa_{\rm ref}.
\end{equation}
For the subsequent Newton corrections, we assume that, for some $\varepsilon\ge0$,
\begin{equation}\label{ckg:iter-tail-assumption}
 |\delta|^{\mathcal L(\alpha)} +|D|_{\infty,\beta_{\rm tail}}^{\mathcal L(\alpha)} \le C\varepsilon.
\end{equation}
At the preliminary step $\beta_{\rm tail}=\ckgcasebeta$ by \eqref{ckg:case-beta-definition}, while at subsequent Newton steps $\beta_{\rm tail}=\ckgnewtonbeta$ by \eqref{ckg:newton-beta-definition}. In the distinct-mode counting argument we set $\theta_0=1-\beta_{\rm tail}$.

Before the first step one has $\varepsilon_{\rm ref}=\kappa_{\rm ref}=0$. After the preliminary step, recalling \eqref{ckg:E0-full}, we may take $\varepsilon_{\rm ref}=\epsilon_0$ and $\kappa_{\rm ref}=C\epsilon_0/\alpha_0=o(1)$, although $\epsilon_0/\alpha_1$ need not be small. 

\subsection{Two-dimensional transversality}

Here we prove the standard transversality result and the measure estimate of the single bad set \eqref{ckg:branch-resonant-strip}.

The tangential twist is the Jacobian of the bare tangential frequencies \eqref{ckg:initial-shifted-coefficients} with respect to $\rho$. In the distinct-mode case, \eqref{ckg:omega-dist} gives
\begin{equation*}
 A_h:=D_\rho\widehat\omega_0^{\rm KG}= \frac1{4\pi}
 \begin{pmatrix}
 \widetilde w_m^{-2}&2(\widetilde w_m\widetilde w_n)^{-1}\\
 2(\widetilde w_m\widetilde w_n)^{-1}&\widetilde w_n^{-2}
 \end{pmatrix}.
\end{equation*}
For a generic elliptic block $B=[j]$, the bare eigenvalues all coincide with $\mu_j^0:=\Lambda_j$ by \eqref{ckg:Lambda-dist}, and hence
\begin{equation*}
 b_j:=D_\rho\mu_j^0= \frac1{2\pi\widetilde w_j} \left(\widetilde w_m^{-1},\widetilde w_n^{-1}\right).
\end{equation*}
By \eqref{ckg:twist-close}, $A_h=A_*+O(h)$; hence, after reducing $h_*$,
\begin{equation}\label{ckg:twist-uniform-inverse}
 \sigma_{\min}(A_h)\ge c_*>0, \qquad \|(A_h^{t})^{-1}\|\le C_*.
\end{equation}
A direct inversion gives the useful identity
\begin{equation}\label{ckg:Ainv-b}
 (A_h^{t})^{-1}b_j= \left(\frac{2\widetilde w_m}{3\widetilde w_j}, \frac{2\widetilde w_n}{3\widetilde w_j}\right).
\end{equation}
For $h=0$ this reduces to $(2/3,2/3)$; the issue below is not a degeneracy of $A_h$, but the possible arithmetic cancellation between the tangential and normal derivatives.

\begin{lemma}
\label{ckg:lem-two-parameter-strip}
There exist $h_*,\varrho_*,\alpha_*,\kappa_*,C>0$, depending only on the fixed tangential sites, such that for $0<h<h_*$, $0<\varrho<\varrho_*$ and
$0<\alpha<\alpha_*$ the following holds. Consider a symmetry-allowed elliptic branch divisor with $|\ell|_1\le2$ which is not retained in the normal form.  Assume \eqref{ckg:frozen-tail-assumption}-\eqref{ckg:iter-tail-assumption} with $\varepsilon_{\rm ref}+\kappa_{\rm ref}+\varepsilon/\alpha<\kappa_*$.  Then either
\begin{equation}\label{ckg:strip-empty-gap}
 |\Phi(\rho)|\ge\kappa_* \quad\text{on }\Xi_0,
\end{equation}
or there is a unit vector $v=v(k,\ell)$ such that, on every line segment contained in $\Xi_0$,
\begin{equation}\label{ckg:strip-transverse}
 |\Phi(\rho+tv)-\Phi(\rho+t'v)| \ge \kappa_*\langle k\rangle |t-t'|.
\end{equation}
Consequently, for the strips defined in \eqref{ckg:branch-resonant-strip},
\begin{equation}\label{ckg:single-strip-measure}
 |\mathcal R_\Phi^{\theta_1}(\alpha)| \le C\varrho\, \frac{\alpha} {\langle k\rangle^{\tau+1} \widetilde w(\Phi)^{\theta_1}}.
\end{equation}
\end{lemma}
\begin{proof}
We first treat the bare reference normal forms, and then add the frozen and Newton corrections. Throughout $0<h<h_*$; since the tangential sites are fixed, $\widetilde w_m,\widetilde w_n=1+O(h)$.

\smallskip
\noindent\emph{Distinct tangential modes.}
For a bare branch $\Phi^0=\langle \omega_0,k\rangle+\sum_j\ell_j\mu_j^0$ as in \eqref{ckg:branch-divisors},
\begin{equation}\label{ckg:reference-gradient-dist}
 D_\rho\Phi^0=A_h^{t}k+\sum_j\ell_jb_j.
\end{equation}
If $|D_\rho\Phi^0|\ge\kappa_*\langle k\rangle$, transversality already holds. Otherwise \eqref{ckg:twist-uniform-inverse} and \eqref{ckg:Ainv-b} imply
\begin{equation}\label{ckg:trans-cancel-dist}
 k=-\sum_j\ell_j\left(\frac{2\widetilde w_m}{3\widetilde w_j},\frac{2\widetilde w_n}{3\widetilde w_j}\right)+O(\kappa_*\langle k\rangle).
\end{equation}
Since $|\ell|_1\le2$ and $\widetilde w_j\ge1$, the sum is uniformly bounded. For sufficiently small $\kappa_*$, \eqref{ckg:trans-cancel-dist} therefore forces $|k|_1\le K_0$, independently of $c$ and of the normal indices. Thus only finitely many $k$ need further consideration. Put $S_h:=\sum_j\ell_j/\widetilde w_j$. For the remaining bounded $k$, \eqref{ckg:trans-cancel-dist} and $\widetilde w_m-\widetilde w_n=O(h)$ give $k_1-k_2=O(h+\kappa_*)$. Since $k_i \in \mathbb{Z}$, for $h$ and $\kappa_*$ small enough, $k_1=k_2=:r$. Since $|S_h|\le2$, the same relation yields $|r|\le4/3+O(h+\kappa_*)$, hence $r\in\{-1,0,1\}$. For $r=0$, a two-normal difference is retained by momentum \eqref{ckg:momentum-kell}, whereas an unretained same-sign pair has $i+j=0$ and a uniform frequency gap $|\Phi^0| \geq 2c^2-O(\varrho)$. For $|r|=1$, near cancellation requires $|S_h|=3/2+O(h+\kappa_*)$. This is impossible with one normal factor, or with opposite normal signs since $|\widetilde w_i^{-1}-\widetilde w_j^{-1}|\le1$. Thus, up to reversing all Melnikov signs, the only possibility is $k=(1,1)$ with two negative normal signs, for which
\begin{equation}\label{ckg:harmonic-relation}
 \frac1{\widetilde w_i}+\frac1{\widetilde w_j}=\frac32+O(h+\kappa_*).
\end{equation}
The arithmetic--harmonic mean inequality then gives $\widetilde w_i+\widetilde w_j\ge 4/(\widetilde w_i^{-1}+\widetilde w_j^{-1})\ge8/3-O(h+\kappa_*)$. Since $\Phi^0=c^2(\widetilde w_m+\widetilde w_n-\widetilde w_i-\widetilde w_j)+ O(\varrho)$ and $\widetilde w_m+\widetilde w_n=2+O(h)$, the leading divisor has modulus at least $c^2/3$. The $O(\varrho)$ quartic corrections preserve this gap. This proves the bare-reference alternative for distinct tangential modes.

\smallskip
\noindent\emph{Identical tangential mode.}
For locally gauge invariant coefficients Lemma~\ref{ckg:lem-selection-same} gives $k=0$. Momentum \eqref{ckg:momentum-kell-same} places every two-normal difference in a retained block; an unretained same-sign pair has indices $m+d,m-d$ by \eqref{ckg:momentum-kell-same}, with $d\ne0$, and Lemma~\ref{ckg:lem-local-large-gap} yields
\[
 \lambda_{m+d}+\lambda_{m-d}-2\lambda_m\ge \kappa_m,
\]
for some $\kappa_m>0$. Thus these coefficients have a uniform gap. In the non-locally-gauge-invariant sector, Lemma~\ref{ckg:lem-selection-same} gives $k=(r,r)$, $r\ne0$. Along $e=(1,1)$, the formulas \eqref{ckg:omega-same}--\eqref{ckg:coeff-same} give
\begin{equation}\label{ckg:same-radial-tangential}
 \partial_e\langle\omega^{(m,m)},k\rangle=\frac{3r}{2\pi\widetilde w_m^2},
\end{equation}
\begin{equation}\label{ckg:same-radial-block}
 \partial_eM_j=a_jI_2+g_j\begin{pmatrix}0&1\\1&0\end{pmatrix},
\end{equation}
where
\begin{equation}\label{ckg:same-radial-ag}
 a_j=\frac1{\pi\widetilde w_m\widetilde w_j}-\frac3{4\pi\widetilde w_m^2},\qquad |g_j|\le\frac1{\sqrt3\,\pi\widetilde w_m\widetilde w_j}.
\end{equation}
Indeed $g_j=(4\pi\widetilde w_m\widetilde w_j)^{-1}(\sqrt{\rho_2/\rho_1}+\sqrt{\rho_1/\rho_2})$ and $\rho_1/\rho_2\in[1/3,3]$. Thus $\|\partial_eM_j\|\le C$ uniformly in $j,c,\rho$. The derivative of every self-adjoint one- or two-normal coefficient operator differs from \eqref{ckg:same-radial-tangential} by at most $2C$ in operator norm. The min--max principle therefore gives transversality of order $|r|$ for all sufficiently large $|r|$, including across eigenvalue crossings.

Only $0<|r|\le R$ remain, with $R$ independent of $c$ and the normal indices. Reversing all Melnikov signs, assume $1\le r\le R$. By \eqref{ckg:coeff-same}, each bare normal eigenvalue equals $\lambda_j-\lambda_m+O(\varrho)$ uniformly. We have $ \langle \omega, r \rangle= 2r \lambda_m+ O(\varrho)$. Without a normal factor in $\Phi^0$, there is a gap of order $c^2$. For one normal factor, momentum \eqref{ckg:momentum-kell-same} fixes its index in terms of $m,r$, giving $\lambda_j-\lambda_m=O(1)$ uniformly in $c$, which cannot cancel the tangential gap of order $c^2$. For a two-normal difference, momentum yields $|i-j|=2|mr|$, up to exchanging indices, and
\[
 |\lambda_i-\lambda_j|\le c\bigl||i|-|j|\bigr|\le c|i-j|\le2c|mr|,
\]
which again cannot compensate the gap of order $c^2$ for $c$ big enough. Two positive normal signs likewise give the lower bound $2(r-1)\lambda_m+\lambda_i+\lambda_j\ge2rc^2$ for the leading part. The $O_R(\varrho)$ shifts preserve these gaps.

It remains to consider two negative normal signs. Fix $i,j$ and an ordered branch $\Phi^0$ of the bare quadratic divisor. If $|\Phi^0(\rho)|\ge1$ throughout $\Xi_0$, the branch has a gap. Otherwise, for some $\rho^0\in\Xi_0$, $|\Phi^0(\rho^0)|<1$. Since the normal eigenvalues differ from $\lambda_j-\lambda_m$ by at most $C\varrho$, recalling $\tilde w_j=\lambda_j/c^2$ this implies
\begin{equation}\label{ckg:same-near-resonance-weights}
 \left|\widetilde w_i+\widetilde w_j-2(r+1)\widetilde w_m\right|\le\frac{C_R(1+\varrho)}{c^2}.
\end{equation}
Momentum \eqref{ckg:momentum-kell-same} gives $i+j=2m(r+1)$. Since $x\mapsto\sqrt{1+x^2/c^2}$ is Lipschitz with constant of order $c^{-1}$, we have
\[
 |\widetilde w_i-\widetilde w_j|\le\frac{\bigl||i|-|j|\bigr|}{c}\le\frac{|i+j|}{c}\le\frac{2|m|(r+1)}c.
\]
Combining the sum and difference estimates gives
\begin{equation}\label{ckg:same-resonance-each-weight}
 \widetilde w_i=(r+1)\widetilde w_m+O_R(c^{-1}),\qquad \widetilde w_j=(r+1)\widetilde w_m+O_R(c^{-1}).
\end{equation}
These estimates are uniform in $i,j,\rho$ and hold throughout $\Xi_0$, since the weights do not depend on $\rho$.

We now establish the required monotonicity directly on the Hermitian Sylvester operator, without differentiating its individual eigenvalues. The eigenvalues of $\partial_eM_j$ in \eqref{ckg:same-radial-block} are $a_j\pm g_j$; hence by \eqref{ckg:same-radial-ag}, using the Loewner order introduced in the Lemma \ref{operatori:hem}, we have
\[
 \partial_eM_j\le\left(\frac{1+1/\sqrt3}{\pi\widetilde w_m\widetilde w_j}-\frac3{4\pi\widetilde w_m^2}\right)I_2.
\]
On the Hilbert--Schmidt space of $2\times2$ coefficient matrices, the two-negative-sign divisor is the self-adjoint Sylvester operator $\mathscr H=-i \mathscr L$ (recall \eqref{ckg:Lplus}), with
\[
 \mathscr H^0_{r,i,j}X:=r(\omega_1^{(m,m)}+\omega_2^{(m,m)})X-M_iX-XM_j^{t}.
\]
Its eigenvalues are the corresponding scalar branch divisors by Lemma~\ref{ckg:lem-spectral-reduction}. The preceding matrix inequality and \eqref{ckg:same-radial-tangential} give, for every coefficient matrix $X$,
\[
 \begin{aligned}
 \langle X,(\partial_e\mathscr H^0_{r,i,j})X\rangle_{\rm HS}
 &\ge\frac1{\pi\widetilde w_m^2}\left[\frac{3(r+1)}2-\left(1+\frac1{\sqrt3}\right)\widetilde w_m\left(\frac1{\widetilde w_i}+\frac1{\widetilde w_j}\right)\right]\|X\|_{\rm HS}^2\\
 &\ge\frac1{\pi\widetilde w_m^2}\left[\frac{3(r+1)}2-\frac{2(1+1/\sqrt3)}{r+1}-O_R(c^{-1})\right]\|X\|_{\rm HS}^2,
 \end{aligned}
\]
where the last inequality uses \eqref{ckg:same-resonance-each-weight}. For every integer $r\ge1$,
\[
 \frac{3(r+1)}2-\frac{2(1+1/\sqrt3)}{r+1}>\frac{r+1}{2},
\]
since $(r+1)^2\ge4>2(1+1/\sqrt3)$. For $c$ sufficiently large, depending only on $m,R$, for some $c_0 ,c_1>0$, this proves
\begin{equation}\label{ckg:same-uniform-Loewner}
 \partial_e\mathscr H^0_{r,i,j}\ge c_0(r+1)I\ge c_1\langle k\rangle I,\qquad k=(r,r),\quad\rho\in\Xi_0,
\end{equation}
where $I$ is the identity on the Hilbert--Schmidt coefficient space. The bound also holds on the invariant symmetry-allowed subspaces. Integration along segments parallel to the unit vector $e/|e|$ and the min--max principle give uniform monotonicity of each ordered branch, including through eigenvalue crossings. This proves the bare-reference alternative in the identical-mode case.

\medskip

It remains to add the frozen and Newton corrections. For each symmetry-allowed coefficient, write the full self-adjoint coefficient operator as $\mathscr H=\mathscr H^0+E$, where $E$ consists of $\langle\delta^{\rm ref}+\delta,k\rangle I$ and at most two signed normal matrix corrections from $D^{\rm ref},D$. In the gapped cases the relevant $k$ are bounded and the bare gap is uniform (or of order $c^2$); hence the supremum bounds in \eqref{ckg:frozen-tail-assumption}--\eqref{ckg:iter-tail-assumption} preserve \eqref{ckg:strip-empty-gap}. In the remaining cases, choose the orientation of a unit direction $v$ so that $\mathscr H^0(\rho+tv)-\mathscr H^0(\rho+t'v)\ge m_0\langle k\rangle(t-t')I$ for $t>t'$, and for some $m_0>0$. The Lipschitz bounds in \eqref{ckg:frozen-tail-assumption}--\eqref{ckg:iter-tail-assumption}, on the same segments, gives
\[
 \|E(\rho+tv)-E(\rho+t'v)\|\le C(\kappa_{\rm ref}+\varepsilon/\alpha)\langle k\rangle|t-t'|.
\]
Choosing the smallness constants so that this slope is below $m_0\langle k\rangle/2$, the operator-monotonicity lemma in the Appendix proves \eqref{ckg:strip-transverse} for every ordered branch, without requiring Lipschitz estimates for differences of individual eigenvalues.

Finally put $\delta_k:=\alpha/(\langle k\rangle^\tau\widetilde w(\Phi)^{\theta_1})$. By the same lemma, on each line parallel to $v$ the section of $\{|\Phi|<2\delta_k\}$ has length at most $C\delta_k/\langle k\rangle$. The orthogonal projection of the parameter square $\Xi_0$ has length at most $\sqrt2\,\varrho$, so Fubini gives \eqref{ckg:single-strip-measure}. The same operator argument applies to the large-index limiting matrices used below.
\end{proof}

The single-strip estimate still has to be summed over the normal labels.  In the locally gauge invariant sector this multiplicity is finite, as recorded next.

By Lemmas~\ref{ckg:lem-selection-same}, \ref{ckg:lem-selection-dist} and \ref{ckg:lem-local-large-gap}, together with the momentum identities \eqref{ckg:momentum-kell}--\eqref{ckg:momentum-kell-same}, we obtain the following finite counting result.
\begin{lemma}
\label{ckg:lem-local-branch-count}
For $0<\alpha\le c_\alpha\varrho$ and sufficiently small $\varrho$, for every $k\in\mathbb Z^2$ the number of locally gauge invariant elliptic branch strips which are not retained in $[P]$ and can intersect $\Xi_0$ is bounded by an absolute constant.  In particular it is bounded by $C\langle k\rangle^2$.  In the distinct-mode case such strips can occur only for $|k|_1\le4$; in the identical-mode case every locally gauge invariant strip is either retained or empty.  The same conclusions hold for the NLS reference normal form.
\end{lemma}

We can now prove the anticipated strong result about the NLS bad sets.

\begin{lemma}\label{ckg:lem-NLS-no-bad-strips}
For both choices of tangential sites, at every NLS KAM step satisfying \eqref{ckg:freq-corrections}, the smallness condition \eqref{ckg:homological-smallness}, and (when present) the NLS finite-resolvent neighbourhood of Lemma~\ref{ckg:lem-hyp-resolvent}, all symmetry-allowed elliptic branch divisors not retained by $[\cdot]$ have absolute value at least $1/2$ on $\Xi_0$. In particular, for $0<\alpha\le c_\alpha\varrho$ and sufficiently small $\varrho$, $\mathcal R_\Phi^{\rm NLS}(\alpha)=\varnothing$ for every symmetry-allowed branch. Coefficients with finite factors require no parameter exclusion.
\end{lemma}
\begin{proof}
The initial NLS Hamiltonian is invariant under the single mass flows. This property is preserved at each KAM step: the normal projection preserves local gauge invariance by \eqref{ckg:normal-projection-LG}, and the homological operator commutes with $\Pi_{\rm LG}$ by \eqref{ckg:block-fast-commute}. Thus the normalized NLS homological solution and the subsequent transformed perturbation remain locally gauge invariant.

Suppose first that $m=n$. Lemma~\ref{ckg:lem-selection-same} gives $k=0$ for every locally gauge invariant coefficient. Momentum excludes a single elliptic factor and forces the two indices in an elliptic difference to coincide; all such differences are retained by Definition~\ref{ckg:def-normal-projection}.
An unretained same-sign quadratic coefficient has spatial indices $m+d$ and $m-d$, with $d\in\mathbb Z\setminus\{0\}$. By \eqref{ckg:coeff-same-NLS}, its divisor after dropping the action-dependent quartic shifts is
\[
 \frac{(m+d)^2-m^2}{2}+
 \frac{(m-d)^2-m^2}{2}=d^2\ge1.
\]

Suppose now that $m\ne n$. For a locally gauge invariant coefficient with no finite factors, the two mass constraints and momentum \eqref{ckg:sel-dist-N}--\eqref{ckg:sel-dist-M} force every zero-normal coefficient to have $k=0$, so it is retained; they exclude all one-normal elliptic coefficients. For a same-species quadratic difference they give $k=0$ and equal spatial indices, so the coefficient is again retained. For the remaining quadratic coefficients, let $D_{\rm lin}$ denote the divisor obtained from the quadratic NLS frequencies only, that is, by omitting the action-dependent shifts in \eqref{ckg:dist-NLS-coeff}:
\[
 D_{\rm lin}:=\frac12\left(m^2k_1+n^2k_2+ \sum_{(j,s)\in\mathcal E}j^2\ell_{j,s}\right).
\]
Up to complex conjugation, the full list, with the Fourier index $k$ fixed by the two mass constraints, is
\[
\begin{array}{c|c|c|c}
\text{elliptic factors}&k&\text{momentum}&D_{\rm lin}\\\hline
 z_{(i,+1)}z_{(j,+1)}&(-2,0)&i+j=2m&(i-m)^2\\
 z_{(i,-1)}z_{(j,-1)}&(0,-2)&i+j=2n&(i-n)^2\\
 z_{(i,+1)}z_{(j,-1)}&(-1,-1)&i+j=m+n&(i-m)(i-n)\\
 z_{(i,+1)}\bar z_{(j,-1)}&(-1,1)&i-j=m-n&(m-n)(i-m).
\end{array}
\]
Here $i,j\notin\{m,n\}$ because these are generic elliptic variables. Consequently every displayed $D_{\rm lin}$ is a nonzero integer. 

In either configuration, denote by $\mathscr H_{\rm lin}$ the elliptic coefficient operator formed from the quadratic frequencies above and restricted to the symmetry-allowed, nonretained coefficient subspace. It is self-adjoint and all its eigenvalues have absolute value at least $1$. The explicit NLS normal forms \eqref{ckg:coeff-same-NLS}, \eqref{ckg:dist-NLS-coeff} and the correction estimate \eqref{ckg:freq-corrections} give, uniformly in the elliptic indices,
\[
 \|\mathscr H^{\rm NLS}-\mathscr H_{\rm lin}\| \le C(\varrho+\varepsilon_{\rm NF}).
\]
Indeed, $k=0$ for $m=n$, while $|k|_1\le4$ for $m\ne n$ by \eqref{ckg:k-local-bound}. For sufficiently small $\varrho+\varepsilon_{\rm NF}$, Weyl's inequality then yields $\operatorname{dist}\bigl(0,\operatorname{spec} \mathscr H^{\rm NLS}\bigr)\ge 1/2$. The argument applies on the current parameter set; for the extensions to $\Xi_0$ used in the measure construction, apply Lemma~\ref{ckg:lem-extension} to the symmetry-allowed matrix entries and keep the forbidden entries zero. The extension changes the bounds only by an absolute factor, so the preceding gap persists after reducing the smallness constants. For $\alpha\le c_\alpha\varrho$ and sufficiently small $\varrho$, we have $2\alpha<1/2$; hence every corresponding NLS branch strip \eqref{ckg:branch-resonant-strip} is empty. Finally, coefficients with exactly one finite factor require no parameter exclusion by Corollary~\ref{ckg:cor-one-finite-factor}. In the distinct-mode case, any locally gauge invariant coefficient with two finite factors has $k=0$ by Lemma~\ref{ckg:lem-selection-dist} and is retained by Definition~\ref{ckg:def-normal-projection}.
\end{proof}

\subsection{The identical-mode case}

We first complete the measure estimate for $m=n$. In this case the selection rules reduce the potentially resonant families to a direct unweighted count, without the singular large-index analysis required for distinct sites.

\begin{lemma}\label{ckg:lem-same-unweighted-measure}
Fix $m=n$ and assume \eqref{ckg:freq-corrections}, with $\varepsilon_{\rm NF}+\kappa_{\rm NF}$ sufficiently small. For $c\ge c_m$, $0<\alpha<c_\alpha\varrho$ and $\tau>1$,
\begin{equation}\label{ckg:same-unweighted-measure}
 \frac1{|\Xi_0|}\left|\bigcup_k\mathcal R_k^0(\alpha)\right|\le C\frac\alpha\varrho,\qquad
 \frac1{|\Xi_0|}\left|\bigcup_{|k|_1>K}\mathcal R_k^0(\alpha)\right|\le C\frac\alpha\varrho(1+K)^{1-\tau},\quad K\ge0.
\end{equation}
Here $\mathcal R_k^0$ is the union of the symmetry-allowed branch strips \eqref{ckg:branch-resonant-strip} at fixed $k$ not retained by $[ \cdot ]$, excluding coefficients containing a finite hyperbolic factor. 
\end{lemma}

\begin{proof}
By the computations in Lemma~\ref{ckg:lem-local-branch-count} for the locally gauge sector, and by the computations in Lemma~\ref{ckg:lem-selection-same} for the non-locally gauge sector, the only potentially resonant family has, up to conjugation and finite species multiplicity, the form
\[
 \Phi=2(r+1)\lambda_m-\lambda_i-\lambda_j+O((r+1)(\varrho+\varepsilon_{\rm NF})),\qquad i+j=2m(r+1).
\]
Let $g(x)=\sqrt{1+x^2}-1$. If its strip is nonempty, $\lambda_j=c^2(1+g(|j|/c))$ and momentum imply
\[
 |g(|i|/c)+g(|j|/c)-2r|\le C(r+1)c^{-2},\qquad \big||i|-|j|\big|\le2|m|(r+1).
\]
Since $g$ is Lipschitz, these inequalities give $|g(|i|/c)-r|\le C(r+1)/c$. For $c\ge c_m$ sufficiently large this implies $g(|i|/c)\ge1/2$, uniformly for $r\ge1$. The inverse
$g^{-1}(u)=\sqrt{u(u+2)}$ has bounded derivative on $[1/2,\infty)$, so
\[
 \big||i|-c\sqrt{r(r+2)}\big|\le C(r+1).
\]
Thus there are at most $C(r+1)$ integer choices for $i$ (including both spatial signs), and momentum determines $j$.  The identical-mode argument in the proof of Lemma~\ref{ckg:lem-two-parameter-strip} applies under the present unweighted bounds: by \eqref{ckg:freq-corrections}, the corrections have size $O((r+1)\varepsilon_{\rm NF})$ and Lipschitz constant $O((r+1)\kappa_{\rm NF})$. Hence the gap/transversality alternative persists for $\varepsilon_{\rm NF}+\kappa_{\rm NF}$ sufficiently small. It follows that \eqref{ckg:single-strip-measure} bounds each strip by $C\varrho\alpha/(r+1)^{\tau+1}$. Summing over the $O(r+1)$ admissible pairs, with $r\ge1$ or $r>K/2$, proves \eqref{ckg:same-unweighted-measure}.
\end{proof}

\subsection{The distinct-mode case: singular large-index counting}

We now assume $m\ne n$ and consider the non-LG elliptic coefficient operators $\mathscr H_{k\ell}:=-\iu\mathscr L_{k\ell}$ from \eqref{ckg:L0}--\eqref{ckg:Lplus}, restricted to their invariant coefficient spaces. Here $\ell$ is a representative signed normal exponent vector, and the fast defect in
\eqref{ckg:fast-defect} $L=\mathfrak L(k,\ell)\ne0$ is constant on each such space, since $q_B$ is constant on each normal block. Define
\begin{equation}\label{ckg:matrix-resonant-set}
\mathcal R_{k\ell}^{\theta_1}(\alpha) := \left\{\rho\in\Xi_0: \operatorname{dist} \bigl(0,\operatorname{spec}\mathscr H_{k\ell}(\rho)\bigr) < \frac{2\alpha} {\langle k\rangle^\tau \widetilde w_{k\ell}^{\theta_1}} \right\},
\end{equation}
where $\widetilde w_{k\ell}$ is the minimum normal-block weight, or $1$ if no normal block occurs.

By Lemma~\ref{ckg:lem-spectral-reduction}, these matrix strips are exactly the unions of the corresponding branch strips \eqref{ckg:branch-resonant-strip}. We denote by $\mathcal R_k^{\theta_1}(\alpha)$ the union of all symmetry-allowed, non-retained elliptic branch strips at fixed $k$, including the LG ones, controlled by Lemma~\ref{ckg:lem-local-branch-count}.

Put $\nu_{\rm int}(h):=(\nu_m(h),\nu_n(h))$. For a coefficient with $n_\ell\in\{0,1,2\}$ elliptic normal factors, counted with multiplicity, denote their spatial Fourier indices by $j_1,\ldots,j_{n_\ell}$ and their Melnikov signs by $\epsilon_s\in\{\pm1\}$. By Lemma~\ref{ckg:lem-intrinsic-tail} and \eqref{ckg:frozen-tail-assumption}--\eqref{ckg:iter-tail-assumption}, extraction of the rest-energy part gives
\begin{equation}\label{ckg:scalar-NG-form}
 \mathscr H_{k\ell} =c^2L I+\langle\nu_{\rm int}(h),k\rangle I +\sum_{s=1}^{n_\ell}\epsilon_s\nu_{j_s}(h)I +\mathscr B_{k\ell}^{\infty}(\rho;c) +\mathscr E_{k\ell}(\rho;c), \qquad n_\ell\le2,
\end{equation}
where $\mathscr B_{k\ell}^{\infty}$ contains the normal limiting matrices, slow quartic shifts and tangential corrections, all independent of the large normal labels. With $E_*:=\varepsilon_{\rm ref}+\varepsilon$, it holds
\begin{equation}\label{ckg:Binfty-bound}
 \|\mathscr B_{k\ell}^{\infty}\|\le C\bigl(1+(\varrho+E_*)\langle k\rangle\bigr).
\end{equation}
Let $J_{\min}:=\min_{1\le s\le n_\ell}|j_s|$. If $J_{\min}$ is sufficiently large, \eqref{ckg:intrinsic-tail} and \eqref{ckg:frozen-tail-assumption}--\eqref{ckg:iter-tail-assumption} yield
\begin{equation}\label{ckg:large-tail-summary}
\|\mathscr E_{k\ell}\| \le C\left( \frac{\varrho}{\widetilde w_{J_{\min}}} +\frac{E_*}{\widetilde w_{J_{\min}}^{1-\theta_0}} \right).
\end{equation}
Lipschitz control of these tails is always taken at matrix level, as in Lemma~\ref{ckg:lem-intrinsic-tail}; Weyl's inequality \eqref{ckg:Weyl-pointwise-ref} is used only pointwise. The following
localization is the block counterpart of \cite[Lemmas~5.8--5.9]{BBG}.

\begin{lemma}
\label{ckg:lem-NG-difference-localization}
There exist $c_0,C_0>0$ such that the following holds for $c\ge c_0$. Consider a non-locally-gauge-invariant two-normal difference coefficient. If its resonant set is nonempty, then
\begin{equation}\label{ckg:c-less-k}
 c\le C_0\langle k\rangle.
\end{equation}
If the two spatial Fourier indices have the same sign, then
\begin{equation}\label{ckg:d-less-k}
 d:=\big||j|-|i|\big|\le C_0\langle k\rangle.
\end{equation}
If they have opposite signs, or if $\min\{|i|,|j|\}\le\frac12\max\{|i|,|j|\}$, then nonemptiness already implies
\begin{equation}\label{ckg:easy-difference-indices}
 |i|+|j|\le C_0\langle k\rangle.
\end{equation}
\end{lemma}

\begin{proof}
The bounds \eqref{ckg:d-less-k}--\eqref{ckg:easy-difference-indices} follow directly from momentum \eqref{ckg:momentum-kell}, as in \cite[Lemma~5.10]{BBG}. The argument of \cite[Lemma~5.9]{BBG} also applies to the block decomposition \eqref{ckg:scalar-NG-form}: the matrix corrections are uniformly controlled, and Weyl's inequality gives the same $c^2$ gap for $\langle k\rangle\ll c$. This yields \eqref{ckg:c-less-k}.
\end{proof}

The following proposition contains the only singular counting argument.  The three difference regimes correspond to the proofs of \cite[Lemmas~5.10--5.12]{BBG}, summarized in \cite[Lemma~5.13]{BBG}; the
three same-sign regimes correspond to \cite[Lemmas~5.14--5.16]{BBG}, summarized in \cite[Lemma~5.17]{BBG}.  The changes here are the finite block multiplicity and the two matrix tails in
\eqref{ckg:large-tail-summary}.

\begin{proposition}
\label{ckg:prop-singular-branch-measure}
Assume $m \neq n$. Let $0\le\theta_0<1$ and $0\le\theta_1<1-\theta_0$, and set
\begin{equation*}
 a_f:=\frac{1-\theta_1}{2-\theta_1-\theta_0}.
\end{equation*}
Define
\begin{equation}\label{ckg:tau-sharp}
 \tau^\sharp:= \min\left\{ \frac{\tau(1-\theta_0)}{3-\theta_1}, \frac{\tau}{2-\theta_1}, \frac{\tau(1-\theta_0)}{2-\theta_1-\theta_0} \right\}-6.
\end{equation}
Assume $\tau^\sharp>0$.  There exist $h_*,\varrho_*,\alpha_*,\kappa_*,C>0$, with $C$ independent of $c,\varrho,\alpha,k$, such that, under \eqref{ckg:frozen-tail-assumption}--\eqref{ckg:iter-tail-assumption} and
\begin{equation*}
 \varepsilon_{\rm ref}+\kappa_{\rm ref}+\varepsilon/\alpha <\kappa_*,
\end{equation*}
for $0<h<h_*$, $0<\varrho<\varrho_*$, $0<\alpha<\min\{\alpha_*,c_\alpha\varrho\}$ (where $c_\alpha$ is as in \eqref{ckg:homological-smallness}) and every $k\in\mathbb Z^2$,
\begin{equation}\label{ckg:NG-union-measure}
 \left|\bigcup_{\ell:\,\mathfrak L(k,\ell)\ne0} \mathcal R_{k\ell}^{\theta_1}(\alpha)\right| \le \frac{C\varrho}{\langle k\rangle^{\tau^\sharp}} \left[ \alpha^{\frac{2}{3-\theta_1}} +\varrho\left(\frac{\alpha}{\varrho}\right)^{\frac1{2-\theta_1}} +E_*^{a_f}\alpha^{1-a_f} \right].
\end{equation}
Moreover, there is $\kappa_0>0$ such that if $|k|_1\le\kappa_0c^{1/2}$ and every elliptic index $j$ in the coefficient satisfies $|j|<c/2$, then
\begin{equation}\label{ckg:NG-c2-lower}
 \operatorname{dist}\bigl(0,\operatorname{spec}\mathscr H_{k\ell}(\rho)\bigr) \ge\frac14c^2, \qquad \rho\in\Xi_0.
\end{equation}
\end{proposition}

\begin{proof}
By Lemma~\ref{ckg:lem-finite-multiplicity}, species and ordered spectral branches contribute only an absolute multiplicity.  We therefore count only the spatial Fourier configurations.

\smallskip
\noindent\emph{The low-index $c^2$ gap.}
If every elliptic index satisfies $|j|<c/2$, then
\[
 0\le\nu_j(h) =c^2\left(\sqrt{1+j^2/c^2}-1\right) \le c^2(\sqrt{5/4}-1)<\frac18c^2.
\]
There are at most two elliptic factors, hence their total slow dispersive contribution is $<c^2/4$.  The remaining slow part in \eqref{ckg:scalar-NG-form} is $O(1+\langle k\rangle+(\varrho+E_*)\langle
k\rangle)$. Since $|L|\ge1$ and $|k|_1\le\kappa_0c^{1/2}$, after reducing $\kappa_0,h_*,\varrho_*$ and the correction smallness constants,
\[
 \operatorname{dist}\bigl(0,\operatorname{spec}\mathscr H_{k\ell}\bigr) \ge c^2-\frac14c^2-o(c^2)\ge\frac14c^2.
\]
This proves \eqref{ckg:NG-c2-lower}; it is the block version of the low--low non-gauge hypothesis used in \cite[Lemma~4.2]{BBG}.

\smallskip
\noindent\emph{Zero or one effective normal index.} By momentum \eqref{ckg:momentum-kell}, coefficients with at most one normal factor, or with two coincident normal indices and equal Melnikov signs, give only $O(1)$ resonant strips for fixed $k$. The coincident-index difference case ($d=0$) is treated below. Thus \eqref{ckg:single-strip-measure} gives a total contribution $C\varrho\alpha/\langle k\rangle^{\tau+1}$, which is absorbed by \eqref{ckg:NG-union-measure}.

\smallskip
\noindent\emph{Two normal frequencies with opposite Melnikov signs.}
Assume $|i|\le|j|$.  If the spatial signs are opposite, or if $|i|\le|j|/2$, then \eqref{ckg:easy-difference-indices} leaves at most $C\langle k\rangle^2$ pairs.  The single-strip estimate gives
\begin{equation}\label{ckg:easy-difference-measure}
 C\varrho\frac{\alpha}{\langle k\rangle^{\tau-1}}.
\end{equation}

Assume now that the spatial signs agree and $|j|/2<|i|\le|j|$.  If $d:=|j|-|i|=0$, then $i=j$; the dispersive difference vanishes. Since the two generic elliptic blocks have equal fast charge,
$L=k_1+k_2$, while momentum gives $mk_1+nk_2=0$. Hence $|k|_1\le C_{m,n}|L|$, and \eqref{ckg:scalar-NG-form} implies
\[
 \operatorname{dist}(0,\operatorname{spec}\mathscr H_{k\ell})\ge c^2|L|-C\langle k\rangle-C\ge\tfrac12c^2
\]
for $c$ sufficiently large. Thus $d=0$ is nonresonant. For $d\ge1$, Lemma~\ref{ckg:lem-NG-difference-localization} gives
\begin{equation}\label{ckg:cd-localization}
 c\lesssim\langle k\rangle, \qquad 1\le d\lesssim\langle k\rangle.
\end{equation}
If $|j|\le2c^3$, then $|i|+|j|\lesssim c^3\lesssim\langle k\rangle^3$; a crude count therefore leaves at most $C\langle k\rangle^6$ pairs, whose total measure is
\begin{equation}\label{ckg:intermediate-difference-measure}
 C\varrho\frac{\alpha}{\langle k\rangle^{\tau-5}}.
\end{equation}
This is the same intermediate estimate as in \cite[Lemma~5.11]{BBG}.

It remains to treat $c^3\le|i|<|j|$ with comparable moduli.  Rationalization gives
\begin{equation*}
 \frac{\lambda_j-\lambda_i}{cd} =\frac{|i|+|j|} {\sqrt{c^2+i^2}+\sqrt{c^2+j^2}} =1+O(\widetilde w_i^{-2}).
\end{equation*}

For fixed $k,d$ and each of the finitely many spatial-sign and species types, let $\mathscr H_{k,d}^{\infty}$ be obtained by replacing $\lambda_j-\lambda_i$ by $cd$ (with the appropriate Melnikov sign) and discarding the decaying normal tails. By Lemma~\ref{ckg:lem-intrinsic-tail}, the generic distinct-mode blocks have $B_B^\infty=0$, so on the symmetry-allowed coefficient spaces
\[
\mathscr H_{k,d}^{\infty}(\rho) =\bigl(\langle\omega(\rho),k\rangle\pm cd\bigr)I.
\]
Moreover,
\begin{equation}\label{ckg:large-difference-matrix-error}
\|\mathscr H_{k\ell}-\mathscr H_{k,d}^{\infty}\| \le C\left( c|d|\widetilde w_i^{-2} +\varrho\widetilde w_i^{-1} +E_*\widetilde w_i^{-(1-\theta_0)} \right).
\end{equation}

Since $d\ge1$, momentum \eqref{ckg:momentum-kell} implies $k\ne0$. The limiting divisor is $\langle\omega(\rho),k\rangle\pm cd$. By \eqref{ckg:twist-uniform-inverse} and the smallness assumptions \eqref{ckg:frozen-tail-assumption}-- \eqref{ckg:iter-tail-assumption}, it has a directional slope bounded below by $C^{-1}\langle k\rangle$. The Fubini argument of Lemma~\ref{ckg:lem-two-parameter-strip} therefore bounds its strip of width $\delta$ by $C\varrho\delta/\langle k\rangle$.

Fix a cutoff $M>c^3$.  For $c^3\le|i|<M$ sum the single-strip estimates. For $|i|\ge M$, both the strip width and the right-hand side of \eqref{ckg:large-difference-matrix-error} decrease with $|i|$;
Weyl's inequality therefore puts all these resonant sets inside one enlarged resonant set for $\mathscr H_{k,d}^{\infty}$.  For a fixed asymptotic type with equal spatial signs and comparable moduli,
$\widetilde w_i\asymp |i|/c$ on $|i|\ge c^3$, and the corresponding union obeys
\begin{equation}\label{ckg:large-difference-union-I}
 \left|\bigcup_{\substack{|i|\ge c^3,\ |j|-|i|=d\\ \operatorname{sgn}(i)=\operatorname{sgn}(j)}} \mathcal R_{k\ell}^{\theta_1}(\alpha)\right| \le C\varrho\frac{c^3|d|}{\langle k\rangle} \left[ \frac{\alpha M^{1-\theta_1}}{\langle k\rangle^\tau} +M^{-2}+\varrho M^{-1} +E_*M^{-(1-\theta_0)} \right].
\end{equation}
The finite number of asymptotic types is absorbed into $C$. Indeed, after factoring out the common factor $C\varrho$ and before taking the common majorant $c^3|d|$, the four contributions are respectively
\[
 \frac{\alpha c^{\theta_1}M^{1-\theta_1}} {\langle k\rangle^{\tau+1}},\qquad \frac{c^3|d|}{\langle k\rangle M^2},\qquad \frac{\varrho c}{\langle k\rangle M},\qquad \frac{E_*c^{1-\theta_0}} {\langle k\rangle M^{1-\theta_0}}.
\]
Since $c\ge1$, $d\ge1$ and $\theta_0,\theta_1<1$, these are bounded by the right-hand side of \eqref{ckg:large-difference-union-I}.  This is precisely the proof of \cite[Lemma~5.12]{BBG}, with the
intrinsic matrix tail $\varrho\widetilde w_i^{-1}$ and the frozen tail $E_*\widetilde w_i^{-(1-\theta_0)}$ kept separately.

Balance the first term in brackets in \eqref{ckg:large-difference-union-I} with each of the remaining three:
\begin{align*}
 M_1\asymp \left(\frac{\langle k\rangle^\tau}{\alpha}
 \right)^{1/(3-\theta_1)},\qquad M_2\asymp \left(\frac{\varrho\langle k\rangle^\tau}{\alpha}
 \right)^{1/(2-\theta_1)},\qquad  M_3\asymp \left(\frac{E_*\langle k\rangle^\tau}{\alpha} \right)^{1/(2-\theta_1-\theta_0)}.
\end{align*}
Choosing $M\asymp\max\{M_1,M_2,M_3\}$, the increasing term is at most the sum of its values at $M_s$, while each decreasing term is at most its value at $M_s$. Thus, before polynomial counting, we obtain
\begin{align*}
 \frac{\alpha^{2/(3-\theta_1)}}
 {\langle k\rangle^{2\tau/(3-\theta_1)}},\qquad \frac{\varrho(\alpha/\varrho)^{1/(2-\theta_1)}}
 {\langle k\rangle^{\tau/(2-\theta_1)}},\qquad \frac{E_*^{a_f}\alpha^{1-a_f}} {\langle k\rangle^{\tau(1-\theta_0)/(2-\theta_1-\theta_0)}}.
\end{align*}
The condition $\tau^\sharp>0$ implies $\tau/(3-\theta_1)>6$; together with $c\lesssim\langle k\rangle$, this ensures $M_1\ge c^3$ for all sufficiently large $\langle k\rangle$.  For the remaining bounded values of $\langle k\rangle$, the same localization also bounds $c$; hence, after decreasing $\alpha_*$ if necessary, $M_1\ge c^3$ holds uniformly there as well. Moreover the prefactor in \eqref{ckg:large-difference-union-I} costs at most three powers of $\langle k\rangle$ by \eqref{ckg:cd-localization}, and the sum over the $O(\langle k\rangle)$ possible values of $d$ costs one more.
After this four-power counting loss, we lower the three decay exponents to the common conservative exponent $\tau^\sharp$ in \eqref{ckg:tau-sharp}.  Since $\tau^\sharp\le\tau-6$, this also absorbs \eqref{ckg:easy-difference-measure}--\eqref{ckg:intermediate-difference-measure}. This gives exactly the three terms in \eqref{ckg:NG-union-measure}.  

\smallskip
\noindent\emph{Two normal frequencies with the same Melnikov sign.}
As in \cite[Lemmas~5.14--5.16]{BBG}, positivity of $\nu_j(h)\ge |j|/3$ or momentum \eqref{ckg:momentum-kell} gives $|i|+|j|\lesssim\langle k\rangle$ unless $L$ has the opposite sign to both normal signs and $i,j$ have opposite spatial signs. The former cases therefore contribute at most $C\varrho\alpha/\langle k\rangle^{\tau-1}$.

In the remaining case, $L$ has the opposite sign to both normal signs and $i,j$ have opposite spatial signs. Setting $q:=|L|\ge1$ and $g(x):=\sqrt{1+x^2}-1$, resonance and momentum give
\begin{equation}\label{ckg:hard-plus-f-relation}
|g(|i|/c)+g(|j|/c)-q| \le C\langle k\rangle c^{-2}, \qquad \big||i|-|j|\big|\le C\langle k\rangle.
\end{equation}

Arguing as in the proof of Lemma~\ref{ckg:lem-same-unweighted-measure}, using the Lipschitz property of $g$ and the uniform relative Lipschitz bound for $g^{-1}$ at $q/2\ge1/2$, we obtain $O(\langle k\rangle)$ admissible pairs, up to finite block/species multiplicity. Thus all same-sign families contribute at most $C\varrho\alpha/\langle k\rangle^{\tau-1}$, which is absorbed by \eqref{ckg:NG-union-measure} because $\tau^\sharp\le\tau-1$ and $0<\alpha<1$. This completes the proof.
\end{proof}

\subsection{Final measure estimates}

We can state the following measure estimates combining for $m=n$ Lemma~\ref{ckg:lem-same-unweighted-measure} and for $m\ne n$ Lemma~\ref{ckg:lem-local-branch-count}, \ref{ckg:lem-two-parameter-strip} and  Proposition~\ref{ckg:prop-singular-branch-measure}.

\begin{theorem}[Block measure estimate]\label{ckg:thm-measure}
For distinct sites assume the hypotheses of
Proposition~\ref{ckg:prop-singular-branch-measure}, the finite-resolvent
neighbourhood \eqref{ckg:finite-resolvent-neighbourhood}, and $\tau^\sharp>2$.
For coincident sites assume instead the unweighted hypotheses of
Lemma~\ref{ckg:lem-same-unweighted-measure} and $\tau>1$. Then, for a constant $C>0$ independent of $c,\varrho,\alpha$,
\begin{equation}\label{ckg:measure-summed-general}
 \frac1{|\Xi_0|} \left|\bigcup_{k\in\mathbb Z^2}\mathcal R_k^{\theta_1}(\alpha)\right| \le
 \begin{cases}
 C\alpha/\varrho,&m=n,\quad \theta_1=0,\\[0.2em]
 C\left[\alpha^{2/(3-\theta_1)}/\varrho+(\alpha/\varrho)^{1/(2-\theta_1)}
 +E_*^{a_f}\alpha^{1-a_f}/\varrho\right],&m\ne n.
 \end{cases}
\end{equation}
\end{theorem}

\begin{corollary}[High Fourier tails of the resonant set]
\label{ckg:cor-summed-measure}
Under the respective assumptions of Theorem~\ref{ckg:thm-measure}, let $K\ge1$ be a Fourier cutoff.  Then, for a constant $C>0$ independent of $c,\varrho,\alpha,K$,
\begin{align*}
 \frac1{|\Xi_0|} \left|\bigcup_{|k|_1>K}\mathcal R_k^{\theta_1}(\alpha)\right| &\le
 \begin{cases}C(1+K)^{1-\tau}\alpha/\varrho,&m=n,\quad\theta_1=0,\\
 CK^{2-\tau^\sharp}\mathfrak M(\alpha,\varepsilon_{\rm ref},\varepsilon),&m\ne n,\end{cases}
 \\
 \mathfrak M(\alpha,\varepsilon_{\rm ref},\varepsilon) &:=\frac{\alpha^{\frac{2}{3-\theta_1}}}{\varrho} +\left(\frac{\alpha}{\varrho}\right)^{\frac1{2-\theta_1}} +\frac{E_*^{a_f}\alpha^{1-a_f}}{\varrho}.
\end{align*}
The same conclusion holds with $K=0$ if the factor $K^{2-\tau^\sharp}$ is replaced by $1$.
\end{corollary}

\section{The KAM iteration}
\label{sec:CKG-KAM-theorem}

The initial Hamiltonians of Theorem~\ref{ckg:thm-initial-KAM-data} have perturbations of size $\epsilon_0\lesssim r_{\rm K}^2$, where the analytic radius is $r_{\rm K}$, not $r_0$; see \eqref{ckg:initial-KAM-radius} and \eqref{ckg:E0-full}. We use the Newton estimates of \cite[Section~6]{BBG} and \cite[Sections~3--4]{poschel1996kam} after verifying their hypotheses with the block homological equation \eqref{ckg:block-homological-eq}.  The analytic Lie estimates themselves are unchanged: the finite number of species and the block structure are absorbed in Theorem~\ref{ckg:thm-block-homological} and Lemma~\ref{ckg:lem-sparse-block-Schur}.  The points that are not scalar consequences of these references are the finite hyperbolic resolvent, the frozen normal tail, and the simultaneous KG--NLS comparison; they are tracked explicitly below.  Throughout this section the constants are independent of $c$, of $\rho\in\Xi_0$, and of the Newton index, but may depend on the fixed sites, $a,p,\tau$, and a fixed initial angular strip.

\subsection{The analytic KAM step}

\begin{lemma}[One block KAM step]\label{ckg:lem-KAM-step}
Let $H=N+P$ be analytic on $\mathfrak D_q(s,r)$ for $q=p,p-4$. Assume all coefficient, finite-resolvent and symmetry hypotheses of Theorem~\ref{ckg:thm-block-homological}, including \eqref{ckg:homological-smallness}, together with \eqref{ckg:block-Melnikov}.  Suppose, for both values of $q$, that
\[
 \|X_P\|_{r,q,\beta_0}^{\mathcal L(\alpha)}\le\epsilon, \qquad 0<5\sigma<s, \qquad \frac{\epsilon}{\alpha\sigma^\mu}\le c_*.
\]
Set
\begin{equation}\label{ckg:eta-one-step}
 \eta^3:=\frac{\epsilon}{\alpha\sigma^\mu},\qquad 0<\eta<1/9.
\end{equation}
Then the normalized solution of \eqref{ckg:block-homological-eq} generates a time-one symplectic map $\phi$ from $\mathfrak D_q(s-5\sigma,\eta r)$ into $\mathfrak D_q(s,r)$, with
\begin{equation*}
 (N+P)\circ\phi=N_++P_+, \qquad N_+=N+[P].
\end{equation*}
There is a constant $C>0$, independent of $c,\rho,r$ and $\epsilon$, such that
\begin{align}
 \|X_{P_+}\|_{\eta r,q,\beta_1}^{\mathcal L(\alpha)} &\le C\eta\epsilon,
 \label{ckg:one-step-P}\\
 \|\phi-\Id\|_{r,q,\beta_1}^{\mathcal L(\alpha)} &\le C\frac{\epsilon}{\alpha\sigma^\mu},
 \label{ckg:one-step-Phi}\\
 \|X_{[P]}\|_{r,q,\beta_0}^{\mathcal L(\alpha)} &\le C\epsilon.
 \label{ckg:one-step-Ncorr}
\end{align}
\end{lemma}

\begin{proof}
The result follows from the standard analytic KAM-step argument of \cite[proof of Lemma~6.1]{BBG} (see also \cite[Section~3]{poschel1996kam}), using the block homological estimates \eqref{ckg:block-homological-est} and \eqref{ckg:Nhat-est}, together with $C\epsilon^2/(\alpha\sigma^\mu\eta^2) =C\eta\epsilon$.
\end{proof}

\subsection{First step}

Choose once and for all an angular width satisfying $0<s_0\le \min\{2,s_{\rm K}/2\}$, and set
\begin{equation}\label{ckg:newton-beta-definition}
 \sigma_0:=\frac{s_0}{40},\qquad 0<\theta<\frac12, \qquad  \ckgnewtonbeta:=
 \begin{cases}
 0,&m=n,\\
 1-2\theta,&m\ne n.
 \end{cases}
\end{equation}
All constants depending on $s_0$ are hence independent of $c,r_0,\rho$. Let $0<\alpha_1<\alpha_0<c_\alpha\varrho$ and define
\begin{equation*}
 \Xi_1:=\Xi_0\setminus \left(\bigcup_k\mathcal R_k^{\ckgcasebeta\theta}(\alpha_0) \right).
\end{equation*}
With the \emph{available} analytic radius $r_{\rm K}$, put
\begin{equation}\label{ckg:eta0}
 \eta_0^3:=\frac{\epsilon_0}{\alpha_0\sigma_0^\mu},\qquad r_1:=\eta_0 r_{\rm K},\qquad s_1:=s_0-5\sigma_0.
\end{equation}

We now apply the KAM step to the initial KG and NLS Hamiltonians \eqref{ckg:H0-block}-- \eqref{ckg:H0-block-NLS}, with the bounds of Theorem~\ref{ckg:thm-initial-KAM-data}.

\begin{lemma}[Preliminary KAM step]\label{ckg:lem-preliminary-step}
For sufficiently small $h,r_0$ and a fixed small numerical $\varepsilon_*>0$, assume
\begin{equation}\label{ckg:prelim-smallness}
 \frac{\epsilon_0}{\alpha_0\sigma_0^\mu}<\varepsilon_*, \qquad \epsilon_0<\varepsilon_*\varrho.
\end{equation}
On $\Xi_1$ there are KG and NLS symplectic maps $\Phi_0,\Phi_0^{\rm NLS}$, defined on $\mathfrak D_q(s_1,r_1)$ and taking values in $\mathfrak D_q(s_0,r_{\rm K})$, such that
\begin{align*}
 (N_0+P_0)\circ\Phi_0&=N_1+P_1,
 \\
 (N_0^{\rm NLS}+P_0^{\rm NLS})\circ\Phi_0^{\rm NLS} &=N_1^{\rm NLS}+P_1^{\rm NLS}.
\end{align*}
For a sufficiently large numerical constant $C_{\rm step}>0$ and a constant $C>0$, both independent of $c,r_0,\rho$,
\begin{align}
 \|X_{P_1}\|_{r_1,q,\ckgnewtonbeta}^{\mathcal L(\alpha_0)}+\|X_{P_1^{\rm NLS}}\|_{r_1,q}^{\mathcal L(\alpha_0)}& \le\epsilon_1:=C_{\rm step}\eta_0\epsilon_0,
 \label{ckg:epsilon1-def}\\
 \|\Phi_0-\Id\|_{r_{\rm K},q}^{\mathcal L(\alpha_0)} +\|\Phi_0^{\rm NLS}-\Id\|_{r_{\rm K},q}^{\mathcal L(\alpha_0)} &\le C\frac{\epsilon_0}{\alpha_0\sigma_0^\mu},
 \notag\\
 \|\Phi_0-\Phi_0^{\rm NLS}\|_{r_{\rm K},p-4} &\le Ch\frac{\epsilon_0}{\alpha_0^2\sigma_0^{2\mu}}. \notag
\end{align}
Moreover the normal correction has size $O(\epsilon_0)$ and the finite resolvent bound \eqref{ckg:hyp-resolvent-allt} persists. The transformed perturbation admits the decomposition
\begin{equation}\label{ckg:P1-decomp}
P_1=P_1^{\rm NLS}+P_1^{\mathcal{NG}}+P_1^R, \qquad P_1^{\mathcal{NG}}:=(I-\Pi_{\rm LG})P_1, \qquad P_1^R:=\Pi_{\rm LG}P_1-P_1^{\rm NLS},
\end{equation}
with
\begin{equation}\label{ckg:P1R-est}
 \|X_{P_1^{\mathcal{NG}}}\|_{r_1,q}^{\mathcal L(\alpha_0)} \leq \epsilon_1, \qquad \|X_{P_1^R}\|_{r_1,p-4}^{\mathfrak{D}_p(s_1,r_1)} \le Ch\frac{\epsilon_1}{\alpha_0}.
\end{equation}
If
\begin{equation}\label{ckg:tau-first}
 \tau_0^\sharp:=\frac{\tau}{3-\theta}-6>2,
\end{equation}
then
\begin{equation}\label{ckg:first-measure}
 \frac{|\Xi_0\setminus\Xi_1|}{|\Xi_0|}\le
 \begin{cases} C\alpha_0/\varrho,&m=n,\\
 C(\theta)[\alpha_0^{2/(3-\theta)}/\varrho+(\alpha_0/\varrho)^{1/(2-\theta)}],&m\ne n.
 \end{cases}
\end{equation}
\end{lemma}

\begin{proof}
Apply Lemma~\ref{ckg:lem-KAM-step} with $(\theta_0,\theta_1)=(1-\ckgcasebeta,\ckgcasebeta\theta)$ on $\mathfrak D_q(s_0,r_{\rm K})$, as allowed by \eqref{ckg:initial-KAM-radius} and Theorem~\ref{ckg:thm-initial-KAM-data}.  Formula \eqref{ckg:eta0} gives \eqref{ckg:epsilon1-def}.  In the distinct-mode case the $O(\epsilon_0)$ correction lies in the finite-resolvent Neumann neighbourhood of Lemma~\ref{ckg:lem-hyp-resolvent} by \eqref{ckg:prelim-smallness}.  The comparison bounds \eqref{ckg:E0-full}--\eqref{ckg:P0R-jet} and \eqref{ckg:N0R-est} allow us to apply \eqref{ckg:block-generator-comparison}. The standard Gronwall comparison of Lie flows (cf. \cite[Appendix]{BBG}), using the homological estimates at both regularities $p$ and $p-4$, gives the claimed bound for $\Phi_0-\Phi_0^{\rm NLS}$. To estimate $P_1^R$, we subtract the two Lie identities of Lemma~\ref{ckg:lem-KAM-step}, projecting the KG identity by $\Pi_{\rm LG}$. The standard composition and commutator estimates of \cite[proof of Lemma~6.1]{BBG}, together with \eqref{ckg:block-generator-comparison}, give \eqref{ckg:P1R-est}.

The measure estimate \eqref{ckg:first-measure} follows from Theorem~\ref{ckg:thm-measure} with $\theta_1=\ckgcasebeta\theta$ and zero coefficient corrections. There is no additional NLS exclusion by Lemma~\ref{ckg:lem-NLS-no-bad-strips}.

The normal correction produced by this step is subsequently included in the reference; only for $m\ne n$ is this a weighted-to-unweighted restart.
It contributes $\varepsilon_{\rm ref}=C\epsilon_0$
and $\kappa_{\rm ref}=C\epsilon_0/\alpha_0$,
rather than a correction measured against $\alpha_1$.

\end{proof}

\subsection{Infinite Newton iteration}

Put $\sigma_1=\sigma_0/2$ and fix $K_1\ge1$ independently of $c,r_0,\rho$.  For $\nu\ge1$ define
\begin{align}
 \sigma_{\nu+1}&:=\tfrac12\sigma_\nu, &\alpha_\nu&:=\tfrac{\alpha_1}{2}(1+2^{-\nu+1}), &K_\nu&:=2^{\nu-1}K_1,
 \label{ckg:iteration-parameters-1}\\
 \epsilon_{\nu+1}&:= \frac{C_1\epsilon_\nu^{4/3}} {(\alpha_\nu\sigma_\nu^\mu)^{1/3}}, &\eta_\nu^3&:=\frac{\epsilon_\nu} {\alpha_\nu\sigma_\nu^\mu}, &r_{\nu+1}&:=\eta_\nu r_\nu,
 \label{ckg:iteration-parameters-2}\\
 s_{\nu+1}&:=s_\nu-5\sigma_\nu. \notag
\end{align}
Here $r_1$ is defined by \eqref{ckg:eta0}; since $\sum_{\nu\ge1}\sigma_\nu=\sigma_0$, one has $s_\infty=s_0-10\sigma_0=3s_0/4>0$.

\begin{lemma}\label{ckg:lem-majorants}
Fix $C_1$ sufficiently large to dominate the constant in \eqref{ckg:one-step-P}.  There exists $\varepsilon_*=\varepsilon_*(C_1,\mu)>0$ such that if
\begin{equation}\label{ckg:majorant-smallness}
 x_1:=\frac{\epsilon_1}{\alpha_1\sigma_1^\mu}<\varepsilon_*,
\end{equation}
then, for all $\nu\ge1$,
\begin{equation*}
 x_{\nu+1}\le \tfrac43C_1 2^\mu x_\nu^{4/3},
\end{equation*}
and, after decreasing $\varepsilon_*$,
\begin{align}
 \eta_\nu&<2^{-\nu-5},\notag\\
 \sum_{j\ge\nu}\epsilon_j&\le2\epsilon_\nu, \qquad \sum_{j\ge\nu}\frac{\epsilon_j}{\alpha_j\sigma_j^\mu} \le Cx_\nu.\label{ckg:epsilon-sums}
\end{align}
The constants are independent of $c,\rho$ and $\nu$.
\end{lemma}

\begin{proof}
By \eqref{ckg:iteration-parameters-1}--\eqref{ckg:iteration-parameters-2},
\[
 x_{\nu+1}=C_1 2^\mu \frac{\alpha_\nu}{\alpha_{\nu+1}}x_\nu^{4/3}, \qquad 1\le\frac{\alpha_\nu}{\alpha_{\nu+1}}\le\frac43.
\]
Choosing $x_1$ below a fixed sufficiently small threshold, the superlinear recurrence gives $x_\nu\le2^{-3\nu-18}$. Since $\eta_\nu=x_\nu^{1/3}$, the first estimate follows; the two series in \eqref{ckg:epsilon-sums} are bounded by geometric majorants of the superexponentially decreasing sequence.
\end{proof}

We now define the parameter exclusions for the Newton iteration. Before each parameter exclusion, the normal-form corrections are extended to $\Xi_0$ by Lemma~\ref{ckg:lem-extension}, preserving their individual weighted supremum and Lipschitz bounds. The tangential KG--NLS corrections are extended in a coupled way, preserving their global $O(h)$ comparison, as explained in Remark~\ref{ckg:lem-autonomous-NLS}.

By Lemma~\ref{ckg:lem-NLS-no-bad-strips}, only KG resonances need to be excluded. For $m\ne n$ we restart with full unweighted Melnikov conditions for $N_1$; for $m=n$ we retain the low Fourier indices and exclude only the new high ones:
\begin{equation}\label{ckg:Xi2-restart}
 \Xi_2:=\Xi_1\setminus
 \begin{cases}
 \displaystyle\bigcup_{|k|_1>K_1}(\mathcal R_{k,1}^{0}(\alpha_1)),&m=n,\\[0.5em]
 \displaystyle\bigcup_{k\in\mathbb Z^2}(\mathcal R_{k,1}^{0}(\alpha_1)),&m\ne n.
 \end{cases}
\end{equation}
Here and below, $\mathcal R_{k,\nu}^{\theta_1}(\alpha)$ denotes the KG resonant set $\mathcal R_k^{\theta_1}(\alpha)$ computed using the extended normal-form coefficients of $N_\nu$.

From the second Newton step onward the Melnikov weight is unchanged. For $\nu\ge2$ define
\begin{equation}\label{ckg:Xi-nu-def}
 \Xi_{\nu+1}:=\Xi_\nu\setminus \bigcup_{|k|_1>K_\nu} \left(\mathcal R_{k,\nu}^{0}(\alpha_\nu)\right).
\end{equation}

\begin{lemma}[Persistence of old Melnikov conditions]
\label{ckg:lem-old-Melnikov}
Let $\nu\ge2$, or $\nu=1$ when $m=n$. Suppose that the normal-form coefficient increments from $N_{\nu-1}$ to $N_\nu$ are bounded by $C\epsilon_{\nu-1}$. Assume the \emph{unweighted} inverse Melnikov estimate with threshold $\alpha_{\nu-1}$ holds on $\Xi_\nu$ for $N_{\nu-1}$. Then, for all $|k|_1\le K_\nu$, the corresponding unweighted estimate for $N_\nu$ holds with threshold $\alpha_\nu$ if
\begin{equation}\label{ckg:old-mel-small}
 C\epsilon_{\nu-1}\langle K_\nu\rangle^{\tau+2} \le \alpha_{\nu-1}-\alpha_\nu.
\end{equation}
For fixed $K_1\ge1$ and sufficiently small $r_0$, with a smallness threshold also depending on $K_1$ and $\tau$, \eqref{ckg:old-mel-small} holds for every $\nu\ge2$. For the NLS operators, the required Melnikov estimates instead follow directly from Lemma~\ref{ckg:lem-NLS-no-bad-strips}, without a persistence argument.
\end{lemma}

\begin{proof}
The difference of homological operators of Fourier index $k$ is at most $C\epsilon_{\nu-1}\langle k\rangle$ in operator norm, since the tangential term contains $\langle k,\omega\rangle$. Writing $L_\nu=L_{\nu-1}+E_\nu$ and using \eqref{ckg:KG-NLS-resolvent-identity},
\[
 L_\nu^{-1}=(I+L_{\nu-1}^{-1}E_\nu)^{-1}L_{\nu-1}^{-1},
\]
we obtain the desired inverse bound whenever
\[
 C\epsilon_{\nu-1}\langle k\rangle^{\tau+2} \le \alpha_{\nu-1}-\alpha_\nu.
\]
For $\nu\ge2$ the factor $\langle K_\nu\rangle^{\tau+2}$ grows at most exponentially in $\nu$, whereas $\epsilon_{\nu-1}$ decays superexponentially by Lemma~\ref{ckg:lem-majorants}. After choosing its initial threshold sufficiently small (depending on the fixed $K_1,\tau$), this dominates the shrinking gap $\alpha_{\nu-1}-\alpha_\nu\asymp\alpha_1 2^{-\nu}$. In particular, for $\nu=2$ it suffices to require $C\epsilon_1\langle2K_1\rangle^{\tau+2}\le\alpha_1/4$. The finite-containing resolvents are controlled directly by Lemma~\ref{ckg:lem-hyp-resolvent}, as they require no scalar Melnikov exclusion.
\end{proof}

\begin{lemma}[Iterative block KAM lemma]\label{ckg:lem-iterative-KAM}
Under \eqref{ckg:prelim-smallness}, \eqref{ckg:majorant-smallness}, \eqref{ckg:old-mel-small} for $\nu\ge2$, and the hypotheses of Theorem~\ref{ckg:thm-block-homological}, there exist single-step symplectic maps $\phi_\nu,\phi_\nu^{\rm NLS}$, $\nu\ge1$, such that
\[
 (N_\nu+P_\nu)\circ\phi_\nu=N_{\nu+1}+P_{\nu+1},
\]
with its NLS counterpart on $\Xi_{\nu+1}$. Define the \emph{cumulative} maps, including the preliminary step, by
\[
 \Phi_\nu:=\Phi_0\circ\phi_1\circ\cdots\circ\phi_\nu, \qquad \Phi_\nu^{\rm NLS}:=\Phi_0^{\rm NLS}\circ \phi_1^{\rm NLS}\circ\cdots\circ\phi_\nu^{\rm NLS}.
\]
They satisfy
\begin{align}
 (N_0+P_0)\circ\Phi_\nu&=N_{\nu+1}+P_{\nu+1},
 \label{ckg:iteration-conjugacy-KG}\\
 (N_0^{\rm NLS}+P_0^{\rm NLS})\circ\Phi_\nu^{\rm NLS} &=N_{\nu+1}^{\rm NLS}+P_{\nu+1}^{\rm NLS}.
 \label{ckg:iteration-conjugacy-NLS}
\end{align}
There is a constant $C>0$ independent of $c,\rho,\nu$ such that, for $q=p,p-4$,
\begin{align}
 \|X_{P_\nu}\|_{r_\nu,q,\ckgnewtonbeta}^{\mathcal L(\alpha_\nu)}+\|X_{P_\nu^{\rm NLS}} \|_{r_\nu,q}^{\mathcal L(\alpha_\nu)} &\le\epsilon_\nu,
 \label{ckg:iter-perturbation-size}\\
 P_\nu^{\mathcal{NG}}=(I-\Pi_{\rm LG})P_\nu, \qquad P_\nu^R=\Pi_{\rm LG}P_\nu-P_\nu^{\rm NLS},
 \label{ckg:iter-P-decomp}\\
\|X_{P_\nu^{\mathcal{NG}}}\|_{r_\nu,q}^{\mathcal L(\alpha_\nu)} \leq \epsilon_\nu, \qquad  \|X_{P_\nu^R}\|_{r_\nu,p-4}^{\mathfrak{D}_p(s_\nu,r_\nu)} &\le Ch\frac{\epsilon_\nu}{\alpha_{\nu-1}},
 \label{ckg:iter-PR-bound}\\
 |\delta_\nu|_\infty+|D_\nu|_{\infty,\ckgnewtonbeta} +\|E_{\cF,\nu}\| &\le C(\epsilon_0+\epsilon_1),
 \label{ckg:iter-freq-size}\\
 |\delta_\nu|_\infty^{\rm lip} +|D_\nu|_{\infty,\ckgnewtonbeta}^{\rm lip} +\|E_{\cF,\nu}\|^{\rm lip} &\le C\left(\frac{\epsilon_0}{\alpha_0} +\frac{\epsilon_1}{\alpha_1}\right),
 \label{ckg:iter-freq-Lip}
\end{align}
with $\delta_\nu$, $D_\nu$ and $E_{\cF,\nu}$ denoting the accumulated normal-form coefficient corrections at stage $\nu$, relative to the initial coefficients in \eqref{ckg:initial-shifted-coefficients}, including the preliminary correction. The analogous bounds hold for NLS.  On the limiting torus, with the scaled norm based on $r_{\rm K}$,
\begin{align}
 \|\Phi_\nu-\Id\|_{r_{\rm K},q}^{\mathcal L(\alpha_1)} +\|\Phi_\nu^{\rm NLS}-\Id\|_{r_{\rm K},q}^{\mathcal L(\alpha_1)} &\le C\left(\frac{\epsilon_0}{\alpha_0\sigma_0^\mu} +\frac{\epsilon_1}{\alpha_1\sigma_1^\mu}\right),
 \label{ckg:iter-Phi-size}\\
 \|\Phi_\nu-\Phi_\nu^{\rm NLS}\|_{r_{\rm K},p-4} &\le Ch\left(\frac{\epsilon_0}{\alpha_0^2\sigma_0^{2\mu}} +\frac{\epsilon_1}{\alpha_1^2\sigma_1^{2\mu}}\right).
 \label{ckg:iter-Phi-diff}
\end{align}
All finite-block resolvent estimates remain valid.  If
\begin{equation}\label{ckg:tau-infinite}
 \tau_*:=\frac{\tau(1-2\theta)}3-6>2,
\end{equation}
put
\begin{equation*}
 a_\theta:=\frac1{2(1-\theta)}
\end{equation*}
and
\begin{equation*}
 \mathfrak M_1:= \frac{\alpha_1^{2/3}}{\varrho} +\left(\frac{\alpha_1}{\varrho}\right)^{1/2} +\frac{(\epsilon_0+\epsilon_1)^{a_\theta} \alpha_1^{1-a_\theta}}{\varrho}.
\end{equation*}
Then
\begin{equation}\label{ckg:iter-measure}
 \frac{|\Xi_0\setminus\Xi_\nu|}{|\Xi_0|}\le
 \begin{cases}
 C\dfrac{\alpha_0}{\varrho}\left(1+\displaystyle\sum_{j=1}^{\nu-1}(1+K_j)^{1-\tau}\right),&m=n,\\[0.5em]
 C(\theta)\left[\dfrac{\alpha_0^{2/(3-\theta)}}{\varrho}
 +(\dfrac{\alpha_0}{\varrho})^{1/(2-\theta)}
 +\mathfrak M_1\left(1+\displaystyle\sum_{j=2}^{\nu-1}2^{-(\tau_*-2)(j-1)}\right)\right],&m\ne n.
 \end{cases}
\end{equation}
\end{lemma}

\begin{proof}
The standard Newton induction follows from Lemmas~\ref{ckg:lem-KAM-step}, \ref{ckg:lem-majorants}, and \ref{ckg:lem-old-Melnikov}, as in \cite[Section~4]{poschel1996kam}. At $\nu=1$, the unweighted KG Melnikov conditions follow from the full restart for $m\ne n$, and from persistence and the high-index exclusion for $m=n$. Subsequent steps use \eqref{ckg:Xi-nu-def}. The finite resolvents remain controlled by Lemma~\ref{ckg:lem-hyp-resolvent}, while NLS requires no exclusions by Lemma~\ref{ckg:lem-NLS-no-bad-strips}. The one-step bounds and \eqref{ckg:epsilon-sums} give \eqref{ckg:iter-perturbation-size}, \eqref{ckg:iter-freq-size}--\eqref{ckg:iter-freq-Lip}, and \eqref{ckg:iter-Phi-size}.

For the KG--NLS comparison, set
\[
R_\nu:= \|X_{P_\nu^R}\|_{r_\nu,p-4}^{\mathfrak D_p(s_\nu,r_\nu)}.
\]
The generator estimates \eqref{ckg:SNG-est} and \eqref{ckg:block-generator-comparison}, together with the analytic Lie-remainder estimates of \cite[proof of Lemma~6.1]{BBG}, give
\[
R_{\nu+1} \le C\eta_\nu R_\nu+ Ch\frac{\epsilon_\nu^2}{\alpha_\nu^2\sigma_\nu^\mu\eta_\nu^2} \lesssim \eta_\nu\left( R_\nu+h\frac{\epsilon_\nu}{\alpha_\nu} \right).
\]
Here the LG projection of brackets between two non-LG terms must be retained, since their fast defects may cancel. By \eqref{ckg:P1R-est} and $\epsilon_{\nu+1}=C_1\eta_\nu\epsilon_\nu$, this closes \eqref{ckg:iter-PR-bound} for $C_1$ sufficiently large. Moreover, $[P_\nu]-[P_\nu^{\rm NLS}]=[P_\nu^R]$, so \eqref{ckg:epsilon-sums} preserves the KG--NLS normal-form comparison, including the finite-block bound.

The standard flow comparison in the weak $p-4$ topology (cf. \cite[Appendix]{BBG}), using the homological estimates at both $p$ and $p-4$, yields
\begin{equation}\label{ckg:one-step-Phi-diff}
\|\phi_\nu-\phi_\nu^{\rm NLS}\|_{r_\nu,p-4} \le Ch\frac{\epsilon_\nu} {\alpha_\nu^2\sigma_\nu^\mu}.
\end{equation}
The available angular Cauchy reserve follows from $\mu=2\tau+6$ and $\tau\ge2$. Telescopic composition, including the preliminary comparison, proves \eqref{ckg:iter-Phi-diff}.

Finally, the initial loss is bounded by \eqref{ckg:first-measure}. The subsequent losses follow from Lemma~\ref{ckg:lem-same-unweighted-measure} for $m=n$, and from Corollary~\ref{ckg:cor-summed-measure} with $(\theta_0,\theta_1)=(2\theta,0)$ for the distinct-mode restart and its high-Fourier tails. Summing with $K_\nu=2^{\nu-1}K_1$ gives \eqref{ckg:iter-measure}.
\end{proof}

\subsection{Convergence}

The convergence of the cumulative KAM maps and their first differentials follows from the standard composition argument of \cite[proof of Theorem~A]{poschel1996kam}, using Lemma~\ref{ckg:lem-majorants}. The KG--NLS comparison passes to the limit from Lemma~\ref{ckg:lem-iterative-KAM}, whose weak-space telescopic estimates follow the argument of \cite[Section~6]{BBG}. Let 
\begin{equation*}
 \Xi_\infty:=\bigcap_{\nu\ge1}\Xi_\nu,
\end{equation*}
We obtain the following result.
\begin{lemma}[Convergence of the KAM maps]\label{ckg:lem-KAM-convergence}
The cumulative maps $\Phi_\nu,\Phi_\nu^{\rm NLS}$ defined in Lemma~\ref{ckg:lem-iterative-KAM}, and their first differentials on a fixed smaller angular strip, converge on the limiting zero section to real analytic embeddings $\Phi_\infty,\Phi_\infty^{\rm NLS}$. The slow normal-form coefficients converge on $\Xi_\infty$ and have uniform Lipschitz extensions to $\Xi_0$.  On the limiting torus, for a constant $C>0$ independent of $c,\rho$,
\begin{align}
 \|\Phi_\infty-\Id\|_{r_{\rm K},p} &\le C\left( \frac{\epsilon_0}{\alpha_0\sigma_0^\mu} +\frac{\epsilon_1}{\alpha_1\sigma_1^\mu}\right),
 \label{ckg:Phi-infty-size}\\
 \|\Phi_\infty-\Id\|_{r_{\rm K},p}^{\rm lip} &\le C\left( \frac{\epsilon_0}{\alpha_0^2\sigma_0^{2\mu}} +\frac{\epsilon_1}{\alpha_1^2\sigma_1^{2\mu}}\right),
 \label{ckg:Phi-infty-Lip}\\
 \|\Phi_\infty-\Phi_\infty^{\rm NLS}\|_{r_{\rm K},p-4} &\le Ch\left( \frac{\epsilon_0}{\alpha_0^2\sigma_0^{2\mu}} +\frac{\epsilon_1}{\alpha_1^2\sigma_1^{2\mu}}\right).
 \label{ckg:Phi-infty-diff}
\end{align}
The limiting perturbation has vanishing KAM jet.
\end{lemma}

Combining Lemmas~\ref{ckg:lem-preliminary-step}, \ref{ckg:lem-iterative-KAM}, and \ref{ckg:lem-KAM-convergence}, we obtain the following KAM theorem. The uniform bi-Lipschitz bounds for the limiting frequency maps follow from \eqref{ckg:twist-close} and \eqref{ckg:iter-freq-Lip}. Their $O(h)$ comparison follows by summing the KG--NLS normal-form differences, using \eqref{ckg:N0R-est}, \eqref{ckg:P0R-jet}, \eqref{ckg:iter-PR-bound}, and the paired extensions of Lemma~\ref{ckg:lem-extension}. Finally, linear stability for $m=n$ follows from the Hermitian normal blocks, whereas for $m\ne n$ the finite resolvent estimate \eqref{ckg:hyp-resolvent-allt} yields hyperbolic normal exponents with $|\operatorname{Re}\lambda|\gtrsim\varrho$.

\begin{theorem}\label{TeoremaGenerale}
Fix the tangential sites $m,n\in\mathbb Z$, $a\ge0$, $p>9/2$, $0<\theta<1/2$, and $\tau$ satisfying \eqref{ckg:tau-first} and \eqref{ckg:tau-infinite}. Let $H_0$ and $H_0^{\rm NLS}$ be the Hamiltonians in \eqref{ckg:H0-block}--\eqref{ckg:H0-block-NLS}, with action-parameter domain $\Xi_0$ and $\varrho=r_0^{4/3}$. There exist $h_*,r_*,\varepsilon_*>0$ such that, for $0<h=c^{-2}<h_*$, $0<r_0<r_*$, and $0<\alpha_1<\alpha_0<c_\alpha\varrho$, if
\begin{equation}\label{ckg:general-smallness}
\frac{\epsilon_0}{\alpha_0\sigma_0^\mu}+\frac{\epsilon_1}{\alpha_1\sigma_1^\mu}+\frac{\epsilon_0}{\varrho}<\varepsilon_*,
\end{equation}
and the Melnikov persistence condition \eqref{ckg:old-mel-small} holds for $\nu\ge2$ (and also for $\nu=1$ when $m=n$), then there exist a Cantor set $\Xi_\infty\subset\Xi_0$, two bi-Lipschitz frequency maps
\[
\widehat\omega_\infty^{\rm KG},\ \omega_\infty^{\rm NLS}:\Xi_0\longrightarrow\mathbb R^2,
\]
and two families of invariant-torus embeddings
\begin{equation}\label{Phi}
\Phi_\infty(\rho)\circ\iota,\quad\Phi_\infty^{\rm NLS}(\rho)\circ\iota:\mathbb T^2\longrightarrow\cP_{\textnormal{KAM}}^{a,p,0},\qquad \rho\in\Xi_\infty,
\end{equation}
where $\iota(\vartheta)=(\vartheta,0,0)$. The embeddings are real analytic for $|\operatorname{Im}\vartheta|<3s_0/8$ and Lipschitz in $\rho$; their internal frequency vectors are, respectively, $c^2(1,1)+\widehat\omega_\infty^{\rm KG}(\rho)$ and $\omega_\infty^{\rm NLS}(\rho)$. The frequency maps have direct and inverse Lipschitz constants independent of $h,r_0$, and the following estimates hold with a constant $C$ independent of $h,r_0,\rho$.

On the limiting zero section, the cumulative KAM maps satisfy
\begin{align}
\|\Phi_\infty-\Id\|_{r_{\rm K},p}+\|\Phi_\infty^{\rm NLS}-\Id\|_{r_{\rm K},p}&\le C\left(\frac{\epsilon_0}{\alpha_0\sigma_0^\mu}+\frac{\epsilon_1}{\alpha_1\sigma_1^\mu}\right),\label{ckg:general-Phi-size}\\
\|\Phi_\infty-\Id\|_{r_{\rm K},p}^{\rm lip}+\|\Phi_\infty^{\rm NLS}-\Id\|_{r_{\rm K},p}^{\rm lip}&\le C\left(\frac{\epsilon_0}{\alpha_0^2\sigma_0^{2\mu}}+\frac{\epsilon_1}{\alpha_1^2\sigma_1^{2\mu}}\right),\label{ckg:general-Phi-Lip}\\
\|\Phi_\infty-\Phi_\infty^{\rm NLS}\|_{r_{\rm K},p-4}&\le Ch\left(\frac{\epsilon_0}{\alpha_0^2\sigma_0^{2\mu}}+\frac{\epsilon_1}{\alpha_1^2\sigma_1^{2\mu}}\right).\label{ckg:general-Phi-diff}
\end{align}
The limiting tangential-frequency and normal-block corrections satisfy
\begin{align}
|\delta_\infty|_\infty+|D_\infty|_{\infty,\ckgnewtonbeta}+\|E_{\cF,\infty}\|&\le C(\epsilon_0+\epsilon_1),\label{ckg:general-normal-corrections}\\
|\delta_\infty|_\infty^{\rm lip}+|D_\infty|_{\infty,\ckgnewtonbeta}^{\rm lip}+\|E_{\cF,\infty}\|^{\rm lip}&\le C\left(\frac{\epsilon_0}{\alpha_0}+\frac{\epsilon_1}{\alpha_1}\right).\label{ckg:general-normal-corrections-lip}
\end{align}
The analogous bounds hold for the NLS corrections $\delta_{\infty}^{\rm NLS}$, $D_\infty^{\rm NLS}$ and $E_{\cF,\infty}^{\rm NLS}$, with $|D_\infty^{\rm NLS}|_{\infty,0}$ in place of $|D_\infty|_{\infty,\ckgnewtonbeta}$; finite-block terms are omitted when $m=n$. Moreover,
\begin{equation}\label{ckg:final-frequency-comparison}
\|\widehat\omega_\infty^{\rm KG}-\omega_\infty^{\rm NLS}\|_\infty \le Ch.
\end{equation}

The excluded parameters satisfy
\begin{equation}\label{ckg:general-measure-final}
\frac{|\Xi_0\setminus\Xi_\infty|}{|\Xi_0|}\le
\begin{cases}
C\dfrac{\alpha_0}{\varrho},&m=n,\\[0.4em]
C(\theta)\left[\dfrac{\alpha_0^{2/(3-\theta)}}{\varrho}+\left(\dfrac{\alpha_0}{\varrho}\right)^{1/(2-\theta)}+\mathfrak M_1\right],&m\ne n,
\end{cases}
\end{equation}
where $\mathfrak M_1$ is defined in Lemma~\ref{ckg:lem-iterative-KAM}. The tori are linearly stable for $m=n$. For $m\ne n$, they are partially hyperbolic, with finite hyperbolic Floquet exponents whose real parts have modulus at least $\kappa\varrho$, where $\kappa>0$ is independent of $h,r_0$.
\end{theorem}

\begin{remark}[Parameter extensions and autonomous NLS family]
\label{ckg:lem-autonomous-NLS}
We briefly clarify the choice of parameter extensions. At each KAM step, the normal-form corrections are initially defined on the retained Cantor set. By Lemma~\ref{ckg:lem-extension}, they admit extensions to $\Xi_0$ preserving their supremum and Lipschitz bounds up to absolute constants. For the tangential KG--NLS comparison, we extend the NLS correction and the difference between the KG and NLS corrections separately, and define the KG extension as their sum. This preserves the $O(h)$ comparison on $\Xi_0$. Moreover, since the initial frequency maps are bi-Lipschitz and the extended corrections have sufficiently small Lipschitz seminorm, the resulting frequency maps remain uniformly bi-Lipschitz, by the elementary perturbation lemma in \cite[Appendix D.]{BBG}.

We also clarify the independence of the NLS family from $h$. Fixing all KAM parameters and cutoffs independently of $h$, Lemma~\ref{ckg:lem-NLS-no-bad-strips} allows the NLS iteration to be performed autonomously on the whole box $\Xi_0$, without Melnikov exclusions. Since the initial NLS Hamiltonian and its preliminary transformations are independent of $h$, so are all subsequent NLS normal forms, homological generators and KAM maps. Indeed, the normalized homological equations determine their solutions uniquely at each parameter value. The NLS objects appearing in the simultaneous KG--NLS construction are precisely the restrictions of this autonomous iteration to $\Xi_\nu(h)$. Choosing the NLS extensions from the autonomous construction, independently of $h$, and extending only the KG--NLS tangential differences on the common Cantor sets, yields the global $O(h)$ frequency comparison without altering the NLS family. Consequently, both $\omega_\infty^{\rm NLS}$ and the NLS torus embeddings are independent of $h$.
\end{remark}

\section{The nonrelativistic results}
In this section we complete the proofs of Theorem~\ref{thm:intro-CKG} and Theorem~\ref{cor:intro-CKG-subsequence}. The preceding construction has been carried out for the defocusing choice ${\rm a}_1=1$ in \eqref{eq:intro-CKG}--\eqref{eq:intro-NLS}. The focusing case ${\rm a}_1=-1$ follows from the same arguments, with the corresponding sign changes in the quartic normal forms. In particular, the leading tangential twist matrix $A_*$ in \eqref{ckg:twist-matrix} is replaced by $-A_*$, preserving its singular values and hence its nondegeneracy. In the distinct-mode case, the finite normal dynamical matrices in \eqref{ckg:LF-def} also change sign, leaving the hyperbolic separation in Lemma~\ref{ckg:lem-hyperbolic} and the uniform resolvent estimates of Lemma~\ref{ckg:lem-hyp-resolvent} unchanged. The remaining perturbative estimates are unaffected at the level of the bounds used in the KAM iteration.

\subsection{Choice of the small parameters and common frequencies}
\label{sec:common-frequencies}

We now select the small parameters and match the two limiting frequency maps.  The overlap argument follows the geometry of \cite[Section~7]{BBG}; we retain a short proof since some arguments used here are different from those presented in that paper. 

Recalling \eqref{ckg:general-smallness} and \eqref{ckg:varrho-r0}, we set
\begin{equation}\label{ckg:varsigma-choice}
 0<\varsigma<\frac1{18},\qquad \theta:=\frac12-3\varsigma, \qquad a_0:=\frac53+3\varsigma, \qquad a_1:=2+\varsigma,
\end{equation}
and choose
\begin{equation}\label{ckg:alpha-choices}
 \alpha_0=r_0^{a_0},\qquad
 \alpha_1=\begin{cases}\alpha_0/2,&m=n,\\r_0^{a_1},&m\ne n,\end{cases}
 \qquad \varrho=r_0^{4/3}.
\end{equation}
Since $\epsilon_0=O(r_0^2)$,
\begin{align}
 \frac{\epsilon_0}{\alpha_0}
 &=O(r_0^{1/3-3\varsigma}),\label{ckg:explicit-smallness-0}\\
 \epsilon_1 &=O\!\left(r_0^{19/9-\varsigma}\right),\qquad
 \frac{\epsilon_1}{\alpha_1}=\begin{cases}
 O(r_0^{4/9-4\varsigma}),&m=n,\\
 O(r_0^{1/9-2\varsigma}),&m\ne n,
 \end{cases}
 \label{ckg:explicit-smallness-1}
\end{align}
that combined with the choices in \eqref{ckg:varsigma-choice} make admissible the condition \eqref{ckg:general-smallness} for $r_0$ sufficiently small. Moreover recalling \eqref{ckg:old-mel-small}, it holds $C\epsilon_0\langle K_1\rangle^{\tau+2}\le\alpha_0-\alpha_1$ when $m=n$, by further reducing \(r_0\), if necessary.

With the choices \eqref{ckg:varsigma-choice}--\eqref{ckg:alpha-choices}, we make explicit the terms in \eqref{ckg:general-measure-final}. For $m \neq n$ we have
\begin{align}
 \frac{\alpha_0^{2/(3-\theta)}}{\varrho} =r_0^{\frac{4\varsigma}{5+6\varsigma}}, \quad \left(\frac{\alpha_0}{\varrho}\right)^{1/(2-\theta)} =r_0^{\frac{2(1+9\varsigma)}{9(1+2\varsigma)}}, \quad  \frac{\alpha_1^{2/3}}{\varrho} =r_0^{2\varsigma/3}, \quad  \left(\frac{\alpha_1}{\varrho}\right)^{1/2} =r_0^{1/3+\varsigma/2},
 \label{ckg:explicit-measure-rates-0}\\
 \frac{(\epsilon_0+\epsilon_1)^{a_\theta} \alpha_1^{1-a_\theta}}{\varrho} \le C r_0^{\frac23+\frac{6\varsigma^2}{1+6\varsigma}}.
 \label{ckg:explicit-measure-rates-frozen}
\end{align}
All exponents of \(r_0\) on the right-hand sides are strictly positive for \eqref{ckg:varsigma-choice} and the slowest decay is $r_0^{2\varsigma/3}$. In the identical-mode case \eqref{ckg:iter-measure} instead gives $C\alpha_0/\varrho=C r_0^{1/3+3\varsigma}$, which is bounded by the same weaker estimate $Cr_0^{2\varsigma/3}$. Thus, in both cases, the retained parameter set has asymptotically full relative measure. Let
\[
 f_h:=\widehat\omega_\infty^{\rm KG},\qquad f_0:=\omega_\infty^{\rm NLS},
\]
introduced in Theorem \ref{TeoremaGenerale}. Both are lipeomorphisms on $\Xi_0$ with uniform bi-Lipschitz constants, and \eqref{ckg:final-frequency-comparison} gives $\|f_h-f_0\|_{\infty}\le Ch$. Define
\begin{equation}\label{Wh}
 \mathcal W_h:=f_h(\Xi_\infty)\cap f_0(\Xi_\infty).
\end{equation}

\begin{lemma}\label{ckg:lem-image-overlap}
Let $Q=[\varrho/2,3\varrho/2]^2$ and let $f,g:Q\to\mathbb{R}^2$ be lipeomorphisms with inverse Lipschitz constants at most $L$.  If $\|f-g\|_\infty\le\delta$, then
\begin{equation}\label{ckg:inner-image-inclusion}
 f(Q_{-L\delta})\subset f(Q)\cap g(Q),
\end{equation}
where $Q_{-d}:=\{x\in Q:\operatorname{dist}(x,\partial Q)>d\}$. Consequently
\begin{equation}\label{ckg:image-boundary-loss}
 |Q\setminus f^{-1}(f(Q)\cap g(Q))|\le C_L\delta\varrho.
\end{equation}
\end{lemma}

\begin{proof}
Since $g$ is a lipeomorphism, invariance of domain implies $g(\partial Q)=\partial g(Q)$.  If $x\in Q_{-L\delta}$, the inverse Lipschitz bound gives
\[
 \operatorname{dist}(g(x),\partial g(Q)) =\inf_{z\in\partial Q}|g(x)-g(z)| \ge L^{-1}\operatorname{dist}(x,\partial Q)>\delta.
\]
Because $|f(x)-g(x)|\le\delta$, it follows that $f(x)\in g(Q)$.  Only at this point do we define $y:=g^{-1}(f(x))\in Q$; then
\[
 |y-x|\le L|f(x)-g(x)|\le L\delta.
\]
Hence \eqref{ckg:inner-image-inclusion} holds. Taking preimages in \eqref{ckg:inner-image-inclusion} gives
\[
Q\setminus f^{-1}(f(Q)\cap g(Q))\subseteq Q\setminus Q_{-L\delta}.
\]
The complement of $Q_{-L\delta}$ is a boundary layer of area at most $C\delta\varrho$ and we obtain \eqref{ckg:image-boundary-loss}.
\end{proof}

\begin{theorem}[Common-frequency Cantor families] \label{ckg:thm-common-frequency}
Under the choices \eqref{ckg:varsigma-choice}--\eqref{ckg:alpha-choices}, there exist $r_*,h_*,\kappa_0,C>0$ such that,for $0<r_0<r_*$ and $0<h<\min\{h_*,\kappa_0\varrho\}$,the set $\mathcal W_h$ defined in \eqref{Wh} has positive measure. For every $\omega\in\mathcal W_h$, the uniquely determined parameters are
\begin{equation}\label{ckg:two-parameters-common}
\begin{aligned}
\rho_h(\omega)&:=f_h^{-1}(\omega)\in\Xi_\infty,\\
\rho_0(\omega)&:=f_0^{-1}(\omega)\in\Xi_\infty,
\end{aligned}
\qquad
f_h(\rho_h(\omega)) =\omega=f_0(\rho_0(\omega)).
\end{equation}
They satisfy
\begin{equation}\label{ckg:rho-distance}
|\rho_h(\omega)-\rho_0(\omega)|\le Ch.
\end{equation}
Furthermore,
\begin{equation}\label{ckg:common-measure}
\frac{|\Xi_0\setminus f_h^{-1}(\mathcal W_h)|}{|\Xi_0|}+\frac{|\Xi_0\setminus f_0^{-1}(\mathcal W_h)|}{|\Xi_0|}\le C\left(r_0^{2\varsigma/3}+\frac{h}{\varrho}\right).
\end{equation}
\end{theorem}

\begin{proof}
Set $Q:=\Xi_0$, $E:=\Xi_\infty$, and $B:=Q\setminus E$. By Theorem~\ref{TeoremaGenerale}, the maps $f_h,f_0$ have uniform direct and inverse Lipschitz constants and $\|f_h-f_0\|_\infty\le Ch$. Applying Lemma~\ref{ckg:lem-image-overlap} in both directions and using the direct Lipschitz bounds gives
\[
|f_h(Q)\setminus f_0(Q)| + |f_0(Q)\setminus f_h(Q)| \le Ch\varrho.
\]

Since $\mathcal W_h=f_h(E)\cap f_0(E)$,
\[
f_h(Q)\setminus\mathcal W_h\subseteq \bigl(f_h(Q)\setminus f_0(Q)\bigr)\cup f_h(B)\cup f_0(B),
\]
and the analogous inclusion holds with $f_h$ and $f_0$ interchanged. The uniform bi-Lipschitz bounds therefore yield
\[
\begin{aligned}
&|Q\setminus f_h^{-1}(\mathcal W_h)|+|Q\setminus f_0^{-1}(\mathcal W_h)|\\
&\qquad\le C\bigl(h\varrho+|B|\bigr).
\end{aligned}
\]
By \eqref{ckg:general-measure-final} and \eqref{ckg:explicit-measure-rates-0}-- \eqref{ckg:explicit-measure-rates-frozen}, together with the identical-mode estimate following them,
\[
|B|\le C\varrho^2r_0^{2\varsigma/3}.
\]
As $|Q|=\varrho^2$, this proves \eqref{ckg:common-measure}.

Choosing $r_*$ and $\kappa_0$ sufficiently small makes the right-hand side of \eqref{ckg:common-measure} smaller than one. Thus $f_h^{-1}(\mathcal W_h)$ has positive measure, and so does $\mathcal W_h$ by the bi-Lipschitz property of $f_h$. The parameters in \eqref{ckg:two-parameters-common}
exist and are unique by the definition of $\mathcal W_h$ and the injectivity of $f_h,f_0$. Finally,
\[
\begin{aligned}
|f_0(\rho_h)-f_0(\rho_0)|=|f_0(\rho_h)-f_h(\rho_h)|\le Ch,
\end{aligned}
\]
and the uniform inverse Lipschitz bound for $f_0$ gives \eqref{ckg:rho-distance}.
\end{proof}

\begin{remark}
The hypothesis $h<\kappa_0\varrho$ ensures a common-frequency set of positive measure for sufficiently small fixed $\kappa_0$. The quantitative bound \eqref{ckg:common-measure} yields asymptotically full relative measure in a joint $r_0,h$ limit when, in addition, $h/\varrho\to0$.
\end{remark}

\subsection{Proof of Theorem \ref{thm:intro-CKG}}
\label{sec:final-comparison}

We now recall \emph{all} the changes of coordinates. Write $\iota(\vartheta):=(\vartheta,0,0)$ for the zero section of the KAM phase space (including $W_{\cF}=0$ when $\mathcal F\ne\varnothing$) and recall \eqref{eq:T0NLS-map}, \eqref{eq:T0-map}, \eqref{psi_mn}, \eqref{psi:m}, \eqref{ckg:C3-map} and \eqref{ckg:C3-NLS-map}. For each $\rho\in\Xi_\infty$, define
\begin{equation}\label{ckg:full-reconstruction-maps}
 \begin{aligned}
 \mathfrak B_h(\rho):= \mathcal T_0(h)\circ\Psi (\rho)\circ\mathcal C_3(h;\rho), \qquad 
 \mathfrak B_0(\rho):= \mathcal T_0^{\rm NLS}\circ\Psi (\rho) \circ\mathcal C_3^{\rm NLS}(\rho),
 \end{aligned}
\end{equation}
recalling that $\mathcal C_3^\sharp$ are the identity for $m=n$ by Proposition~\ref{ckg:prop-cubic-elimination} and denoting by $\Psi$ the angle elimination $\Psi_m$ or $\Psi_{m,n}$. The complete embeddings in the diagonal KG variables and in the coupled NLS variables are, respectively,
\begin{equation}\label{ckg:full-torus-embeddings}
 \begin{aligned}
 \Psi_h(\vartheta,\rho) &:=\mathfrak B_h(\rho)\circ\Phi_\infty(\rho)
   \circ\iota(\vartheta),\\
 \Psi_0(\vartheta,\rho) &:=\mathfrak B_0(\rho)\circ\Phi_\infty^{\rm NLS}(\rho) \circ\iota(\vartheta).
 \end{aligned}
\end{equation}
Here $\Phi_\infty$ and $\Phi_\infty^{\rm NLS}$ are those in \eqref{Phi} and include \emph{both} the preliminary map $\Phi_0$ and every Newton map, as specified in Lemma~\ref{ckg:lem-iterative-KAM}.  Thus \eqref{ckg:full-torus-embeddings} is exactly the inverse reconstruction associated with the pull-backs \eqref{ckg:H0-def} and \eqref{ckg:iteration-conjugacy-KG}--\eqref{ckg:iteration-conjugacy-NLS}. In particular the dependence on the reference parameter $\rho$ is retained through the angle maps and the subsequent KAM maps. By Remark~\ref{ckg:lem-autonomous-NLS},the NLS embedding $\Psi_0$ is defined on the whole box $\Xi_0$ and is independent of $h$.

\begin{lemma}[Embedding comparison at common frequency]
\label{ckg:lem-embedding-common}
Under \eqref{ckg:varsigma-choice}--\eqref{ckg:alpha-choices}, for every $\omega\in\mathcal W_h$ and the two parameters of \eqref{ckg:two-parameters-common}, there is a constant $C>0$ independent of $c,r_0,\omega$ such that
\begin{align}
 \sup_{\vartheta\in\mathbb T^2}\|\Psi_h(\vartheta,\rho_h(\omega))-\Psi_0(\vartheta,\rho_0(\omega))\|_{a,p-4}&\le Chr_0^{-3/2},\label{ckg:embedding-p4}\\
 \sup_{\vartheta\in\mathbb T^2}\|\Psi_h(\vartheta,\rho_h(\omega))-\Psi_0(\vartheta,\rho_0(\omega))\|_{a,p}&\le Cr_0^{2/3}.\label{ckg:embedding-p}
\end{align}
Consequently, for $0\le s\le1$,
\begin{equation}\label{ckg:embedding-interpolation}
 \sup_{\vartheta\in\mathbb T^2}\|\Psi_h(\vartheta,\rho_h(\omega))-\Psi_0(\vartheta,\rho_0(\omega))\|_{a,p-4s}\le Ch^s r_0^{\frac23(1-s)-\frac32s}.
\end{equation}
\end{lemma}

\begin{proof}
Fix $\omega\in\mathcal W_h$, write $\rho_h=\rho_h(\omega)$, $\rho_0=\rho_0(\omega)$, and set
\[
 A:=\frac{\epsilon_0}{\alpha_0^2\sigma_0^{2\mu}}+\frac{\epsilon_1}{\alpha_1^2\sigma_1^{2\mu}}.
\]
For $\vartheta\in\mathbb T^2$, put
\[
 w_h:=\Phi_\infty(\rho_h)\circ\iota(\vartheta),\qquad w_0:=\Phi_\infty^{\rm NLS}(\rho_0)\circ\iota(\vartheta).
\]
By \eqref{ckg:general-Phi-size}, both points belong to a fixed smaller KAM domain, uniformly in $\vartheta$ and the parameters. Moreover, \eqref{ckg:general-Phi-diff}, \eqref{ckg:general-Phi-Lip}, and \eqref{ckg:rho-distance}, with $r_{\rm K}=3r_0/4$, give
\[
 \|w_h-w_0\|_{r_0,p-4}\le ChA+CA|\rho_h-\rho_0|\le ChA.
\]

Set
\[
 v_h:=\mathcal C_3(h;\rho_h)(w_h),\qquad v_0:=\mathcal C_3^{\rm NLS}(\rho_0)(w_0).
\]
The bounds \eqref{ckg:C3-close}, for $q=p-4$, and Cauchy estimates on a fixed smaller domain imply
\[
 \|D\mathcal C_3^{\rm NLS}\|_{r_0,p-4\to r_0,p-4}\le1+Cr_0^{1/3},\qquad
 \|\mathcal C_3^{\rm NLS}-\mathrm{Id}\|_{r_0,p-4}^{\rm lip}\le C\frac{r_0^{1/3}}{\varrho}=Cr_0^{-1}.
\]
Splitting the difference first at the common parameter $\rho_h$ and then at the common point $w_0$, we obtain from \eqref{ckg:C3-diff}
\[
\begin{aligned}
 \|v_h-v_0\|_{r_0,p-4}
 &\le Chr_0^{1/3}+(1+Cr_0^{1/3})\|w_h-w_0\|_{r_0,p-4}+Cr_0^{-1}|\rho_h-\rho_0|\\
 &\le Ch(r_0^{1/3}+A+r_0^{-1}).
\end{aligned}
\]

Next put $U_h:=\Psi(\rho_h)(v_h)$ and $U_0:=\Psi(\rho_0)(v_0)$, where $\Psi$ denotes $\Psi_m$ or $\Psi_{m,n}$. By \eqref{ckg:same-AA}, \eqref{ckg:Psi-same}, and \eqref{ckg:Psi-dist}, the scaled norm \eqref{ckg:scaled-X-norm} yields
\[
\begin{aligned}
 \|U_h-U_0\|_{a,p-4}&\le C\left(\varrho^{-1/2}|\rho_h-\rho_0|+\sqrt\varrho\,\|v_h-v_0\|_{r_0,p-4}\right) \\
 &\le Ch\left(\varrho^{-1/2}+\sqrt\varrho(r_0^{1/3}+A+r_0^{-1})\right).
\end{aligned}
\]

By \eqref{ckg:general-Phi-size}, \eqref{ckg:C3-close}, and \eqref{psi:m}--\eqref{psi_mn}, both $U_h$ and $U_0$ belong to a ball of radius $C\sqrt\varrho$ in $\ell^{a,p}$. Applying Theorem~\ref{Teorema0} on a fixed larger ball and using \eqref{sti.tra.1} with $q=p-4$, Cauchy estimates give
\[
 \sup_{\|u\|_{a,p-4}\le C\sqrt\varrho}
 \|D\mathcal T_0^{\rm NLS}(u)\|_{\mathcal L(\ell^{a,p-4})}\le1+C\varrho.
\]
Thus \eqref{res.tr.0} and \eqref{ckg:full-torus-embeddings} imply
\[
\begin{aligned}
 \|\Psi_h(\vartheta,\rho_h)-\Psi_0(\vartheta,\rho_0)\|_{a,p-4} &\le\|(\mathcal T_0(h)-\mathcal T_0^{\rm NLS})(U_h)\|_{a,p-4} +\|\mathcal T_0^{\rm NLS}(U_h)-\mathcal T_0^{\rm NLS}(U_0)\|_{a,p-4}\\
 &\le Ch\varrho^{3/2}+(1+C\varrho)\|U_h-U_0\|_{a,p-4}\\
 &\le Ch\left(\varrho^{-1/2}+\sqrt\varrho\,A\right).
\end{aligned}
\]

By \eqref{ckg:varsigma-choice}--\eqref{ckg:alpha-choices}, \eqref{ckg:explicit-smallness-0}--\eqref{ckg:explicit-smallness-1}, and the fact that $\sigma_0,\sigma_1$ are fixed,
\[
 A\le
 \begin{cases}
 Cr_0^{-4/3-6\varsigma},&m=n,\\
 Cr_0^{-17/9-3\varsigma},&m\ne n.
 \end{cases}
\]
Since $\varrho=r_0^{4/3}$ and $0<\varsigma<1/18$, we conclude that
\[
 \sup_{\vartheta\in\mathbb T^2}\|\Psi_h(\vartheta,\rho_h)-\Psi_0(\vartheta,\rho_0)\|_{a,p-4}
 \le
 \begin{cases}
 Chr_0^{-2/3-6\varsigma},&m=n,\\
 Chr_0^{-11/9-3\varsigma},&m\ne n
 \end{cases}
 \le Chr_0^{-3/2},
\]
proving \eqref{ckg:embedding-p4}.

For the strong norm, \eqref{ckg:general-Phi-size}, \eqref{ckg:C3-close}, \eqref{ckg:Psi-same}--\eqref{ckg:Psi-dist}, and \eqref{sti.tra.1} give
\[
 \sup_{\vartheta\in\mathbb T^2}\left(\|\Psi_h(\vartheta,\rho_h)\|_{a,p}+\|\Psi_0(\vartheta,\rho_0)\|_{a,p}\right)\le C\sqrt\varrho=Cr_0^{2/3}.
\]
The triangle inequality proves \eqref{ckg:embedding-p}. Finally, the interpolation inequality
\[
 \|{\rm f}\|_{a,p-4s}\le\|{\rm f}\|_{a,p-4}^s\|{\rm f}\|_{a,p}^{1-s},\qquad 0\le s\le1,
\]
gives \eqref{ckg:embedding-interpolation}.
\end{proof}

For $\rho\in\Xi_\infty$, the formulas \eqref{ckg:Psi-same}--\eqref{ckg:Psi-dist} give $\Psi(\rho)\circ\iota(\vartheta)=\mathcal I_{m,n}(\vartheta,\rho)$, where $\mathcal I_{m,n}$ is defined in \eqref{eq:intro-reference-torus}. On the smaller KAM domains, the same formulas and the scaled norm \eqref{ckg:scaled-X-norm} yield
\[
 \|\Psi(\rho)(w)-\Psi(\rho)(w')\|_{a,p}\le C\sqrt\varrho\,\|w-w'\|_{r_0,p}.
\]
Thus \eqref{ckg:general-Phi-size}, \eqref{ckg:C3-close}, and \eqref{sti.tra.1}, imply
\begin{equation}\label{ckg:individual-torus-nearness}
\begin{aligned}
\sup_{\rho\in\Xi_\infty}&\sup_{\vartheta\in\mathbb T^2}\Bigl(\|\Psi_h(\vartheta,\rho)-\mathcal I_{m,n}(\vartheta,\rho)\|_{a,p}+\|\Psi_0(\vartheta,\rho)-\mathcal I_{m,n}(\vartheta,\rho)\|_{a,p}\Bigr) \\
&\le C\left[\varrho^{3/2}+\sqrt\varrho\left(r_0^{1/3}+\frac{\epsilon_0}{\alpha_0\sigma_0^\mu}+\frac{\epsilon_1}{\alpha_1\sigma_1^\mu}\right)\right]\le Cr_0^{7/9-2\varsigma}.
\end{aligned}
\end{equation}
Here the last inequality follows from \eqref{ckg:varrho-r0} and \eqref{ckg:explicit-smallness-0}--\eqref{ckg:explicit-smallness-1}; the cubic maps are the identity when $m=n$.

The NLS normal form and its generating Hamiltonians commute with the two local masses defined in \eqref{eq:NQM-zeta}. Their flows, the angle maps, and the Birkhoff and cubic maps therefore intertwine the corresponding phase actions. With the conventions \eqref{eq:Hamiltonian-convention} and \eqref{ckg:XF-KAM}, set
\begin{equation*}
 \mathscr R_\phi(\xi,\eta):= (e^{-\iu\phi}\xi,e^{-\iu\phi}\eta).
\end{equation*}
The complete NLS embedding, including all the maps in \eqref{ckg:full-torus-embeddings}, satisfies the exact equivariance
\begin{equation}\label{ckg:NLS-equivariance}
 \Psi_0(\vartheta-\phi(1,1),\rho) =\mathscr R_\phi\Psi_0(\vartheta,\rho).
\end{equation}
The point of retaining this identity through the full composition is that the relativistic rest-energy rotation can be removed without estimating a secular phase error.

\begin{proof}[Conclusion of the proof of Theorem~\ref{thm:intro-CKG}]
The choice of parameters \eqref{ckg:varsigma-choice}--\eqref{ckg:alpha-choices} satisfies the hypotheses of Theorem~\ref{TeoremaGenerale} by \eqref{ckg:explicit-smallness-0}--\eqref{ckg:explicit-smallness-1} and the Melnikov persistence estimate following them. Theorem~\ref{ckg:thm-common-frequency} gives $\mathcal W_h$, the inverse frequency maps $\rho_h,\rho_0$, and \eqref{eq:intro-action-comparison}--\eqref{eq:intro-common-measure}, by \eqref{ckg:two-parameters-common}, \eqref{ckg:rho-distance}, and \eqref{ckg:common-measure}. The maps $\rho_h,\rho_0$ are bi-Lipschitz as restrictions of the inverses of the frequency maps in Theorem~\ref{TeoremaGenerale}.

The full embeddings are given by \eqref{ckg:full-torus-embeddings}. Their approximation \eqref{eq:intro-torus-closeness} to the two-mode tori follows by \eqref{ckg:individual-torus-nearness}.The stability statements follow from Theorem~\ref{TeoremaGenerale}: the symplectic reconstruction preserves the normal linearized dynamics, and the finite hyperbolic block in the case $m\ne n$ satisfies \eqref{ckg:hyp-resolvent-allt}. Finally, write $\Theta(t):=\vartheta_0-(c^2(1,1)+\omega)t$. The NLS equivariance \eqref{ckg:NLS-equivariance} gives
\[
 \mathscr R_{-c^2t}\Psi_0(\Theta(t),\rho_0(\omega))=\Psi_0(\vartheta_0-\omega t,\rho_0(\omega)).
\]
Since $\mathscr R_\phi$ is an isometry in every $\ell^{a,q}$, \eqref{ckg:embedding-interpolation} proves \eqref{eq:intro-uniform-comparison} uniformly for $t\in\mathbb R$. Reconstruction \eqref{eq:xi-eta-inverse} and the two-carrier definition \eqref{eq:intro-two-carrier} give \eqref{eq:intro-physical-comparison} by the triangle inequality, using \eqref{eq:intro-uniform-comparison} and the fact that$I-A_c^{-1}$ has symbols bounded by $h^s \langle j \rangle^{2s}$.
\end{proof}

\subsection{Proof of Theorem \ref{cor:intro-CKG-subsequence}}\label{proofcorollary}

Put $Q:=\Xi_0=[\varrho/2,3\varrho/2]^2$. The uniform inverse Lipschitz bound for the frequency maps, together with \eqref{ckg:common-measure}, gives
\[
|\mathcal W_{h_j}|\ge L^{-2}\varrho^2\left(1-Cr_0^{2\varsigma/3}-C\frac{h_j}{\varrho}\right).
\]
All sets $\mathcal W_{h_j}$ are measurable and contained in a fixed bounded region of slow frequencies. Since $\mathcal W_* =\bigcap_{N\ge1}\bigcup_{j\ge N}\mathcal W_{h_j}$, continuity of finite Lebesgue measure from above yields
\[
|\mathcal W_*|=\lim_{N\to\infty}\left|\bigcup_{j\ge N}\mathcal W_{h_j}\right|\ge\limsup_{j\to\infty}|\mathcal W_{h_j}|.
\]
As $h_j\to0$, this proves \eqref{eq:intro-fixed-NLS-measure}.

By the definition of $\mathcal W_*$, for any element $\omega$ there exists a subsequence $j_k\to\infty$ such that $\omega\in\mathcal W_{h_{j_k}}$ for every $k$. By Remark~\ref{ckg:lem-autonomous-NLS}, the NLS frequency map $f_0=\omega_\infty^{\rm NLS}$ and the full embedding $\Psi_0$ are independent of $h$.
Since
\[
f_0\bigl(\rho_0^{(h_{j_k})}(\omega)\bigr)=\omega,
\]
injectivity of $f_0$ gives
\[
\rho_0^{(h_{j_k})}(\omega) =\rho_\omega:=f_0^{-1}(\omega)
\]
for every $k$. Define the fixed NLS embedding $\Psi_\omega(\vartheta):=\Psi_0(\vartheta,\rho_\omega)$ and, for a fixed initial tangential angle $\vartheta_0$, the corresponding solution $(\xi_0(t),\eta_0(t)):=\Psi_\omega(\vartheta_0-\omega t).$ The identification of the simultaneous and autonomous NLS embeddings on the common Cantor sets, combined with \eqref{ckg:embedding-interpolation}, gives, for every $0<s\le1$,
\begin{equation}\label{ckg:fixed-NLS-embedding-limit}
\lim_{k\to\infty}\sup_{\vartheta\in\mathbb T^2}\|\Psi_{h_{j_k}}(\vartheta,\rho_{h_{j_k}}(\omega))-\Psi_\omega(\vartheta)\|_{a,p-4s}=0,
\end{equation}
The exact equality of the slow frequencies and the NLS equivariance \eqref{ckg:NLS-equivariance} then yield \eqref{eq:intro-fixed-NLS-limit} from \eqref{ckg:fixed-NLS-embedding-limit}, by the same argument used in the proof of Theorem~\ref{thm:intro-CKG}. Finally, using the reconstruction \eqref{eq:xi-eta-inverse} and \eqref{eq:intro-physical-comparison} we also obtain the physical-field convergence
\begin{equation}\label{ckg:fixed-NLS-two-carrier-limit}
\lim_{k\to\infty}\sup_{t\in\mathbb R}\left\|u_{h_{j_k}}(t)-\left(e^{-\iu c_{j_k}^2t}\xi_0(t)+e^{\iu c_{j_k}^2t}\overline{\eta_0(t)}\right)\right\|_{a,p-4s}=0,\qquad 0<s\le1.
\end{equation}

\begin{remark}[Generalizations of the limit approach]\label{gen}
The independence of the NLS family from \(h\) presented in Remark~\ref{ckg:lem-autonomous-NLS} does not rely on the NLS Melnikov non-resonance conditions being satisfied on the entire parameter domain. More generally, one may first construct the NLS KAM $\nu$-iteration on its own $h$-independent Cantor sets $\Xi_\nu^{\rm NLS}$ over which Melnikov conditions are satisfied, and then perform the simultaneous KG--NLS iteration on a common set $\Xi_\nu(h)\subseteq\Xi_\nu^{\rm NLS}$. Uniqueness of the normalized homological solutions identifies the NLS transformations on these sets with restrictions of the autonomous ones. One may then extend the NLS frequency corrections from $\Xi_\nu^{\rm NLS}$ and the KG--NLS differences from$\Xi_\nu(h)$, as in Remark~\ref{ckg:lem-autonomous-NLS}, while preserving the $h$-independence of the NLS frequency map and the global $O(h)$ comparison.. This observation is expected to be relevant when exciting more tangential modes, since additional NLS quasi-resonances may arise. Finally, combining this remark with the KAM construction of the real cubic KG in \cite{BBG} we strongly believe that the approach presented in this remark could be applied also to the real model and its KAM scheme, obtaining a subsequential nonrelativistic limit towards a single NLS solution, similarly to the result presented in Theorem~\ref{cor:intro-CKG-subsequence}.
\end{remark}

\appendix
\refstepcounter{section}
\section*{Appendix A: Technical results}
\phantomsection
\addcontentsline{toc}{section}{Appendix A}
\addtocontents{toc}{\protect\setcounter{tocdepth}{1}}
\subsection{Convolution and multilinear analytic estimates}
For scalar sequences on $\mathbb Z$ we use the standard convolution $(u\star v)_j:={\sum}_{k\in\mathbb Z}u_{j-k}v_k.$ For $a\ge0$, $p>1/2$ and $\beta\ge0$ by Lemma A.1 in \cite{BBG} it satisfies
\begin{equation}\label{algebra.unif}
 \|u\star v\|_{a,p,\beta} \le C_{a,p,\beta}\|u\|_{a,p,\beta}\|v\|_{a,p,\beta},
\end{equation}
with a constant independent of $c\ge1$.

\medskip

In particular, for $n\ge2$ and scalar sequences $z^{(k)}$ we write $\bigstar_{k=1}^n z^{(k)}:=z^{(1)}\star\cdots\star z^{(n)}$; associativity and \eqref{algebra.unif} yield
\begin{equation}\label{ckg:multiconvolution-uniform}
 \left\|\mathop{\bigstar}_{k=1}^n z^{(k)}\right\|_{a,p,\beta} \le C_{n,a,p,\beta}\prod_{k=1}^n\|z^{(k)}\|_{a,p,\beta}.
\end{equation}

For $r=1,\ldots,n$ let $b^{(r)}=(b_j^{(r)})_{j\in\mathbb Z}$ be a positive diagonal Fourier multiplier. For our space of indexes $\cZ=\mathbb{Z} \times \{\pm 1\}$ we have the following adaptation of Lemma A.2 in \cite{BBG}.

\begin{lemma}[Multilinear Hamiltonian estimate]\label{sti.abs.2}
Let $n\ge2$ and let $F$ be a homogeneous Hamiltonian of degree $n$, initially defined on finitely supported variables, of the form
\begin{equation}\label{sti.abs.1}
 F(\zeta,\bar\zeta) =\sum_{\boldsymbol\alpha,\boldsymbol\sigma} F_{\boldsymbol\alpha,\boldsymbol\sigma} \prod_{r=1}^n b_{\jmath(\alpha_r)}^{(r)} \zeta_{\alpha_r}^{\sigma_r},
\end{equation}
where $|F|_\infty:=\sup_{\boldsymbol\alpha,\boldsymbol\sigma} |F_{\boldsymbol\alpha,\boldsymbol\sigma}|<\infty$ and $F_{\boldsymbol\alpha,\boldsymbol\sigma}=0$ unless $\sum_{r=1}^n\sigma_r\jmath(\alpha_r)=0$. For each $r$ define 
\begin{equation}\label{wk}
 \begin{aligned}
 w_j^{(r)}:=b_j^{(r)}\sum_{s\in\{\pm1\}}
 (|\zeta_{(j,s)}|+|\bar\zeta_{(j,s)}|), \qquad  W_j^{(r)}:=w_j^{(r)}+w_{-j}^{(r)},\qquad j\in\mathbb Z.
 \end{aligned}
\end{equation}
Then, on finitely supported variables,
\begin{equation}\label{sti.abs.6}
 \|X_F(\zeta,\bar\zeta)\|_{a,p,\beta} \le C_n|F|_\infty\sum_{r=1}^n \left\|b^{(r)} \left(\mathop{\bigstar}_{k\ne r}W^{(k)}\right) \right\|_{a,p,\beta},
\end{equation}
where $C_n$ is independent of $c$ and of the multipliers. The same estimate and an analytic extension hold on any open domain of the phase space on which the right-hand side is locally bounded.
\end{lemma}

\subsection{Quartic resonances and small divisors}

\begin{lemma}[Exact quartic resonances]\label{APP:lem:quartic-resonances}
For $(\boldsymbol\alpha,\boldsymbol\sigma)\in \mathcal Z^4\times\{\pm1\}^4$, the conditions
\begin{equation}\label{APP:eq:res-assumptions}
\sum_{r=1}^4\sigma_r\jmath(\alpha_r)=0, \qquad \Omega(\boldsymbol\alpha,\boldsymbol\sigma)=0
\end{equation}
hold if and only if $(\boldsymbol\alpha,\boldsymbol\sigma) \in\mathcal{IR}$ (see Definition~\ref{resonant}). If the corresponding monomial is physically gauge invariant, it is also locally gauge invariant. Indeed, after reordering so that $j_1=j_2$, $j_3=j_4$ and $\boldsymbol\sigma=(+,-,+,-)$, the species labels satisfy either
\begin{align}
&s_1=s_2,\qquad s_3=s_4,
\label{APP:eq:res-action-type}\\
&(s_1,s_2,s_3,s_4)\in
\{(+,-,-,+),(-,+,+,-)\}.
\label{APP:eq:res-exchange-type}
\end{align}
The first case gives action monomials, while the second gives exchange monomials
\begin{equation}\label{APP:eq:exchange-monomial}
\xi_j\bar\eta_j\eta_k\bar\xi_k \quad\text{and its complex conjugate}.
\end{equation}
\end{lemma}

\begin{proof}
Since the species labels do not enter the divisor, the first assertion follows from the scalar resonance classification \cite{berti2013kam}. Under the above reordering, physical gauge invariance \eqref{GaugeG} reads $s_1-s_2+s_3-s_4=0$, whose solutions for $s_r\in\{\pm1\}$ are precisely \eqref{APP:eq:res-action-type} and \eqref{APP:eq:res-exchange-type}. Both types have vanishing local species charges.
\end{proof}

\begin{lemma}[Quartic small divisors]\label{lem:small-divisors-CKG}
Let $J$ be the set in \eqref{eq:Jfourier} and $\Omega, \Omega^{\rm NLS}$ as in \eqref{eq:OmegaKG}. Then there exists a constant $\gamma_J>0$ such that, for every $c>1$,
\begin{equation}\label{eq:small-divisor-KG}
  |\Omega(\boldsymbol\alpha,\boldsymbol\sigma)| \ge \gamma_J
\end{equation}
whenever $(\boldsymbol\alpha,\boldsymbol\sigma)\in \mathcal L_{\mathcal J}\setminus\mathcal{IR}$ (see definition \ref{resonant}). If moreover $\sum_r\sigma_r\ne0$, then
\begin{equation}\label{eq:small-divisor-NG}
  |\Omega(\boldsymbol\alpha,\boldsymbol\sigma)| \ge \gamma_J c^2.
\end{equation}
The analogous NLS divisor is non-zero integer on $\mathcal L_{\mathcal J}\setminus\mathcal{IR}$ under the local gauge condition; in particular $|\Omega^{\rm NLS}|\ge1/2$.
\end{lemma}

\begin{proof}
After forgetting the species labels, the Klein--Gordon estimates \eqref{eq:small-divisor-KG}--\eqref{eq:small-divisor-NG} are exactly those used in \cite{BBG}, following \cite{berti2013kam}. Indeed, $\mathcal L_{\mathcal J}$ contains a Fourier index in the fixed set $J$, and Lemma~\ref{APP:lem:quartic-resonances} gives the same scalar resonances; the species labels do not affect the divisor. For NLS, local gauge invariance gives $\sum_r\sigma_r=0$. After permuting the signs, momentum and a vanishing divisor would give $j_1+j_2=j_3+j_4$ and $j_1^2+j_2^2=j_3^2+j_4^2$, hence $(j_1-j_3)(j_1-j_4)=0$. This is a paired configuration in $\mathcal{IR}$. Outside $\mathcal{IR}$ the divisor is therefore a nonzero half-integer.
\end{proof}

\begin{lemma}[Lipschitz extensions]\label{ckg:lem-extension}
Let $(X,d)$ be a compact metric space and $E\subset X$ be nonempty and closed. For $\lambda\ge0$, write $\|u\|_E^{(\lambda)}:=\|u\|_{\infty,E}+\lambda\operatorname{Lip}_E(u)$. Every real-valued Lipschitz function $f$ on $E$ admits an extension $F$ to $X$ preserving its minimum, maximum and Lipschitz seminorm. Moreover, if $f,g:E\to\mathbb R$ are Lipschitz, there exist extensions $F,G$ such that
\[
\|F-G\|_{\infty,X}\le\|f-g\|_{\infty,E},\qquad \|F\|_X^{(\lambda)}=\|f\|_E^{(\lambda)},\qquad \|G\|_X^{(\lambda)}\le2\|f\|_E^{(\lambda)}+\|g\|_E^{(\lambda)}
\]
for every $\lambda\ge0$. Let $V\subset\mathbb C^{m\times m}$ be a real linear subspace, with $m\le m_*$, endowed with the operator norm. Every Lipschitz map $A:E\to V$ admits an extension $\widetilde A:X\to V$ satisfying
\[
\|\widetilde A\|_{\infty,X}+\lambda\operatorname{Lip}_X(\widetilde A)\le C_{m_*}\bigl(\|A\|_{\infty,E}+\lambda\operatorname{Lip}_E(A)\bigr),\qquad \lambda\ge0.
\]
The constant depends only on $m_*$ and is uniform for families of such matrix spaces. The same estimate holds for weighted families with arbitrary positive scalar weights, one for each matrix block, independent of $x\in X$.
\end{lemma}

\begin{proof}
For a real-valued Lipschitz function $f$, by \cite[Lemma~A.2]{poschel1996kam} there exists $F_0: X \to \mathbb{R}$ preserving the Lipschitz seminorm and the supremum norm. For the paired assertion, apply this construction to $f$ and $r:=g-f$, obtaining $F$ and $R$, and set $G:=F+R$. Then $G|_E=g$, $\|F-G\|_{\infty,X}=\|r\|_{\infty,E}$ and $\|G\|_X^{(\lambda)}\le\|f\|_E^{(\lambda)}+\|r\|_E^{(\lambda)}\le2\|f\|_E^{(\lambda)}+\|g\|_E^{(\lambda)}$.

For the matrix-valued assertion, choose an orthonormal basis of $V$ with respect to the real Hilbert--Schmidt inner product $\operatorname{Re}\operatorname{Tr}(A^*B)$ and extend the corresponding real coordinate functions by the scalar construction. Reconstructing the matrix preserves membership in $V$. Since $\dim_{\mathbb R}V\le2m^2$ and $\|A\|\le\|A\|_{\rm HS}\le\sqrt m\,\|A\|$, the supremum and Lipschitz bounds follow with a constant depending only on $m_*$. For weighted families, apply the same construction to each matrix multiplied by its fixed positive scalar weight, and divide by that weight afterwards.
\end{proof}

\subsection{Spectral and sparse-block estimates}

\begin{lemma}
\label{ckg:lem-Weyl-branches}
Let $A(\rho)$ be Hermitian of dimension at most four and let $\mu_1(\rho)\le\cdots\le\mu_d(\rho)$ be its ordered eigenvalues.  Then
\begin{equation}\label{ckg:Weyl-Lip-full}
 |\mu_a(\rho)-\mu_a(\rho')| \le\|A(\rho)-A(\rho')\|, \qquad a=1,\dots,d.
\end{equation}
If $A=A^{\rm ref}+D$, then, after ordering both spectra,
\begin{equation}\label{ckg:Weyl-pointwise-ref}
 |\mu_a(A(\rho))-\mu_a(A^{\rm ref}(\rho))| \le\|D(\rho)\|.
\end{equation}
\end{lemma}

\begin{remark}
By \eqref{ckg:Weyl-Lip-full}, each ordered eigenvalue branch of a Lipschitz Hermitian family is Lipschitz, even at crossings. However, \eqref{ckg:Weyl-pointwise-ref} does not imply a bound on $\operatorname{Lip}(\mu_a(A)-\mu_a(A^{\rm ref}))$ in terms of $\operatorname{Lip}(D)$ alone; such a bound may fail at crossings. Accordingly, the weighted tail Lipschitz estimates used below are imposed at matrix level, while Weyl's inequality transfers matrix tails to pointwise spectral tails.
\end{remark}

\begin{lemma}[Sparse block Schur test]
\label{ckg:lem-sparse-block-Schur}
Let $\mathfrak B$ be a countable family of finite-dimensional Hilbert spaces $(E_B)_{B\in\mathfrak B}$ and let
\[
 M=(M_{BC})_{B,C\in\mathfrak B}, \qquad M_{BC}\in\mathcal L(E_C,E_B).
\]
Assume that there exist integers $D_{\rm row},D_{\rm col}\ge1$ such that
\begin{equation}\label{ckg:Schur-sparsity}
 \#\{C:M_{BC}\neq0\}\le D_{\rm row}, \qquad \#\{B:M_{BC}\neq0\}\le D_{\rm col},
\end{equation}
uniformly in $B,C$. Let $(u_B)_{B\in\mathfrak B}$ and $(v_B)_{B\in\mathfrak B}$ be positive weights and define
\[
 \|x\|_v^2 :=\sum_{C\in\mathfrak B}v_C^2|x_C|_{E_C}^2, \qquad \|y\|_u^2 :=\sum_{B\in\mathfrak B}u_B^2|y_B|_{E_B}^2.
\]
If
\begin{equation}\label{ckg:Schur-block-bound}
 A:= \sup_{\substack{B,C\in\mathfrak B\\M_{BC}\neq0}} \frac{u_B}{v_C}\, \|M_{BC}\|_{\mathcal L(E_C,E_B)} <\infty,
\end{equation}
then $M$ defines a bounded operator
\begin{equation}\label{ckg:Schur-conclusion}
 M:\ell^2_v((E_C))\longrightarrow\ell^2_u((E_B)), \qquad \|M\|_{\ell^2_v\to\ell^2_u} \le A\sqrt{D_{\rm row}D_{\rm col}}.
\end{equation}

In particular, if $B=[i]$, $C=[j]$ and $u_B=e^{a|i|}\langle i\rangle^q\widetilde w_B^\beta,$ and $v_C=e^{a|j|}\langle j\rangle^q,$ then
\begin{equation}\label{ckg:Schur-weighted-specialization}
 \|M\|_{\ell^{a,q}\to\ell^{a,q,\beta}} \le \sqrt{D_{\rm row}D_{\rm col}}\, \sup_{B,C} \|M_{BC}\|_{q,\beta;B\leftarrow C}.
\end{equation}
For $\beta=0$ this gives the corresponding estimate on $\ell^{a,q}$.
\end{lemma}

\begin{proof}
Conjugating by the weights reduces the estimate to the nonnegative scalar matrix $a_{BC}:=(u_B/v_C)\|M_{BC}\|$. By \eqref{ckg:Schur-sparsity}--\eqref{ckg:Schur-block-bound},
\[
 \sup_B\sum_Ca_{BC}\le AD_{\rm row},\qquad \sup_C\sum_Ba_{BC}\le AD_{\rm col}.
\]
The scalar Schur test gives \eqref{ckg:Schur-conclusion}. For the weights in the statement, $(u_B/v_C)\|M_{BC}\|=\|M_{BC}\|_{q,\beta;B\leftarrow C}$, which gives \eqref{ckg:Schur-weighted-specialization}.
\end{proof}

\begin{lemma}[Operator monotonicity and strip measure]\label{operatori:hem}

Let $I\subset\mathbb R$ be an interval and let $H,E:I\to\operatorname{Herm}(d)$, with $d$ uniformly bounded.  For Hermitian matrices we write $A\ge B$ if $A-B$ is positive semidefinite, equivalently if $\langle\bar u,(A-B)u\rangle\ge0$ for every $u\in\mathbb C^d$. Assume that
\begin{equation}\label{ckg:Loewner-transversality}
 H(t)-H(t')\ge m(t-t')I, \qquad t>t',
\end{equation}
for some $m>0$, and that
\begin{equation*}
 \|E(t)-E(t')\|\le\eta|t-t'|, \qquad \eta<m.
\end{equation*}
If $\lambda_1(t)\le\cdots\le\lambda_d(t)$ are the ordered eigenvalues of $H(t)+E(t)$, then for $t>t'$,
\begin{equation}\label{ckg:eigen-increment-monotone}
 \lambda_a(t)-\lambda_a(t') \ge(m-\eta)(t-t'), \qquad a=1,\ldots,d.
\end{equation}
If $-H$ satisfies \eqref{ckg:Loewner-transversality}, the same conclusion holds with the opposite sign.  Consequently, for every $\delta>0$,
\begin{equation}\label{ckg:one-dimensional-strip}
 \left|\left\{t\in I: \min_{1\le a\le d}|\lambda_a(t)|<\delta\right\}\right| \le \frac{2d\,\delta}{m-\eta}.
\end{equation}
The same statement applies to a Hermitian Sylvester operator after identifying its finite-dimensional coefficient space with a Hilbert space endowed with the Hilbert--Schmidt inner product.
\end{lemma}

\begin{proof}
For $t>t'$, the perturbation bound implies $E(t)-E(t')\ge-\eta(t-t')I$. Hence \eqref{ckg:Loewner-transversality} gives
\[
 [H(t)+E(t)]-[H(t')+E(t')]\ge(m-\eta)(t-t')I.
\]
The min--max principle yields \eqref{ckg:eigen-increment-monotone}; the opposite orientation gives the corresponding decreasing inequality. In either case each ordered eigenvalue is strictly monotone, and the set $\{t\in I:|\lambda_a(t)|<\delta\}$ has diameter at most $2\delta/(m-\eta)$. Summing over $a=1,\ldots,d$ proves \eqref{ckg:one-dimensional-strip}.
\end{proof}


\begin{thebibliography}{99}
\small

\bibitem[BBHM18]{baldiBertiHausMontalto2018WW}
P.~Baldi, M.~Berti, E.~Haus and R.~Montalto,
\newblock Time quasi-periodic gravity water waves in finite depth,
\newblock \emph{Invent. Math.} \textbf{214} (2018), 739--911.

\bibitem[BBM14]{baldiBertiMontalto2014Airy}
P.~Baldi, M.~Berti and R.~Montalto,
\newblock KAM for quasi-linear and fully nonlinear forced perturbations of Airy equation,
\newblock \emph{Math. Ann.} \textbf{359} (2014), 471--536.

\bibitem[BBM16]{baldiBertiMontalto2016KdV}
P.~Baldi, M.~Berti and R.~Montalto,
\newblock KAM for autonomous quasi-linear perturbations of KdV,
\newblock \emph{Ann. Inst. H. Poincar\'e Anal. Non Lin\'eaire} \textbf{33} (2016), 1589--1638.

\bibitem[BM21]{baldiMontalto2021Euler}
P.~Baldi and R.~Montalto,
\newblock Quasi-periodic incompressible Euler flows in 3D,
\newblock \emph{Adv. Math.} \textbf{384} (2021), 107730.

\bibitem[BBG26]{BBG}
D.~Bambusi, A.~Belloni and F.~Giuliani,
\newblock Non-relativistic limit of the nonlinear Klein-Gordon equation:
uniform in time approximation of KAM solutions,
\newblock preprint, arXiv:2501.17691v2; \emph{to appear in Archive for Rational Mechanics and Analysis}.

\bibitem[BG01]{bambusiGraffi2001}
D.~Bambusi and S.~Graffi,
\newblock Time quasi-periodic unbounded perturbations of Schr\"odinger operators and KAM methods,
\newblock \emph{Comm. Math. Phys.} \textbf{219} (2001), 465--480.

\bibitem[BBM11]{bambusiBertiMagistrelli2011}
D.~Bambusi, M.~Berti and E.~Magistrelli,
\newblock Degenerate KAM theory for partial differential equations,
\newblock \emph{J. Differential Equations} \textbf{250} (2011), 3379--3397.

\bibitem[BLZ24]{baoLuZhang2024}
W.~Bao, Y.~Lu and Z.~Zhang,
\newblock Convergence rates in the nonrelativistic limit of the
cubic Klein--Gordon equation,
\newblock \emph{SIAM J. Math. Anal.} \textbf{56} (2024),
6822--6860.

\bibitem[BMS04]{bechouche2004nonrelativistic}
P.~Bechouche, N.~J.~Mauser and S.~Selberg,
\newblock Nonrelativistic limit of Klein--Gordon--Maxwell
to Schr\"odinger--Poisson,
\newblock \emph{Amer. J. Math.} \textbf{126} (2004),
31--64.

\bibitem[BGR24]{bernier2024infinite}
J.~Bernier, B.~Gr\'ebert and T.~Robert,
\newblock Infinite dimensional invariant tori for nonlinear
Schr\"odinger equations,
\newblock preprint, arXiv:2412.11845v2 (2024); \emph{to appear in Memoirs of the European Mathematical Society.}

\bibitem[BP26]{bernierProcesi2026}
J.~Bernier and M.~Procesi,
\newblock Lagrangian invariant tori for completely
resonant nonlinear Schr\"odinger equations,
\newblock preprint, arXiv:2610.26044v1 (2026).

\bibitem[BB11]{bertiBiasco2011}
M.~Berti and L.~Biasco,
\newblock Branching of Cantor manifolds of elliptic tori
and applications to PDEs,
\newblock \emph{Comm. Math. Phys.} \textbf{305} (2011),
741--796.

\bibitem[BBP13]{berti2013kam}
M.~Berti, L.~Biasco and M.~Procesi,
\newblock KAM theory for the Hamiltonian derivative wave equation,
\newblock \emph{Ann. Sci. \`Ec. Norm. Sup\'er.} \textbf{46} (2013), 301-373.

\bibitem[BB13]{bertiBolle2013}
M.~Berti and P.~Bolle,
\newblock Quasi-periodic solutions with Sobolev regularity of NLS on $\mathbb T^d$ with a multiplicative potential,
\newblock \emph{J. Eur. Math. Soc.} \textbf{15} (2013), 229--286.

\bibitem[BFM21]{bertiFranzoiMaspero2021WW}
M.~Berti, L.~Franzoi and A.~Maspero,
\newblock Traveling quasi-periodic water waves with constant vorticity,
\newblock \emph{Arch. Ration. Mech. Anal.} \textbf{240} (2021), 99--202.

\bibitem[BFPT24]{bertiFeolaProcesiTerracina2024}
M.~Berti, R.~Feola, M.~Procesi and S.~Terracina,
\newblock Reducibility of Klein--Gordon equations with maximal order perturbations,
\newblock preprint:  arXiv:2402.11377v1 (2024); \emph{to appear in Mem. Eur. Math. Soc.}.

\bibitem[BFM24]{bertiFranzoiMaspero2024WW}
M.~Berti, L.~Franzoi and A.~Maspero,
\newblock Pure gravity traveling quasi-periodic water waves with constant vorticity,
\newblock \emph{Comm. Pure Appl. Math.} \textbf{77} (2024), 990--1064.

\bibitem[BKM18]{bertiKappelerMontalto2018}
M.~Berti, T.~Kappeler and R.~Montalto,
\newblock Large KAM tori for perturbations of the defocusing NLS equation,
\newblock \emph{Ast\'erisque} \textbf{403} (2018),
1--160.

\bibitem[BKM21]{bertiKappelerMontalto2021KdV}
M.~Berti, T.~Kappeler and R.~Montalto,
\newblock Large KAM tori for quasi-linear perturbations of KdV,
\newblock \emph{Arch. Ration. Mech. Anal.} \textbf{239} (2021), 1395--1500.

\bibitem[BM20]{bertiMontalto2020WW}
M.~Berti and R.~Montalto,
\newblock Quasi-periodic standing wave solutions of gravity-capillary water waves,
\newblock \emph{Mem. Amer. Math. Soc.} \textbf{263} (2020), no.~1273.

\bibitem[BFMT25]{bianchiniFranzoiMontaltoTerracina2025}
R.~Bianchini, L.~Franzoi, R.~Montalto and S.~Terracina,
\newblock Large amplitude quasi-periodic traveling waves in two dimensional forced rotating fluids,
\newblock \emph{Comm. Math. Phys.} \textbf{406(3)} (2025), 66.

\bibitem[BMP23]{biasco2023small}
L.~Biasco, J.~E.~Massetti and M.~Procesi,
\newblock Small amplitude weak almost periodic solutions for the 1D NLS,
\newblock \emph{Duke Math. J.} \textbf{172} (2023), 2643--2714.

\bibitem[BCGP24]{corsi2024almost}
L.~Biasco, L.~Corsi, G.~Gentile and M.~Procesi,
\newblock Asymptotically full measure sets of almost-periodic solutions for the NLS equation,
\newblock preprint, arXiv:2412.02648v2 (2024).

\bibitem[Bou98]{bourgain1998Hamiltonian}
J.~Bourgain,
\newblock Quasi-periodic solutions of Hamiltonian perturbations of 2D linear Schr\"odinger equations,
\newblock \emph{Ann. of Math.} \textbf{148} (1998), 363--439.

\bibitem[CY00]{chierchiaYou2000}
L.~Chierchia and J.~You,
\newblock KAM tori for 1D nonlinear wave equations with periodic boundary conditions,
\newblock \emph{Comm. Math. Phys.} \textbf{211(2)} (2000),
497--525.

\bibitem[CLW25]{cong2025existencefulldimensionalkamtori}
H.~Cong, S.~Li and X.~Wu,
\newblock The existence of full-dimensional KAM tori for one-dimensional nonlinear Klein--Gordon equation,
\newblock preprint, arXiv:2505.06526v2 (2025).

\bibitem[CY21]{congYuan2021}
H.~Cong and X.~Yuan,
\newblock The existence of full dimensional invariant tori
for 1-dimensional nonlinear wave equation,
\newblock \emph{Ann. Inst. H. Poincar\'e Anal. Non Lin\'eaire}
\textbf{Vol. 38.No. 3.} (2021), 759--786.

\bibitem[CM18]{corsiMontalto2018Kirchhoff}
L.~Corsi and R.~Montalto,
\newblock Quasi-periodic solutions for the forced Kirchhoff equation on $\mathbb T^d$,
\newblock \emph{Nonlinearity} \textbf{31(11)} (2018), 5075--5109.

\bibitem[CMP21]{corsiMontaltoProcesi2021Airy}
L.~Corsi, R.~Montalto and M.~Procesi,
\newblock Almost-periodic response solutions for a forced quasi-linear Airy equation,
\newblock \emph{J. Dynam. Differential Equations} \textbf{33(3)} (2021), 1231--1267.

\bibitem[EGK16]{EGK}
L.~H.~Eliasson, B.~Gr\'ebert and S.~B.~Kuksin,
\newblock KAM for the nonlinear beam equation,
\newblock \emph{Geom. Funct. Anal.} \textbf{26(6)} (2016), 1588-1715.

\bibitem[EK10]{eliassonKuksin2010}
L.~H.~Eliasson and S.~B.~Kuksin,
\newblock KAM for the nonlinear Schr\"odinger equation,
\newblock \emph{Ann. of Math. (2)} \textbf{172} (2010), 371--435.

\bibitem[FS14]{faouSchratz2014}
E.~Faou and K.~Schratz,
\newblock Asymptotic preserving schemes for the Klein--Gordon equation in the non-relativistic limit regime,
\newblock \emph{Numer. Math.} \textbf{126(3)} (2014), 441--469.

\bibitem[FGP20]{feolaGiulianiProcesi2020DP}
R.~Feola, F.~Giuliani and M.~Procesi,
\newblock Reducible KAM tori for the Degasperis--Procesi equation,
\newblock \emph{Comm. Math. Phys.} \textbf{377(3)} (2020), 1681--1759.

\bibitem[FG24]{feolaGiuliani2024WW}
R.~Feola and F.~Giuliani,
\newblock Quasi-periodic traveling waves on an infinitely deep perfect fluid under gravity,
\newblock \emph{Mem. Amer. Math. Soc.} \textbf{295} (2024), No. 1471.

\bibitem[FMT25]{feolaMontaltoTerracina2025WW}
R.~Feola, R.~Montalto and S.~Terracina,
\newblock Time quasi-periodic three-dimensional traveling gravity water waves,
\newblock preprint: arXiv:2509.10318v1 (2025); \emph{to appear in Invent. Math.}

\bibitem[FP15]{feolaProcesi2015}
R.~Feola and M.~Procesi,
\newblock Quasi-periodic solutions for fully nonlinear forced reversible Schr\"odinger equations,
\newblock \emph{J. Differential Equations} \textbf{259} (2015), 3389--3447.

\bibitem[FP17]{feolaProcesi2017autonomousNLS}
R.~Feola and M.~Procesi,
\newblock KAM for quasi-linear autonomous NLS,
\newblock preprint, arXiv:1705.07287 (2017).


\bibitem[FM24]{franzoi2024kam}
L.~Franzoi and R.~Montalto,
\newblock A KAM approach to the inviscid limit for the 2D Navier--Stokes equations,
\newblock \emph{Ann. Henri Poincar\'e} \textbf{25(12)} (2024), 5231--5275.

\bibitem[Giu17]{giuliani2017generalizedKdV}
F.~Giuliani,
\newblock Quasi-periodic solutions for quasi-linear generalized KdV equations,
\newblock \emph{J. Differential Equations} \textbf{262(10)} (2017), 5052--5132.

\bibitem[GVR18]{GVdR}
B.~Gr\'ebert and V.~Vila\c ca da Rocha,
\newblock Stable and unstable time quasi-periodic solutions for a system of coupled NLS equations,
\newblock \emph{Nonlinearity} \textbf{31(10)} (2018), 4776-4811.

\bibitem[JY26]{JiaYang2026}
R.~Jia and Z.~Yang,
\newblock Non-relativistic limit from complex Klein-Gordon equation to a system of coupled nonlinear Schr\"odinger equations in 2D and 3D,
\newblock \emph{Chinese Ann. Math. Ser. B} \textbf{47} (2026), 753-770.

\bibitem[Kuk87]{Kuk87}
S.~B. Kuksin.
\newblock Hamiltonian perturbations of infinite-dimensional linear systems with an imaginary spectrum.
\newblock {\em  Functional Analysis and Its Applications}, \textbf{21} (1987), 192--205.

\bibitem[KP96]{kuksin1996invariant}
S.~Kuksin and J.~P\"oschel,
\newblock Invariant Cantor manifolds of quasi-periodic oscillations for a nonlinear Schr\"odinger equation,
\newblock \emph{Ann. of Math. (2)} \textbf{143} (1996), 149--179.

\bibitem[LW23]{leiWu2023}
Z.~Lei and Y.~Wu,
\newblock Non-relativistic limit for the cubic nonlinear Klein--Gordon equations,
\newblock preprint, arXiv:2309.10235 (2023).

\bibitem[LW10]{linWu2010}
C.-K.~Lin and K.-C.~Wu,
\newblock Singular limits of the Klein--Gordon equation,
\newblock \emph{Arch. Ration. Mech. Anal.} \textbf{197} (2010), 689--711.

\bibitem[LZ26]{luZheng2026}
Y.~Lu and Z.~Zheng,
\newblock Uniform-in-time error estimates for the Klein--Gordon--Schr\"odinger system in the nonrelativistic limit regime,
\newblock preprint, arXiv:2611.20646 (2026).

\bibitem[Mac01]{machihara2001nonrelativistic}
S.~Machihara,
\newblock The nonrelativistic limit of the nonlinear Klein-Gordon equation,
\newblock \emph{Funkcial. Ekvac.} \textbf{44} (2001), 243-252.

\bibitem[MNO02]{machihara2002nonrelativistic}
S.~Machihara, K.~Nakanishi and T.~Ozawa,
\newblock Nonrelativistic limit in the energy space for nonlinear
Klein-Gordon equations,
\newblock \emph{Math. Ann.} \textbf{322} (2002), 603-621.

\bibitem[MN02]{masmoudi2002nonlinear}
N.~Masmoudi and K.~Nakanishi,
\newblock From nonlinear Klein-Gordon equation to a system of coupled
nonlinear Schr\"odinger equations,
\newblock \emph{Math. Ann.} \textbf{324} (2002), 359-389.

\bibitem[Mon17]{montalto2017Kirchhoff}
R.~Montalto,
\newblock Quasi-periodic solutions of forced Kirchhoff equation,
\newblock \emph{NoDEA Nonlinear Differential Equations Appl.} \textbf{24} (2017), article 9.

\bibitem[Naj90]{najman1990nonrelativistic}
B.~Najman,
\newblock The nonrelativistic limit of the nonlinear Klein-Gordon equation,
\newblock \emph{Nonlinear Anal.} \textbf{15} (1990), 217-228.

\bibitem[Nak02]{nakanishi2002scattering}
K.~Nakanishi,
\newblock Nonrelativistic limit of scattering theory
for nonlinear Klein--Gordon equations,
\newblock \emph{J. Differential Equations} \textbf{180}
(2002), 453--470.

\bibitem[Pas18]{pasquali2018almost}
S.~Pasquali,
\newblock Almost global existence for the nonlinear Klein-Gordon equation
in the nonrelativistic limit,
\newblock \emph{J. Math. Phys.} \textbf{59} (2018), 011502.

\bibitem[Pas19]{pasquali2019dynamics}
S.~Pasquali,
\newblock Dynamics of the nonlinear Klein-Gordon equation in the
nonrelativistic limit,
\newblock \emph{Ann. Mat. Pura Appl.} \textbf{198} (2019), 903-972.

\bibitem[P\"os96]{poschel1996kam}
J.~P\"oschel,
\newblock A KAM-theorem for some nonlinear partial differential equations,
\newblock \emph{Ann. Scuola Norm. Sup. Pisa Cl. Sci.} \textbf{23} (1996),
119-148.

\bibitem[P\"os96b]{poschel1996quasi}
J.~P\"oschel,
\newblock Quasi-periodic solutions for a nonlinear wave equation,
\newblock \emph{Comment. Math. Helv.} \textbf{71} (1996), 269--296.

\bibitem[PP15]{procesiProcesi2015}
C.~Procesi and M.~Procesi,
\newblock A KAM algorithm for the resonant non-linear Schr\"odinger equation,
\newblock \emph{Adv. Math.} \textbf{272} (2015), 399--470.

\bibitem[SWZ26]{shenWangZheng2026}
J.~Shen, Y.~Wang and H.~Zheng,
\newblock High-order asymptotic expansion for the nonlinear Klein--Gordon equation in the non-relativistic limit regime,
\newblock \emph{J. Differential Equations} \textbf{453}
(2026), 113832.

\bibitem[Tsu84]{tsutsumi1984nonrelativistic}
M.~Tsutsumi,
\newblock Nonrelativistic approximation of nonlinear Klein-Gordon equations
in two space dimensions,
\newblock \emph{Nonlinear Anal.} \textbf{8} (1984), 637-643.

\bibitem[Way90]{way90}
C.~E. Wayne.
\newblock Periodic and quasi-periodic solutions of nonlinear wave equations via
  {KAM} theory.
\newblock {\em Comm. Math. Phys.}, 127(3):479--528, 1990.

\end{thebibliography}
\end{document}